%% file: main.tex
\documentclass[a4paper,review]{cas-sc}

\usepackage[section]{placeins} 
\input{preamble}

\begin{document}

\title[mode=title]{A positivity-preserving and conservative high-order flux reconstruction method for the polyatomic Boltzmann--BGK equation}
\shorttitle{High-order flux reconstruction method for the polyatomic Boltzmann--BGK equation}
\shortauthors{T. Dzanic \textit{et al.}}

\author[1]{T. Dzanic}[orcid=0000-0003-3791-1134]
\cormark[1]
\cortext[cor1]{Corresponding author}
\ead{tdzanic@tamu.edu}
\author[1]{F. D. Witherden}[orcid=0000-0003-2343-412X]
\author[2]{L. Martinelli}[orcid=0000-0001-6319-7646]

\address[1]{Department of Ocean Engineering, Texas A\&M University, College Station, TX 77843, USA}
\address[2]{Department of Mechanical and Aerospace Engineering, Princeton University, Princeton, NJ 08544, USA}

\begin{abstract}
In this work, we present a positivity-preserving high-order flux reconstruction method for the polyatomic Boltzmann--BGK equation augmented with a discrete velocity model that ensures the scheme is discretely conservative. Through modeling the internal degrees of freedom, the approach is further extended to polyatomic molecules and can encompass arbitrary constitutive laws. The approach is validated on a series of large-scale complex numerical experiments, ranging from shock-dominated flows computed on unstructured grids to direct numerical simulation of three-dimensional compressible turbulent flows, the latter of which is the first instance of such a flow computed by directly solving the Boltzmann equation. The results show the ability of the scheme to directly resolve shock structures without any \textit{ad hoc} numerical shock capturing method and correctly approximate turbulent flow phenomena in a consistent manner with the hydrodynamic equations. 
\end{abstract}




\maketitle

\input{introduction}
\input{preliminaries}
\input{methodology}
\input{implementation}

\input{results}

\input{conclusions}

\section*{Acknowledgements}
\label{sec:ack}
This work was supported in part by the U.S. Air Force Office of Scientific Research via grant FA9550-21-1-0190 ("Enabling next-generation heterogeneous computing for massively parallel high-order compressible CFD") of the Defense University Research Instrumentation Program (DURIP) under the direction of Dr. Fariba Fahroo.

\bibliographystyle{unsrtnat}
\bibliography{reference}



\end{document}

%% file: preamble.tex
\usepackage{amssymb}
\usepackage{amsthm}
\usepackage{amsmath}
\usepackage{upgreek}
\usepackage{pdflscape}
\usepackage{listings}
\usepackage{multirow}
\usepackage[percent]{overpic}
\usepackage{multicol}
\usepackage{color}
\usepackage{stmaryrd}
\usepackage{xfrac}
\usepackage[capitalise]{cleveref}
\usepackage{booktabs}
\usepackage[section]{placeins} 
\usepackage[section]{algorithm}
\usepackage{calc}
\usepackage[titletoc]{appendix}
\usepackage{siunitx}
\usepackage{mathtools}
\usepackage{tabularx}
\usepackage{amsmath}
\usepackage{diagbox}
\usepackage[sort&compress,square,numbers]{natbib}

\usepackage{graphicx}
\usepackage{graphics}
\usepackage{wrapfig}
\usepackage{float}
\usepackage{subfig}
\usepackage[percent]{overpic}
\usepackage{varwidth}
\usepackage{tikz}
\usepackage{pgfplots}
\usepackage{adjustbox}
\usetikzlibrary{arrows,matrix,positioning,fit}
\usetikzlibrary{shapes,positioning}
\usetikzlibrary{backgrounds}
\usepackage{tikz-layers}
\usepgfplotslibrary{fillbetween}
\usetikzlibrary{intersections}
\pgfplotsset{compat=1.14}
\usepgfplotslibrary{colorbrewer}
\usepgfplotslibrary{patchplots}
\usepgfplotslibrary[colorbrewer]
\usetikzlibrary{pgfplots.colorbrewer}
\usetikzlibrary[pgfplots.colorbrewer]
\usepgfplotslibrary{units}
\usetikzlibrary{spy}
\usepackage{pgfplotstable}
\usepackage{arrayjobx}
\graphicspath{ {./figs/} }
\usetikzlibrary{external}

\newlength\myheight
\newlength\mydepth
\settototalheight\myheight{Xygp}
\newcommand*\inlinegraphics[1]{%
  \settototalheight\myheight{Xygp}%
  \settodepth\mydepth{Xygp}%
  \raisebox{-\mydepth}{\includegraphics[height=\myheight]{#1}}%
}
\newcommand\orcid[1]{\href{https://orcid.org/#1}{\inlinegraphics{orcid_16x16.png}}}

\makeatletter
\def\BState{\State\hskip-\ALG@thistlm}
\makeatother

\newdefinition{definition}{Definition}[section]

\newcommand{\shalf}{{\sfrac{1}{2}}}

%% file: introduction.tex
\section{Introduction}
\label{sec:intro}
The development and application of computational fluid dynamics methods have, for the most part, focused on solutions of the Euler and Navier--Stokes equations to give insight into the dynamics of fluid flow. These approaches naturally rely on the notion that the continuum assumption can be regarded as valid for the given fluid, such that transport phenomena can be adequately described solely through the macroscopic flow variables (e.g., density, velocity, etc). However, for many problems of interest, such as rarefied gases and hypersonic flows, the continuum assumption starts to break down as collisions between molecules become less frequent and regions of the flow significantly deviate from thermodynamic equilibrium. In such cases, it becomes necessary to revert to the governing equations of molecular gas dynamics which underpin the macroscopic behavior of the fluid. These methods, which derive from the kinetic theory of gases, can offer a more detailed description of non-equilibrium systems and flows outside of the continuum regime while seamlessly recovering the hydrodynamic equations in the asymptotic limit.

To characterize gas dynamics at the molecular level, the Boltzmann equation can be used to provide a statistical description of particle transport and collision \citep{Cercignani1988}. In comparison to the Euler and Navier--Stokes equations, which solve for the nonlinear transport of the conserved flow variables, the Boltzmann equation simply solves for a scalar quantity, the \textit{distribution function} or \textit{probability density function}, for which the transport can be represented linearly. The complexity of this governing equation instead comes from two main sources: dimensionality and suitable modeling of the interactions among molecules. With regard to the former, the distribution function can be considered as a probability measure for particles existing at a given location with a certain velocity. As such, it must be defined over physical space and velocity space (and time), which can require discretization of up to 6-dimensional spaces. This problem is further compounded when attempting to model additional internal degrees of freedom  (e.g., rotation, vibration, etc.), which can require discretizations of even higher dimensionality \citep{Baranger2020}. Furthermore, the proper discretization of the velocity space, which is, in theory, unbounded, is an open problem, with varying methods offering differing advantages in terms of accuracy, stability, and computational cost \citep{Heyningen2021}. For the latter, the approximation of particle interactions in velocity space can be the most costly and algorithmically complex part of the Boltzmann equation, requiring integration over spaces of even higher dimension than the velocity space itself. This collision approximation also can, in many scenarios, cause the governing equation to become exceedingly stiff, further complicating attempts at efficient computation. 

To reduce the computational cost associated with the direct computation of the collision term, an approximate model of the collision process can instead be used. One such approach, the BGK (Bhatnagar–Gross–Krook) \citep{Bhatnagar1954} model, simplifies this process to a relaxation of a gas in non-equilibrium tending towards thermodynamic equilibrium, avoiding much of the computational cost associated with computing the collision integrals. This simplified model retains much of the behavior of the full collision operator, including its entropy-satisfying properties embedded in Boltzmann's H-theorem and convergence to the continuum approximation in the asymptotic limit, albeit with some corrections necessary to recover the proper transport coefficients. One additional drawback of the BGK model is that for general implementations, conservation of the macroscopic flow variables is only guaranteed if the integration in velocity space can be carried out exactly, which is usually not the case for discrete representations \citep{Mieussens2000JCP}. A common approach to remedy this issue is to resolve the velocity space finely enough such that the integration errors are negligible \citep{Evans2011}, but this comes at the expense of even further increasing the computational cost of solving the system. 

A potential improvement in the cost of solving the Boltzmann equation can be obtained through the use of a higher-fidelity spatial discretization. The typical schemes of choice for implementations of the Boltzmann equation on general unstructured grids are low-order finite volume methods \citep{Mieussens2000, Bernard2014} which are generally simple to implement and numerically robust. However, these benefits come with the drawback of low accuracy per degree of freedom and excessive numerical dissipation that makes them prohibitively expensive for scale-resolving simulations. An alternative class of schemes that are more suited for these types of problems are high-order discontinuous spectral element methods. These schemes offer many advantages in the simulation of multi-scale phenomena due to their ability to achieve high-order accuracy on general unstructured grids while still retaining a compact data structure suited for modern massively-parallel computing architectures. As a result, the amount of resolution required by these schemes to achieve a certain level of accuracy can be significantly lower than for low-order finite volume methods, which opens up the possibility of simulating complex flow problems that may otherwise be computationally intractable. The use of high-order discontinuous spectral element methods for the approximation of the Boltzmann equation has been attempted in the recent works of \citet{Xiao2021} and \citet{Jaiswal2022} with very promising results. The drawback of these methods is that they are more algorithmically complex and typically far less robust, which can limit their potential as a method of reducing the computational cost associated with the Boltzmann equation. 

The objective of this work is therefore to propose, develop, and validate an efficient  numerical approach for accurately and robustly solving the Boltzmann--BGK equation with the goal of drastically reducing the computational cost and enabling its application to three-dimensional problems that were hitherto intractable. We utilize the flux reconstruction scheme of \citet{Huynh2007}, a high-order discontinuous spectral element method, augmented with a positivity-preserving limiter that ensures the positivity of the distribution function as well as the macroscopic density/pressure. This spatial discretization is paired with a nodal discretization of the velocity space utilizing a modified form of the discrete velocity model of \citet{Mieussens2000} which guarantees the scheme remains conservative and well-balanced regardless of the resolution. Furthermore, we extend the approach to polyatomic molecules and general constitutive laws through an appropriate model for the internal degrees of freedom. The applicability of the proposed approach is shown in a series of complex numerical experiments, ranging from shock-dominated flows computed on unstructured grids without any numerical shock capturing algorithm to direct numerical simulation of three-dimensional compressible turbulent flows, the latter of which is, to the authors' knowledge, the first instance of such a flow computed by directly solving the Boltzmann equation. 

To this end, the organization of this paper is as follows. Some preliminaries on the various properties of the Boltzmann--BGK equation and the associated discretization methods are presented in \cref{sec:preliminaries}. The methodology of the proposed approach is then presented in \cref{sec:methodology}, followed by the details of the numerical implementation in \cref{sec:implementation}. The results of various numerical experiments are then shown in \cref{sec:results}, and conclusions are drawn in \cref{sec:conclusion}. 

%% file: preliminaries.tex
\section{Preliminaries}\label{sec:preliminaries}
This work pertains to approximations of the \textit{polyatomic} Boltzmann equation with arbitrary constitutive laws, given as
\begin{equation}\label{eq:boltzmann}
    \partial_t f (\mathbf{x}, \mathbf{u}, \zeta, t) + \mathbf{u} {\cdot} \nabla f = \mathcal C(f, f'),
\end{equation}
where $\mathbf{x} \in \Omega^{\mathbf{x}}$ is the physical space for some physical domain $\Omega^{\mathbf{x}} \subseteq \mathbb R^d$ and spatial dimension $d$, $\mathbf{u} \in \Omega^{\mathbf{u}}$ is the associated velocity space for some velocity domain $\Omega^{\mathbf{u}} \in \mathbb R^m$ and velocity dimension $m \geq d$, $\zeta \in \mathbb R^+$ is the internal energy, $f (\mathbf{x}, \mathbf{u}, \zeta, t)\in \mathbb R$ is a scalar particle distribution function, and $\mathcal C(f, f')$ is some collision operator that accounts for intermolecular interactions \citep{Cercignani1988}. This distribution function, describing the microscopic state of the system, gives a measure of the probability density of particles existing at a given location $\mathbf{x}$ traveling at a given velocity $\mathbf{u}$ (the combination of which is referred to as a point in \emph{phase space}) with an internal energy $\zeta$. The internal energy is given as a continuous variable that takes into account any additional degrees of freedom such as rotational and vibrational modes or unresolved velocity components \citep{Baranger2020}. From this distribution, one can gather information about the unique macroscopic state of the system (i.e., the conserved flow variables $\mathbf{Q}(\mathbf{x}, t)$) through its moments,
\begin{equation}\label{eq:moments}
    \mathbf{Q}(\mathbf{x}, t) = \left[\rho, \rho \mathbf{U}, E \right]^T = 
    \int_{\mathbb R^d} \int_{0}^{\infty} f (\mathbf{x}, \mathbf{u}, \zeta, t)\ \boldsymbol{\psi} (\mathbf{u}, \zeta) \ \mathrm{d}\zeta\ \mathrm{d}\mathbf{u} 
\end{equation}
where $\rho$ is the density, $\rho \mathbf{U}$ is the momentum vector, $E$ is the total energy, and $\boldsymbol{\psi} (\mathbf{u}, \zeta) \coloneqq [1, \mathbf{u}, (\mathbf{u}\cdot\mathbf{u})/2 + \zeta]^T$ is the vector of collision invariants. We utilize the notation that the lowercase symbol $\mathbf{u}$ refers to the microscopic velocity whereas the uppercase symbol $\mathbf{U}$ refers to the macroscopic velocity.

\subsection{Bhatnagar–Gross–Krook Operator}
Due to its high-dimensional nature, the computational cost of solving the scalar Boltzmann equation is considerably higher than the associated conservation laws for the conserved variables $\mathbf{Q}(\mathbf{x}, t)$. This is further exacerbated by the complexity of the collision operator $C(f, f')$, which in its full form can require integration over spaces of even higher dimension. To ameliorate this cost, an approximate form of the collision operator was introduced by \citet{Bhatnagar1954} by modeling the collision of particles as a multidimensional relaxation process of a gas tending towards thermodynamic equilibrium. The collision operator for the BGK (Bhatnagar–Gross–Krook) model can be given as 
\begin{equation}
    \mathcal C (f, f') \approx \frac{g(\mathbf{x}, \mathbf{u}, \zeta, t) - f(\mathbf{x}, \mathbf{u}, \zeta, t)}{\tau},
\end{equation}
where $g(\mathbf{x}, \mathbf{u}, \zeta, t)$ is the equilibrium distribution function and $\tau$ is the collision time scale. It can be shown that if the following compatibility condition is satisfied, 
\begin{equation}\label{eq:compatibility}
    \int_{\mathbb R^d} \int_{0}^{\infty} g (\mathbf{x}, \mathbf{u}, \zeta, t)\ \boldsymbol{\psi} (\mathbf{u}, \zeta) \ \mathrm{d}\zeta\ \mathrm{d}\mathbf{u}  = \int_{\mathbb R^d} \int_{0}^{\infty} f (\mathbf{x}, \mathbf{u}, \zeta, t)\ \boldsymbol{\psi} (\mathbf{u}, \zeta) \ \mathrm{d}\zeta\ \mathrm{d}\mathbf{u},
\end{equation}
the resulting system is conservative with respect to the conserved flow variables. 

For the polyatomic case, the equilibrium distribution function is given as
\begin{equation}\label{eq:equilibriummoments}
    g(\mathbf{x}, \mathbf{u}, \zeta, t) = g_{\mathbf{u}}(\mathbf{x}, \mathbf{u}, t) \times g_\zeta (\mathbf{x}, \zeta, t),
\end{equation}
where $g_{\mathbf{u}}$ is the equilibrium distribution function of the monatomic case and $g_\zeta (\mathbf{x}, \zeta, t)$ is an additional term used to augment the equilibrium distribution function to account for internal degrees of freedom. We will assume that $g_\zeta (\mathbf{x}, \zeta, t)$ is properly normalized such that its inclusion does not affect the moments, i.e., 
\begin{equation}\label{eq:normdist}
    \int_{\mathbb R^d} \int_{0}^{\infty} g (\mathbf{x}, \mathbf{u}, \zeta, t)\ \boldsymbol{\psi} (\mathbf{u}, \zeta) \ \mathrm{d}\zeta\ \mathrm{d}\mathbf{u}  = \int_{\mathbb R^d}  g_\mathbf{u} (\mathbf{x}, \mathbf{u}, t)\ \boldsymbol{\psi} (\mathbf{u}) \ \mathrm{d}\mathbf{u} = \mathbf{Q}(\mathbf{x}, t).
\end{equation}
From the H-theorem, the monatomic equilibrium state is generally taken to be the state that minimizes the entropy $H(z)$, i.e.,
\begin{equation}
    g_{\mathbf{u}} = \underset{z}{\mathrm{arg\ min}}\ H(z),
\end{equation}
where
\begin{equation}
    H(z) = \int_{\mathbb R^d} z \log(z) \ \mathrm{d} \mathbf{u}.
\end{equation}
In the \emph{continuous} case, the distribution function in velocity space that minimizes this entropy and satisfies \cref{eq:moments,eq:equilibriummoments} is a Maxwellian of the form
\begin{equation}
    g_{\mathbf{u}}(\mathbf{x}, \mathbf{u}, t) = \frac{\rho (\mathbf{x}, t)}{\left[2 \pi \theta(\mathbf{x}, t) \right]^{d/2}}\exp \left [-\frac{ \|\mathbf{u} - \mathbf{U}(\mathbf x, t) \|_2^2}{2 \theta (\mathbf x, t)} \right],
\end{equation}
where $\theta = P/\rho$ is a scaled temperature, $P = (\gamma - 1)(E - \rho(\mathbf{U}\cdot\mathbf{U})/2)$ is the pressure, and $\gamma$ is the ratio of specific heats. We utilize the notation $\|\cdot \|_2^2$ to denote the squared norm along the spatial dimension and $g(\mathbf{Q}(\mathbf{x}, t))$ to denote the Maxwellian corresponding to the conserved variables $\mathbf{Q}(\mathbf{x}, t)$.

\subsection{Internal Degrees of Freedom}
In its standard form, the Boltzmann--BGK approach assumes a monatomic molecule in which the only degrees of freedom are its translational components. This limits the number of degrees of freedom, $n$, to the dimensionality of the velocity space (i.e., $n = m$), which in turn, constrains the specific heat ratio by the relation $\gamma = 1 + 2/n$. This constraint on the degrees of freedom imposes several difficulties for the application of the approach to practical problems. First, there is significant interest in the simulation of real gases, many of which are polyatomic in nature, for which the standard approach cannot be used. Secondly, even for monatomic molecules, certain flow conditions such as high temperatures can induce additional kinetic effects such as vibration which cannot be accounted for without additional degrees of freedom. Lastly, although this drawback is less critical, the simulation of monatomic flows with spatial homogeneity (e.g., one-dimensional flows) requires the discretization of higher-dimensional velocity spaces to recover the proper specific heat ratio, unnecessarily increasing the computational cost. 

To simultaneously account for the additional degrees of freedom of polyatomic gases, ancillary kinetic modes such as rotation and vibration, and arbitrary constitutive laws, \citet{Baranger2020} introduced a model for internal degrees of freedom for the molecule. In this approach, an additional dimension is considered to account for the internal energy of the molecule $\zeta$, and the equilibrium distribution function is modified by the factor 
\begin{equation}
    g_\zeta (\mathbf{x}, \zeta, t) = \Lambda(\delta) \left (\frac{\zeta}{\theta(\mathbf{x}, t)} \right)^{\frac{\delta}{2} - 1} \frac{1}{\theta(\mathbf{x}, t)}\exp \left(-\frac{\zeta}{\theta(\mathbf{x}, t)} \right),
\end{equation}
where $\delta \geq 0$ is the number of internal degrees of freedom and $\Lambda (\delta) = 1/\Gamma(\delta/2)$ is a normalization factor to ensure \cref{eq:normdist} is satisfied. With this model, the total degrees of freedom are then $n = m + \delta$, such that a specific heat ratio of $\gamma = 1 + 2/(m + \delta)$ can be recovered. For the monatomic case, the internal degrees of freedom are simply set as $\delta = 0$ and the $\zeta$ domain is neglected. Note that $\delta$ does not necessarily have to be an integer nor a constant (see \citet{Baranger2020}, Section 3), which allows for the implementation of arbitrary constitutive laws for multi-physics applications (e.g., temperature-dependent internal degrees of freedom). 

\subsection{Velocity and Internal Energy Spaces}
Discretizations for the velocity space generally fall into one of two approaches: nodal form or moment form. In the former approach, the infinite velocity space is truncated onto a finite domain and represented by a set of discrete points. The associated transport equation, \cref{eq:boltzmann}, is then solved for each one of these discrete points. In the latter approach, the distribution function and the velocity space are represented by a set of global basis functions, and the transport equation is solved for the entire velocity space using an appropriate treatment of the inner product in \cref{eq:boltzmann}. (e.g., orthogonal projection). Due to their orthogonality with respect to the measure $\exp(-x^2)$ on the real domain, Hermite polynomials are a natural choice for the basis functions, typically allowing for a more computationally efficient representation but at the expense of robustness, particularly for flows deviating from equilibrium \cite{Heyningen2021}. However, for the internal energy domain, finding a set of suitable global basis functions becomes less trivial. 

The focus of this work is on the former approach -- nodal discretizations of the velocity and internal energy spaces -- as this method is the most general, encompassing arbitrary distribution functions and enabling the straightforward implementation of additional physics such as the inclusion of body forces. For this approach, the discretization for the resulting system of equations becomes greatly simplified, reducing to a set of linear advection equations that are essentially independent across the spatial and velocity/internal energy domains and share a dependency only through the nonlinear source term. However, to accurately resolve the underlying physics, it must be ensured that the resolution is fine enough and the truncated velocity domain encompasses enough of the predominant behavior of the distribution function. Therefore, the proper discretization of the velocity/internal energy domains becomes a critical factor in the accuracy of the Boltzmann approach, and, as will be shown in \cref{ssec:dvm}, can also affect the numerical properties of the scheme such as conservation and well-balancing.

\subsection{Limiting Behavior}
Although the Boltzmann--BGK model is a significantly more general approach than the equations governing fluid flow in the continuum limit, it can be shown to recover these hydrodynamic equations in the asymptotic limit. The primary non-dimensional quantity characterizing the Boltzmann equation is the Knudsen number, defined as the ratio
\begin{equation}
    Kn = \frac{\lambda}{L_{\mathrm{ref}}},
\end{equation}
where $L_{\mathrm{ref}}$ is the characteristic length scale and $\lambda$ is the particle mean free path. In the limit as $Kn \to 0$, the continuum approximation can be regarded as valid, whereas for higher $Kn$, the flow enters the rarefied regime. The Knudsen number can be related to Mach number $M$ and Reynolds number $Re$ as 
\begin{equation}
    Kn = \sqrt{\frac{\gamma \pi}{2}} \frac{M}{Re},
\end{equation}
where $\gamma$ is the specific heat ratio. The Knudsen number can similarly be expressed in terms of the macroscopic flow variables as 
\begin{equation}
    Kn = \sqrt{\frac{\gamma \pi}{2}} \frac{\mu}{\rho c_s L_{\mathrm{ref}}},
\end{equation}
where $\mu$ is the dynamic viscosity and $c_s = \sqrt{\gamma P/\rho}$ is the speed of sound.

In the low Knudsen number limit, the BGK approximation converges to the hydrodynamic limit governed by the Euler/Navier--Stokes equations with the dynamic viscosity given as 
\begin{equation}
    \mu = \tau P,
\end{equation}
which defines the collision time in terms of the Knudsen number and the macroscopic flow quantities as
\begin{equation}
    \tau = Kn \sqrt{\frac{2}{\gamma \pi}} \frac{\rho c_s L_{\mathrm{ref}}}{P} = \sqrt{\frac{2\gamma}{\pi}} \frac{Kn L_{\mathrm{ref}}}{c_s}.
\end{equation}
In the low Mach number limit, the variations in pressure become small, such that for a fixed collision time, the dynamic viscosity remains essentially constant. However, for higher Mach numbers where the pressure, and therefore the viscosity, is expected to vary more strongly, the collision time can be set adaptively based on the macroscopic flow state to recover a more physically consistent viscosity coefficient that varies as a function of temperature instead of pressure \citep{Mieussens2000}. This method is briefly explored in this work. Furthermore, due to the use of a single relaxation time, momentum relaxation and thermal relaxation are assumed to be identical, such that the model recovers a unit Prandtl number ($Pr = 1$), similarly to the assumption of the Reynolds analogy. Relatively simple modifications to the equilibrium distribution function can be performed to recover other values for the Prandtl number \citep{Shakhov1972, Wang2019}, but they are omitted in this work. 

%% file: methodology.tex
\section{Methodology}\label{sec:methodology}

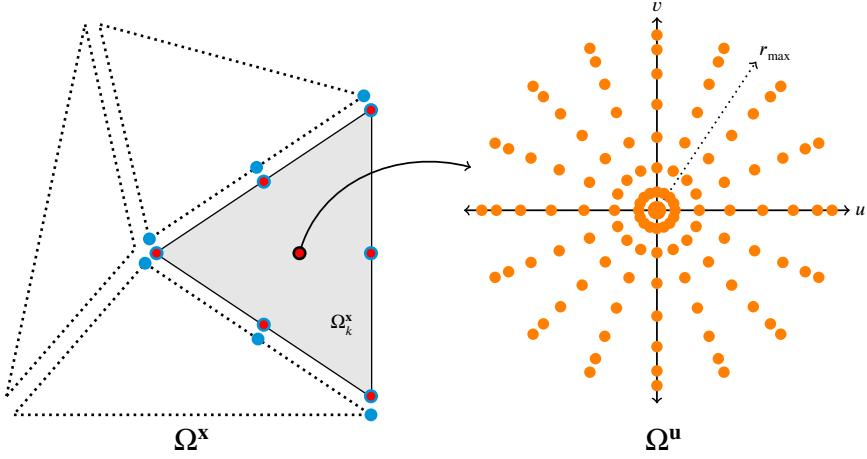
\begin{figure}[tbhp]
        \adjustbox{width=0.7\linewidth, valign=b}{\input{figs/scheme}}
        \newline
        \caption{\label{fig:scheme} Schematic of a two-dimensional phase space discretization using an unstructured spatial domain $\Omega^{\mathbf{x}}$ with $\mathbb P_2$ elements and a velocity domain $\Omega^{\mathbf{v}}$ with $N_r = 8$, $N_\phi = 16$ and no internal degrees of freedom. Circles denote the spatial solution nodes (red), interface flux nodes (blue), and velocity space nodes (orange), respectively.}
    \end{figure}

A description of the proposed numerical approach for discretizing the Boltzmann--BGK equation is presented in this section, including details on the positivity-preserving high-order spatial discretization, the discretely conservative velocity and internal energy discretizations, boundary conditions for the system, and a guideline for the required resolution levels. A schematic of this discretization is presented in \cref{fig:scheme}.

\subsection{Spatial Discretization}
For a fixed location $\mathbf{u}_0$ in velocity space, the evolution of the particle distribution function is equivalent to a linear advection equation with a nonlinear source term $S$. With a slight abuse of notation, we let $f(\mathbf{x},\zeta, t) = f(\mathbf{x}, \mathbf{u}_0, \zeta, t)$ for some arbitrary $\mathbf{u}_0 \in \Omega^{\mathbf{u}}$, which yields the transport equation
\begin{equation}\label{eq:transport}
    \partial_t f (\mathbf{x}, \zeta, t) +  \boldsymbol{\nabla} {\cdot} \mathbf{F}(f) = S, \quad \quad \mathbf{F}(f) = \mathbf{u}_0 f.
\end{equation}
Here, $\mathbf{u}_0$ represents some constant advection velocity and $\mathbf{F}(f)$ is the flux. The left-hand side of this transport equation is discretized using the flux reconstruction (FR) method of \citet{Huynh2007}, a generalization of the nodal discontinuous Galerkin method \citep{Hesthaven2008DG}. 

In this approach, the spatial domain $\Omega^{\mathbf{x}}$ is partitioned into $N_e$ elements $\Omega^{\mathbf{x}}_k $ such that $\Omega^{\mathbf{x}} = \bigcup_{N_e}\Omega^{\mathbf{x}}_k$ and $\Omega^{\mathbf{x}}_i\cap\Omega^{\mathbf{x}}_j=\emptyset$ for $i\neq j$, shown on the left-hand side of \cref{fig:scheme}. Within each element $\Omega^{\mathbf{x}}_k$, the distribution function is approximated by a polynomial of the form
\begin{equation}
    f (\mathbf{x}) = \sum_{i = 1}^{N_s} f (\mathbf{x}^s_i) {\phi}_i (\mathbf{x}),
\end{equation}
where $\mathbf{x}^s_i \in \Omega^{\mathbf{x}}_k \ \forall \ i \in \{1,..., N_s\}$ is a set of $N_s$ solution nodes and ${\phi}_i (\mathbf{x})$ is a set of nodal basis functions with the property ${\phi}_i (\mathbf{x}^s_j) = \delta_{ij}$. The order of the approximation, represented by $\mathbb P_p$ for some order $p$, is defined as the maximal order of $f (\mathbf{x})$. Due to the linearity of the flux, the discontinuous flux $\mathbf{F}^D(\mathbf{x})$ can be simply represented as 
\begin{equation}
    \mathbf{F}^D(\mathbf{x}) = \mathbf{u}_0 f (\mathbf{x}).
\end{equation}
To allow communication between elements and to ensure consistency between the polynomial spaces of the solution and the divergence of the flux, the corrected flux is formed by amending the discontinuous flux with additional correction terms.
\begin{equation}
    \mathbf{F}^C(\mathbf{x}) = \mathbf{F}^D(\mathbf{x}) + \sum_{i = 1}^{N_f} \left[F^I_i - \mathbf{F}^D(\mathbf{x}^f_i)\cdot \mathbf{n}_i \right] \mathbf{g}_i (\mathbf{x}).
\end{equation}
Here, $\mathbf{x}^f_i \in \partial \Omega^{\mathbf{x}}_k \ \forall \ i \in \{1,..., N_f\}$ is a set of $N_f$ interface flux nodes, $\mathbf{n}_i$ is their associated outward-facing normal vector, $F^I_i$ is a common interface flux yet to be defined, and $\mathbf{g}_i$ is the correction function associated with the given flux node. These correction functions have the properties that
\begin{equation}
    \mathbf{n}_i \cdot\mathbf{g}_j (\mathbf{x}^f_i) = \delta_{ij} \quad \mathrm{and} \quad \sum_{i = 1}^{N_f} \mathbf{g}_i (\mathbf{x}) \in RT_p,
\end{equation}
where $RT_p$ is the Raviart--Thomas space of order $p$. In this work, these correction functions are chosen such as to recover the nodal discontinuous Galerkin approach \citep{Huynh2007, Hesthaven2008DG, Trojak2021}.

At each interface flux point $\mathbf{x}^f_i$, there exists a pair of solution values $(f^-_i, f^+_i)$, denoting the inner value of the solution (from the element of interest) and the outer value of the solution (from the face-adjacent element), respectively. To ensure the correction direction of information propagation, the common interface flux $F^I_i$ for systems of equations is generally computed using a Riemann solver \citep{Rusanov1962, Roe1981} based on this solution pair and their associated normal vector. However, due to the point-wise equivalency of the Boltzmann equation to the linear advection equation, the common interface flux can instead be computed using a simple upwinding approach as
\begin{equation}
    F_i^I = \begin{cases}
    u_n f_i^-, \quad \quad \mathrm{if} \ u_n > 0\\
    u_n f_i^+, \quad \quad \mathrm{else},
    \end{cases}
\end{equation}
where $u_n = \mathbf{u}_0 \cdot \mathbf{n}_i$. With this discretization and an approximation of the source term $S$, the system can be evolved using a suitable temporal integration approach. 

\subsection{Velocity Discretization}\label{ssec:velocity}
For the velocity discretization, the velocity space was truncated onto a finite domain $\Omega^{\mathbf{u}} \subset \mathbb R^m$ and represented nodally by $N_v$ discrete points, shown on the right-hand side of \cref{fig:scheme}. As the distribution functions tend to have a radially-symmetric nature in many cases, particularly for flows near equilibrium, it is convenient for $m > 1$ to represent the space in polar/spherical coordinates as 
\begin{equation}    
    \begin{bmatrix}
           u \\
           v
    \end{bmatrix} 
    = 
    \begin{bmatrix}
           r \cos \phi + U_0 \\
           r \sin \phi + V_0
    \end{bmatrix} 
    \quad \mathrm{and} \quad 
    \begin{bmatrix}
           u \\
           v \\
           w
    \end{bmatrix} 
    = 
    \begin{bmatrix}
           r \sin \psi \cos \phi  + U_0 \\
           r \sin \psi \sin \phi   + V_0\\
           r \cos \psi  + W_0
    \end{bmatrix}, 
\end{equation}
where $r \in (0, r_{\max}]$ is the radial distance for some maximum radial extent $r_{\max}$, $\phi \in [0, 2\pi)$ is the polar angle, $\psi \in [0, \pi)$ is the azimuthal angle, and $\mathbf{U}_0 = [U_0, V_0, W_0]$ is a vector of velocity offsets. The offsets were computed as the component-wise average of the minima and maxima of the initial macroscopic velocities in the domain. 
\begin{equation}
    \mathbf{U}_0 = \frac{1}{2} \left (\underset{\mathbf{x}\in\Omega^{\mathbf{x}}}{\max}\left [ \mathbf{U} (\mathbf{x}, 0) \right] + \underset{\mathbf{x}\in\Omega^{\mathbf{x}}}{\min}\left [ \mathbf{U} (\mathbf{x}, 0) \right]  \right).
\end{equation}

The two primary parameters in the velocity space discretization are the radial extent of the domain and the distribution of the nodes. For the former, the choice of the radial extent plays a large role in the properties of the scheme -- too small of an extent can cause inaccurate predictions of the underlying physics due to the loss of information regardless of the resolution of the discretization, whereas too large of an extent can degrade the efficiency of the discretization while also restricting the maximum permissible explicit time step due to unnecessarily large particle velocities. As described in \citet{Evans2011}, a standard approach is to prescribe the radial extent as some factor of the thermal velocity (i.e., speed of sound). In this work, we take a similar approach but apply an \textit{a priori} case-dependent metric to automatically account for variation in the velocity distributions between problems. The radial extent is computed  as
\begin{equation}
    r_{\max} = k\underset{\mathbf{x}\in\Omega^{\mathbf{x}}}{\max} \left [c_s (\mathbf{x}, 0)\right] + \frac{1}{2}\|  \delta \mathbf{U} \|_2,
\end{equation}
where $c_s = \sqrt{\gamma P/\rho}$ is the speed of sound and $\| \delta \mathbf{U} \|_2$ is the component-wise maximum difference in the macroscopic velocities in the domain, i.e., 
\begin{equation}
    \delta \mathbf{U} = \underset{\mathbf{x}\in\Omega^{\mathbf{x}}}{\max}\left [ \mathbf{U} (\mathbf{x}, 0) \right] - \underset{\mathbf{x}\in\Omega^{\mathbf{x}}}{\min}\left [ \mathbf{U} (\mathbf{x}, 0) \right].
\end{equation}
The parameter $k$ is chosen such that the contribution of the distribution function outside of the domain is negligible. Therefore, $k$ is computed such that the relative magnitude of the distribution function outside of the velocity domain compared to the maximum is of $\mathcal O(\epsilon_{\mathbf{u}})$, where $\epsilon_{\mathbf{u}}$ is some small user-defined constant to be presented in \cref{sec:implementation}. If the distribution function within the domain is assumed to be initially Maxwellian, then $k$ can be analytically computed as 
\begin{equation}
    k = \sqrt{-\frac{2}{\gamma}\log(\epsilon_{\mathbf{u}})},
\end{equation}
which guarantees that the relative magnitude of the initial distribution function outside of velocity domain is \textit{at least} a factor of $\epsilon_{\mathbf{u}}$ smaller than the maxima everywhere in the domain. Note that this choice of domain extent does not take into account flows whose temperature/velocity are expected to drastically differ from their initial conditions. In those scenarios, it may be necessary to verify the sensitivity of the solution on the domain extent \textit{a posteriori}. 

The distribution of nodes was chosen with polar/spherical integration in mind without $\textit{a priori}$ knowledge of the expected distribution function. Since there should be no bias for any particular angle in the polar/azimuthal distribution, a constant spacing with equal weighting was used. For $N_\phi$ polar and $N_\psi$ azimuthal nodes, the nodal locations and weights were given as 
\begin{equation}
    \phi_i = i\Delta \phi \ \forall\ i \in \{0, ..., N_\phi-1 \}, \quad \quad \Delta \phi = \frac{2\pi}{N_\phi}, \quad \quad w_i^\phi = \Delta \phi,
\end{equation}
\begin{equation}
    \psi_i = (i + \shalf)\Delta \psi \ \forall\ i \in \{0, ..., N_\psi-1 \}, \quad \quad \Delta \psi = \frac{\pi}{N_\psi}, \quad \quad w_i^\psi = \Delta \psi.
\end{equation}
respectively. Note that one end is open for the polar distribution and both ends are open for the azimuthal distribution to account for periodicity in the respective components and to avoid coincident nodal points. For the radial distribution, the Gauss-Legendre nodes $x^g$ and weights $w^g$ were chosen as they provide satisfactory integration properties, avoid the singularity at $r = 0$, and cluster nodes at both the head and tail of the distribution. By applying a suitable normalization, the radial nodal locations and weights can be given as 
\begin{equation}
    r_i = r_{\max} \frac{x_i^g + 1}{2}, \quad \quad w_i^r = \frac{r_{\max}}{2}w_i^g  J(r_i),
\end{equation}
where $J(r_i)$ is the Jacobian of the transformation from polar/spherical coordinates to Cartesian coordinates. For a velocity domain dimension $m$, the Jacobian is calculated as the surface area of an $m$-ball of radius $r_i$, e.g.,
\begin{equation}
    J(r_i) = \begin{cases}
    2\pi r_i, \quad \quad \mathrm{\ if\ } m = 2, \\
    4 \pi r_i^2, \quad \quad \mathrm{\, if\ } m = 3.\\
    \end{cases}
\end{equation}
The dimensionality of the velocity space can then be calculated as $N_v = N_r N_\phi$ for $m=2$ and $N_v = N_r N_\phi N_\psi$ for $m=3$. 

In the case of $m=1$, the nodal distribution can be given across $[-r_{\max}, r_{\max}]$. To maintain consistency with the polar/spherical nodal distributions and to retain the clustering of the nodes around the head and tail, the nodal distribution was computed as the concatenation of two Gauss-Legendre nodal distributions of size $N_v/2$ on $[-r_{\max}, 0]$ and $[0, r_{\max}]$, respectively, which is identical to the two-dimensional case with $N_\phi = 2$. As such, for the one-dimensional cases in this work, $N_v$ is set to an even integer. 

\subsection{Internal Energy Discretization}
For the internal energy discretization with $\delta > 0$, a similar methodology as with the radial velocity discretization was applied. The internal energy nodal locations and weights are given in terms of the normalized Gauss-Legendre quadrature as
\begin{equation}
    \zeta_i = \zeta_{\max} \frac{x_i^g + 1}{2} \ \forall\ i \in \{1, ..., N_\zeta \}, \quad \quad w_i^r = \frac{\zeta_{\max}}{2} w_i^g.
\end{equation}
For the choice of $\zeta_{\max}$, similar arguments as with $r_{\max}$ can be made, but in this case, the upper bound does not have an impact on the maximum permissible explicit time step. To calculate $\zeta_{\max}$, a similar approach was taken as with $r_{\max}$ by solving for a value such that the relative contribution of the internal energy distribution outside of the domain is of $\mathcal O(\epsilon_{\zeta})$, where $\epsilon_{\zeta}$ is again some small user-defined constant to be presented in \cref{sec:implementation}. However, unlike for the velocity domain, there does not exist a closed-form solution to this problem for general values of $\delta$, and therefore it must be computed numerically. For several values of $\delta$ and $\epsilon_{\zeta}$, the ratio $\zeta_{\max}/\theta_{\max}$ was computed and presented in \cref{tab:zeta_max}.

\begin{figure}[tbhp]
    \centering
    \begin{tabularx}{\textwidth}{r | @{\extracolsep{\fill}} c c c c c c c}
        $\delta$ & $\epsilon_{\zeta} = 10^{-2}$ & $\epsilon_{\zeta} = 10^{-4}$ & $\epsilon_{\zeta} = 10^{-6}$ & $\epsilon_{\zeta} = 10^{-8}$ & $\epsilon_{\zeta} = 10^{-10}$ & $\epsilon_{\zeta} = 10^{-12}$ & $\epsilon_{\zeta} = 10^{-14}$ \\ \midrule
        $2$ & 4.605 & 9.210 & 13.816 & 18.421 & 23.026 & 27.631 & 32.236 \\
        $3$ & 5.453 & 10.380 & 15.175 & 19.916 & 24.628 & 29.320 & 33.999 \\
        $4$ & 6.471 & 11.667 & 16.627 & 21.488 & 26.295 & 31.067 & 35.815 \\
        $5$ & 7.656 & 13.065 & 18.165 & 23.133 & 28.026 & 32.870 & 37.680 \\
    \end{tabularx}
    \captionof{table}{\label{tab:zeta_max} Numerically computed values of the relative domain extent $\zeta_{\max}/\theta_{\max}$ for varying values of the tolerance $\epsilon_{\zeta}$ and internal degrees of freedom $\delta$.}
\end{figure}

\subsection{Discrete Velocity Model}\label{ssec:dvm}
For the Boltzmann--BGK equation, accurate integration over the velocity and internal energy spaces is of key importance to the properties of the numerical scheme such as conservation and well-balancing. However, issues can arise when these infinite, continuous spaces are represented discretely. If $\mathbf{M}$ is a discrete nodal integration operator over $\Omega^{\mathbf{u}} \times \Omega^\zeta$ with strictly non-negative entries and $\mathbf{x}$ represents some discrete nodal values, i.e.,
\begin{equation}
   \mathbf{M} \left ( \mathbf{x} \right) \approx \int_{\mathbb R^d} \int_{0}^{\infty} x\ \mathrm{d}\zeta\ \mathrm{d}\mathbf{u},
\end{equation}
then it can be seen that in the limit of infinite resolution and infinite domain size, the discrete integral will converge towards the continuous integral if the integration operator is consistent. However, for finite approximations, this equality does not hold in general, and therefore for some discrete equilibrium distribution function $\mathbf{g}$,
\begin{equation}
    \mathbf{M}  \left(\boldsymbol{\psi} \otimes \mathbf{g} (\mathbf{Q})\right) \neq \mathbf{Q}.
\end{equation}
Without this equality, the moments of the equilibrium distribution function do not yield the same conserved variables used to create the equilibrium distribution function, which results in a scheme that is not  conservative and is not well-balanced (i.e., a constant solution in thermodynamic equilibrium does not remain constant nor in thermodynamic equilibrium). 

For nodal discretizations of the velocity and internal energy spaces, a standard approach is to discretize these spaces with enough resolution such that the integration errors are negligible \citep{}. However, due to the high-dimensional nature of the governing equations, this combinatorial explosion results in computational costs that severely limit the applicability of these methods. As an alternative approach, \citet{Mieussens2000} introduced a discrete velocity model (DVM) for the Boltzmann--BGK equation to address these issues. In the DVM approach, a discrete equilibrium distribution function $\mathbf{g}$ is sought that satisfies the \textit{discrete} compatibility condition, 
\begin{equation}
\mathbf{M}  \left(\boldsymbol{\psi} \otimes \mathbf{g} \right) = \mathbf{M}  \left(\boldsymbol{\psi} \otimes \mathbf{f} \right) = \mathbf{Q},
\end{equation}
and minimizes the discrete entropy $H'$,
\begin{equation}
    \mathbf{g} = \underset{\mathbf{z}}{\mathrm{arg\ min}}\ H'(\mathbf{z}), \quad \quad H'(\mathbf{z}) = \mathbf{M}  \left(\mathbf{z} \log \mathbf{z}\right).
\end{equation}
It was shown in \citet{Mieussens2000} that there exists a unique solution $\mathbf{g}'$ to this problem, given by a \textit{modified} Maxwellian of the form
\begin{equation}
    \mathbf{g}' = \mathbf{g}(\mathbf{Q}'),
\end{equation}
where $\mathbf{g}(\mathbf{Q}')$ is a discrete Maxwellian corresponding to a modified set of conserved variables $\mathbf{Q}' \neq \mathbf{Q}$ that converge to $\mathbf{Q}$ in the limit of infinite resolution. 
However, in general, the values of $\mathbf{Q}'$ that satisfy the discrete compatibility condition do not have a closed form solution, and therefore they must be found iteratively. This can be cast as a root-finding problem, where $\mathbf{Q}'$ is the root of
\begin{equation}
    \mathbf{M}  \left(\boldsymbol{\psi} \otimes \mathbf{g} (\mathbf{Q}') \right) - \mathbf{Q} = 0.
\end{equation}
Note that due to the dimensional reduction caused by the discrete integration operator, this nonlinear optimization problem is only of the same dimension as $\mathbf{Q}$ (i.e., $d+2$), irrespective of the dimensionality of $\mathbf{g}$. Furthermore, the Jacobian $\mathbf{J}$ of the root-finding problem can be analytically calculated, allowing for the use of efficient root-finding algorithms such as Newton's method with minimal overhead. 

To simplify the presentation, we will use the same notation as \citet{Mieussens2000},
where the equilibirum distribution function is given with respect to a vector $\boldsymbol{\alpha} (\mathbf{Q}) = [\alpha_1, \alpha_2, ..., \alpha_{d+2}]^T$ instead of $\mathbf{Q}$ as 
\begin{equation}
    g(\boldsymbol{\alpha}) = \alpha_1 \exp \left [-\alpha_2 \left(\sum_{i=1}^{d} \left (\mathbf{u}_i -  \alpha_{i+2}  \right)^2 \right) \right],
\end{equation}
where for $d = 3$,
\begin{equation}
    \alpha_1 = \frac{\rho}{\left(2 \pi \theta\right)^{d/2}}, \quad
    \alpha_2 = \frac{1}{2 \theta}, \quad
    \alpha_3 = U, \quad
    \alpha_4 = V, \quad
    \alpha_5 = W.
\end{equation}
The Jacobian can then be given as 
\begin{equation}
    \mathbf{J}(\boldsymbol{\alpha}) = \mathbf{M} \left(\boldsymbol{\Theta} \otimes \boldsymbol{\psi} \otimes \mathbf{g} (\boldsymbol{\alpha}) \right),
\end{equation}
where
\begin{equation}
    \boldsymbol{\Theta} = \frac{\partial g(\boldsymbol{\alpha})}{\partial \boldsymbol{\alpha}} \frac{1}{g(\boldsymbol{\alpha})},
\end{equation}
such that for $d = 3$,
\begin{align}
    \boldsymbol{\Theta}_1 &= \frac{1}{\alpha_1}, \\
    \boldsymbol{\Theta}_2 &= -(\mathbf{u} - \alpha_3)^2 -(\mathbf{v} - \alpha_4)^2-(\mathbf{w} - \alpha_5)^2 + \frac{\delta - 4\zeta\alpha_2}{2\alpha_2},\\
    \boldsymbol{\Theta}_3 &= 2\alpha_2(\mathbf{u} - \alpha_3),\\
    \boldsymbol{\Theta}_4 &= 2\alpha_2(\mathbf{v} - \alpha_4), \\
    \boldsymbol{\Theta}_5 &= 2\alpha_2(\mathbf{w} - \alpha_5).
\end{align}
Note that due to the inclusion of the internal energy components, $\boldsymbol{\Theta}$ differs from the expression in \citet{Mieussens2000} by the quantity
\begin{equation}
\frac{\delta - 4\zeta\alpha_2}{2\alpha_2} = \frac{\partial g_\zeta(\boldsymbol{\alpha})}{\partial \boldsymbol{\alpha}} \frac{1}{g_\zeta(\boldsymbol{\alpha})}.
\end{equation}
For the case of $\delta = 0$, this term is simply neglected as $g_\zeta = 1$. The root can then be found using Newton's method through the relation
\begin{equation}
    \boldsymbol{\alpha}_{n+1} = \boldsymbol{\alpha}_{n} - \mathbf{J}(\boldsymbol{\alpha}_{n})^{-1} \left [\mathbf{M}  \left(\boldsymbol{\psi} \otimes \mathbf{g} (\boldsymbol{\alpha}_{n}) \right) - \mathbf{Q} \right],
\end{equation}
with the initial guess $\boldsymbol{\alpha}_0$ taken as
\begin{equation}
    \boldsymbol{\alpha}_0 =  \boldsymbol{\alpha}(\mathbf{Q}).
\end{equation}
In practice, since $\boldsymbol{\alpha}_0$ is generally a good guess for even moderately resolved velocity and internal energy spaces, very few iterations are required to bring the residual down to machine precision levels, with 1-2 iterations generally being sufficient.

\subsection{Positivity-Preserving Limiter}\label{ssec:pp}
As the distribution function $f$ represents the number density of particles, it is inherently a non-negative quantity. However, as the high-order spatial scheme does not preserve a maximum principle on its own, a distribution function that is non-negative at $t = 0$ may not remain non-negative for $t > 0$. In addition to the fact that negative values of $f$ are non-physical, this effect can result in stability problems. It is trivial to show that if $f$ is strictly non-negative, the resulting density and specific internal energy are guaranteed to be non-negative given a discrete integration operator $\mathbf{M}$ with all positive values. The issue arises when $f$ is no longer non-negative -- then the density or specific internal energy values may become negative for which the equilibrium  Maxwellian distribution is severely ill-behaved, causing the scheme to diverge. Therefore, for the purpose of robustness, it is necessary to modify the spatial scheme to ensure that $f$ remains non-negative.

In this work, this property is enforced similarly to the work of \citet{Jaiswal2022} by using the "squeeze" limiter of \citet{Zhang2010}. Within each element $\Omega_k$, the spatial average of the distribution function for each velocity and internal energy node $\mathbf{u}_q$ and $\zeta_r$ is computed as
\begin{equation}
    \bar{f}_{q,r} = \sum_{i=1}^{N_s} m_i f_{q,r}(\mathbf{x}_i),
\end{equation}
where the subscripts $q,r$ denote the $q$-th velocity node and $r$-th internal energy node, respectively, and $m_i$ are a set of quadrature weights such that
\begin{equation}
    \sum_{i=1}^{N_s} m_i z(\mathbf{x}_i) = \frac{\int_{\Omega_k} z(\mathbf{x})\ \mathrm{d}\mathbf{x}}{\int_{\Omega_k} \ \mathrm{d}\mathbf{x}},
\end{equation}
for any polynomial $z(\mathbf{x})$ of degree $\leq p$. In \citet{Zhang2010} and \citet{Zhang2011}, it was shown that the maximum principle is preserved for $\bar{f}$ under a standard CFL condition with strong-stability preserving temporal integration and an upwind interface flux. Due to the positivity of the equilibrium distribution function $g$ in the source term of the Boltzmann--BGK equation, this property can be extended to the element-wise mean of high-order flux reconstruction approximations of the Boltzmann--BGK equation under the condition $\Delta t \leq \tau$ as this results in a convex combination of two positive states $g$ and $f$. With the maximum principle property on $\bar{f}$, a conservative, positivity-preserving, and formally high-order reconstruction of $f$ can be given as 
\begin{equation}
    \hat{f}_{q,r}(\mathbf{x}_i) = \bar{f}_{q,r}(\mathbf{x}_i) + \beta_{q,r}\left[ f_{q,r}(\mathbf{x}_i) - \bar{f}_{q,r} \right],
\end{equation}
where
\begin{equation}
    \beta_{q,r} = \min \left [ \left | \frac{\bar{f}_{q,r}}{\bar{f}_{q,r} - f_{q,r}^{\min}}\right |, 1\right] \quad \mathrm{and} \quad  f_{q,r}^{\min} = \min f_{q,r}(\mathbf{x}_i)\ \forall \ i \in \{1, ..., N_s\}.
\end{equation}
This limiting is performed on a per-node basis in the velocity and internal energy spaces prior to each evaluation of the time derivative of the distribution function. Note that unlike some systems where an arbitrarily small positive constant must be given as a minimum to prevent numerical complications (e.g., Euler equations), $f$ is allowed to attain a zero value, such that there is no need to prescribe an arbitrary minimum. 

\subsection{Discontinuity Capturing}\label{ssec:disc}
Although the limiting approach in \cref{ssec:pp} is guaranteed to preserve the positivity of density and internal energy, the use of a high-order spatial discretization may still result in the introduction of spurious oscillations in the vicinity of discontinuities without proper treatment. For the Boltzmann equation, this is generally less of a problem as the discontinuous structures can in fact be resolved with some finite thickness at the particle mean free path level \citep{Xiao2021}. However, if the numerical resolution is not fine enough such that these structures are resolved, they can behave as discontinuities which can lead to numerical instabilities in the solution. For many high-order numerical schemes for both continuum gas dynamics and kinetics, it is commonplace to modify the parameters of the system in the vicinity of discontinuities such that they can be resolved on the order of the mesh scale (e.g., artificial viscosity \citep{Xiao2021}). However, as the focus of this work is on \textit{directly resolving} flow physics, we instead pose some resolution requirements of the scheme such that the structures can be resolved without any additional treatment. 

For a given numerical resolution, there exists some maximal mesh scale $h_{\max}$ such that features of a smaller scale cannot be resolved. For a high-order discontinuous spectral element method, computing this scale is more ambiguous due to the presence of multiple solution points within an element and inhomogeneity in the distribution of solution points. In this work, we compute this scale as the maximal distance between any solution point and its Voronoi neighbors within the element (see Section 2.4 of \citet{Guermond2011}). To resolve "discontinuous" structures, it is necessary for their thickness $\Delta$ to be greater than this maximal mesh scale, i.e.,
\begin{equation}\label{eq:shockmesh}
    \Delta \geq h_{\max}.
\end{equation}
To form an \textit{a priori} estimate of this thickness, we follow the methodology of \citet{Xiao2021}, where the thickness is calculated in terms of the molecular mean free path $\lambda$ as 
\begin{equation}\label{eq:shockthick}
    \Delta \sim 10 \lambda,
\end{equation}
using the relation that the thickness of a weak shock wave is on the order of 10 mean free paths. Note that for strong shocks or other discontinuities, this relation may vary but it is generally of the same order of magnitude \citep{Alsmeyer1976}. A mesh Knudsen number may then be calculated as
\begin{equation}\label{eq:meshknudsen}
    Kn_h = \frac{\lambda}{h_{\max}} = \frac{Kn L_0}{h_{\max}}. 
\end{equation}
From \cref{eq:shockmesh} and \cref{eq:shockthick}, an approximate resolution requirement can be posed in terms of this mesh Knudsen number as  
\begin{equation}\label{eq:resreq}
    Kn_h \geq 1/10.
\end{equation}
In practice, this requirement is checked \textit{a priori} based on the initial Knudsen number but may also be observed over the course of the simulation based on the local mesh scale and local Knudsen number. Note that this requirement does not \textit{guarantee} that the resulting scheme will be monotone in the vicinity of discontinuities, particularly due to the ambiguity in calculating $h_{\max}$ and the shock thickness, but it does form a suitable estimate such that the predicted results are generally well-behaved. 

\subsection{Boundary Conditions}
The enforcement of boundary conditions in the solution of the Boltzmann equation is nontrivial, and there exist a variety of approaches for dealing with this problem. In this work, the boundary conditions are \textit{weakly} enforced via the standard FR approach of forming a boundary state $f^+(\mathbf{u})$ that is used in conjunction with the interior solution values $f^-(\mathbf{u})$ to calculate the interface flux at the boundary. As the interface flux is computed using an upwinding approach, the boundary state can only affect the incoming particles at the boundary while the outgoing particles are essentially left unmodified. 

The presentation of these boundary conditions will generally follow the work of \citet{Evans2011}. The most straightforward boundary condition to enforce is a Neumann-type boundary condition, which can be implemented by simply setting the boundary state as 
\begin{equation}
    f^+(\mathbf{u}) = f^-(\mathbf{u}).
\end{equation}
This boundary condition essentially assumes that the boundary has no effect on the solution, such that the correction terms associated with the interface flux are zero. For Dirichlet-type boundary conditions, the implementation becomes slightly more involved. It can be assumed in most scenarios that the boundary state can be represented in terms of the macroscopic variables $\mathbf{Q}$. The simplest case for this boundary condition is to assume that the boundary state is in thermodynamic equilibrium, such that the distribution function can be represented as a Maxwellian. Therefore, the boundary state can be set as 
\begin{equation}
    f^+(\mathbf{u}) = g(\mathbf{Q}'),
\end{equation}
where $\mathbf{Q}'$ is the modified Maxwellian used to ensure a conservative scheme. Since the boundary state is typically fixed, this modified Maxwellian can be pre-computed. 

Wall boundary conditions for the Boltzmann equation are significantly more complex, and in most cases, it is still an open problem as to what the proper approach is for representing interactions of the particles with the wall. Some theoretical development on wall boundary conditions for the Boltzmann equation is presented in \citet{Williams2001}, and numerical applications of various approaches for wall boundary conditions are shown in \citet{Evans2011}. The focus of this work is strictly on specular wall boundary conditions for convex surfaces which model the wall interaction as a reflection of the incoming particles, mimicking an effect that is somewhat similar to a slip-wall boundary condition for continuum approximations. For this boundary condition, the boundary state is set as 
\begin{equation}\label{eq:bcspecular}
    f^+(\mathbf{u}) = f^-(\mathbf{u} - 2(\mathbf{u} \cdot \mathbf{n}) \mathbf{n}),
\end{equation}
where $\mathbf{n}$ is the outward-facing normal direction of the wall. In general, the implementation of this boundary condition requires interpolation in velocity space, but for simple geometries, the implementation can be vastly simplified by use of a velocity space that is symmetric with respect to the wall-tangent directions. 

%% file: figs/scheme.tex
     \begin{tikzpicture}[spy using outlines={rectangle, height=3cm,width=2.3cm, magnification=3, connect spies}]
		\begin{axis}[name=plot1,
		    axis line style={draw=none},
		    tick style={draw=none},
		    axis x line=left,
            axis y line=left,
            axis equal image,
            clip mode=individual,
    		xmin=-2,
    		xmax=1,
    		xticklabels={,,},
    		ymin=-1,
    		ymax=1,
    		yticklabels={,,},
    		style={font=\Large},
    		scale = 1]
    		
            \draw [black, very thick] plot [] coordinates {(0, -1)  (0, 1) (-1.5, 0) (0, -1)};
            \fill[fill=black!10] (0, -1)--(0, 1)--(-1.5, 0);
            
            \draw [black, very thick, dotted] plot [] coordinates {(-.05, 1.1) (-1.55, 0.1) (-1.9, 1.6) (-.05, 1.1)};
            \draw [black, very thick, dotted] plot [] coordinates {(-1.58, -0.07) (0, -1.13) (-2.5, -1.13) (-1.58, -0.07)};
            \draw [black, very thick, dotted] plot [] coordinates {(-2.55, -1.0) (-2, 1.6) (-1.65, 0.03) (-2.55, -1.0)};
		  
            \draw[-,very thick, color=cyan!80!blue,fill=red] (0, -1) circle[radius=0.040];
            \draw[-,very thick, color=cyan!80!blue,fill=red] (0, 0) circle[radius=0.040];
            \draw[-,very thick, color=cyan!80!blue,fill=red] (0, 1) circle[radius=0.040];
            \draw[-,very thick, color=cyan!80!blue,fill=red] (-0.75, -0.5) circle[radius=0.040];
            \draw[-,very thick, color=cyan!80!blue,fill=red] (-1.5, 0) circle[radius=0.040];
            \draw[-,very thick, color=cyan!80!blue,fill=red] (-0.75, 0.5) circle[radius=0.040];
            \draw[->, thick](-0.5, 0) to [out=75,in=165] (0.7, 0.6);
            \draw[-,very thick, color=black, fill=red] (-0.5, 0) circle[radius=0.040];
            
            \fill[-,fill=cyan!80!blue] (-.05, 1.1)  circle[radius=0.045];
            \fill[-,fill=cyan!80!blue] (-1.55, 0.1)  circle[radius=0.045];
            \fill[-,fill=cyan!80!blue] (-.8,0.6)  circle[radius=0.045];
            
            \fill[-,fill=cyan!80!blue] (-1.58, -0.07)  circle[radius=0.045];
            \fill[-,fill=cyan!80!blue] (-0.79, -0.6)  circle[radius=0.045];
            \fill[-,fill=cyan!80!blue] (0, -1.13) circle[radius=0.045];

            \draw [->, black, thick] plot [] coordinates {(2.0, 0.3) (3.35, 0.3)};
            \draw [->, black, thick] plot [] coordinates {(2.0, 0.3) (2.0, 1.65)};
            \draw [->, black, thick] plot [] coordinates {(2.0, 0.3) (0.65, 0.3)};
            \draw [->, black, thick] plot [] coordinates {(2.0, 0.3) (2.0, -1.05)};
            \draw [->, black, thick, dotted] plot [] coordinates {(2.0, 0.3) (2.69446279, 1.33933702)};
            
            \node[] at (-1.25, -1.3) {$\Omega^{\mathbf{x}}$};
            \node[] at (-0.2, -0.5) {\small$\Omega^{\mathbf{x}}_k$};
            \node[] at (2.05, -1.3) {$\Omega^{\mathbf{u}}$};
            \node[] at (3.425, 0.3) {\small$u$};
            \node[] at (2.0, 1.725) {\small$v$};
            \node[] at (2.83, 1.4) {\small$r_{\max}$};

            \fill[-,fill=orange] (2.02481883968904, 0.3) circle[radius=0.04];
            \fill[-,fill=orange] (2.0229296180093828, 0.3094977587595207) circle[radius=0.04];
            \fill[-,fill=orange] (2.017549569845302, 0.3175495698453019) circle[radius=0.04];
            \fill[-,fill=orange] (2.0094977587595206, 0.32292961800938275) circle[radius=0.04];
            \fill[-,fill=orange] (2.0, 0.3248188396890399) circle[radius=0.04];
            \fill[-,fill=orange] (1.9905022412404794, 0.32292961800938275) circle[radius=0.04];
            \fill[-,fill=orange] (1.9824504301546981, 0.3175495698453019) circle[radius=0.04];
            \fill[-,fill=orange] (1.9770703819906172, 0.3094977587595207) circle[radius=0.04];
            \fill[-,fill=orange] (1.9751811603109601, 0.3) circle[radius=0.04];
            \fill[-,fill=orange] (1.9770703819906172, 0.2905022412404793) circle[radius=0.04];
            \fill[-,fill=orange] (1.9824504301546981, 0.28245043015469806) circle[radius=0.04];
            \fill[-,fill=orange] (1.9905022412404794, 0.27707038199061723) circle[radius=0.04];
            \fill[-,fill=orange] (2.0, 0.27518116031096007) circle[radius=0.04];
            \fill[-,fill=orange] (2.0094977587595206, 0.27707038199061723) circle[radius=0.04];
            \fill[-,fill=orange] (2.017549569845302, 0.28245043015469806) circle[radius=0.04];
            \fill[-,fill=orange] (2.0229296180093828, 0.2905022412404793) circle[radius=0.04];
            \fill[-,fill=orange] (2.127083451616483, 0.3) circle[radius=0.04];
            \fill[-,fill=orange] (2.1174097998693573, 0.3486327314613986) circle[radius=0.04];
            \fill[-,fill=orange] (2.089861570414608, 0.38986157041460784) circle[radius=0.04];
            \fill[-,fill=orange] (2.0486327314613986, 0.41740979986935733) circle[radius=0.04];
            \fill[-,fill=orange] (2.0, 0.4270834516164833) circle[radius=0.04];
            \fill[-,fill=orange] (1.9513672685386014, 0.41740979986935733) circle[radius=0.04];
            \fill[-,fill=orange] (1.910138429585392, 0.38986157041460784) circle[radius=0.04];
            \fill[-,fill=orange] (1.8825902001306427, 0.34863273146139867) circle[radius=0.04];
            \fill[-,fill=orange] (1.8729165483835166, 0.3) circle[radius=0.04];
            \fill[-,fill=orange] (1.8825902001306427, 0.25136726853860136) circle[radius=0.04];
            \fill[-,fill=orange] (1.910138429585392, 0.21013842958539214) circle[radius=0.04];
            \fill[-,fill=orange] (1.9513672685386012, 0.1825902001306427) circle[radius=0.04];
            \fill[-,fill=orange] (2.0, 0.1729165483835167) circle[radius=0.04];
            \fill[-,fill=orange] (2.0486327314613986, 0.1825902001306427) circle[radius=0.04];
            \fill[-,fill=orange] (2.089861570414608, 0.2101384295853921) circle[radius=0.04];
            \fill[-,fill=orange] (2.1174097998693573, 0.25136726853860125) circle[radius=0.04];
            \fill[-,fill=orange] (2.2965422438022944, 0.3) circle[radius=0.04];
            \fill[-,fill=orange] (2.2739693095739115, 0.41348180369950727) circle[radius=0.04];
            \fill[-,fill=orange] (2.2096870315008768, 0.5096870315008768) circle[radius=0.04];
            \fill[-,fill=orange] (2.113481803699507, 0.5739693095739118) circle[radius=0.04];
            \fill[-,fill=orange] (2.0, 0.5965422438022944) circle[radius=0.04];
            \fill[-,fill=orange] (1.8865181963004927, 0.5739693095739118) circle[radius=0.04];
            \fill[-,fill=orange] (1.7903129684991232, 0.5096870315008768) circle[radius=0.04];
            \fill[-,fill=orange] (1.7260306904260883, 0.4134818036995073) circle[radius=0.04];
            \fill[-,fill=orange] (1.7034577561977056, 0.30000000000000004) circle[radius=0.04];
            \fill[-,fill=orange] (1.7260306904260883, 0.18651819630049274) circle[radius=0.04];
            \fill[-,fill=orange] (1.7903129684991232, 0.09031296849912321) circle[radius=0.04];
            \fill[-,fill=orange] (1.8865181963004924, 0.0260306904260883) circle[radius=0.04];
            \fill[-,fill=orange] (2.0, 0.003457756197705608) circle[radius=0.04];
            \fill[-,fill=orange] (2.1134818036995076, 0.0260306904260883) circle[radius=0.04];
            \fill[-,fill=orange] (2.2096870315008768, 0.09031296849912315) circle[radius=0.04];
            \fill[-,fill=orange] (2.2739693095739115, 0.18651819630049252) circle[radius=0.04];
            \fill[-,fill=orange] (2.5103533484402187, 0.3) circle[radius=0.04];
            \fill[-,fill=orange] (2.471505012972519, 0.4953037711001196) circle[radius=0.04];
            \fill[-,fill=orange] (2.3608743134833396, 0.6608743134833397) circle[radius=0.04];
            \fill[-,fill=orange] (2.1953037711001198, 0.7715050129725192) circle[radius=0.04];
            \fill[-,fill=orange] (2.0, 0.8103533484402188) circle[radius=0.04];
            \fill[-,fill=orange] (1.8046962288998805, 0.7715050129725192) circle[radius=0.04];
            \fill[-,fill=orange] (1.6391256865166604, 0.6608743134833397) circle[radius=0.04];
            \fill[-,fill=orange] (1.5284949870274809, 0.49530377110011964) circle[radius=0.04];
            \fill[-,fill=orange] (1.4896466515597813, 0.30000000000000004) circle[radius=0.04];
            \fill[-,fill=orange] (1.5284949870274807, 0.10469622889988045) circle[radius=0.04];
            \fill[-,fill=orange] (1.6391256865166604, -0.06087431348333966) circle[radius=0.04];
            \fill[-,fill=orange] (1.8046962288998802, -0.17150501297251913) circle[radius=0.04];
            \fill[-,fill=orange] (2.0, -0.21035334844021886) circle[radius=0.04];
            \fill[-,fill=orange] (2.1953037711001198, -0.1715050129725192) circle[radius=0.04];
            \fill[-,fill=orange] (2.3608743134833396, -0.06087431348333977) circle[radius=0.04];
            \fill[-,fill=orange] (2.471505012972519, 0.10469622889988009) circle[radius=0.04];
            \fill[-,fill=orange] (2.739646651559781, 0.3) circle[radius=0.04];
            \fill[-,fill=orange] (2.683344402666589, 0.5830505193562425) circle[radius=0.04];
            \fill[-,fill=orange] (2.523009162999845, 0.8230091629998446) circle[radius=0.04];
            \fill[-,fill=orange] (2.2830505193562427, 0.9833444026665892) circle[radius=0.04];
            \fill[-,fill=orange] (2.0, 1.039646651559781) circle[radius=0.04];
            \fill[-,fill=orange] (1.7169494806437575, 0.9833444026665892) circle[radius=0.04];
            \fill[-,fill=orange] (1.4769908370001554, 0.8230091629998446) circle[radius=0.04];
            \fill[-,fill=orange] (1.3166555973334109, 0.5830505193562427) circle[radius=0.04];
            \fill[-,fill=orange] (1.260353348440219, 0.3000000000000001) circle[radius=0.04];
            \fill[-,fill=orange] (1.3166555973334109, 0.016949480643757453) circle[radius=0.04];
            \fill[-,fill=orange] (1.4769908370001552, -0.22300916299984458) circle[radius=0.04];
            \fill[-,fill=orange] (1.7169494806437569, -0.3833444026665889) circle[radius=0.04];
            \fill[-,fill=orange] (1.9999999999999998, -0.43964665155978105) circle[radius=0.04];
            \fill[-,fill=orange] (2.2830505193562427, -0.38334440266658903) circle[radius=0.04];
            \fill[-,fill=orange] (2.523009162999845, -0.2230091629998448) circle[radius=0.04];
            \fill[-,fill=orange] (2.683344402666589, 0.016949480643756953) circle[radius=0.04];
            \fill[-,fill=orange] (2.9534577561977056, 0.3) circle[radius=0.04];
            \fill[-,fill=orange] (2.8808801060651965, 0.664872486756855) circle[radius=0.04];
            \fill[-,fill=orange] (2.6741964449823077, 0.9741964449823075) circle[radius=0.04];
            \fill[-,fill=orange] (2.364872486756855, 1.1808801060651966) circle[radius=0.04];
            \fill[-,fill=orange] (2.0, 1.2534577561977056) circle[radius=0.04];
            \fill[-,fill=orange] (1.635127513243145, 1.1808801060651966) circle[radius=0.04];
            \fill[-,fill=orange] (1.3258035550176925, 0.9741964449823075) circle[radius=0.04];
            \fill[-,fill=orange] (1.1191198939348035, 0.664872486756855) circle[radius=0.04];
            \fill[-,fill=orange] (1.0465422438022944, 0.3000000000000001) circle[radius=0.04];
            \fill[-,fill=orange] (1.1191198939348033, -0.06487248675685481) circle[radius=0.04];
            \fill[-,fill=orange] (1.3258035550176923, -0.3741964449823075) circle[radius=0.04];
            \fill[-,fill=orange] (1.6351275132431446, -0.5808801060651965) circle[radius=0.04];
            \fill[-,fill=orange] (1.9999999999999998, -0.6534577561977055) circle[radius=0.04];
            \fill[-,fill=orange] (2.3648724867568554, -0.5808801060651965) circle[radius=0.04];
            \fill[-,fill=orange] (2.6741964449823072, -0.3741964449823077) circle[radius=0.04];
            \fill[-,fill=orange] (2.8808801060651965, -0.06487248675685553) circle[radius=0.04];
            \fill[-,fill=orange] (3.1229165483835164, 0.3) circle[radius=0.04];
            \fill[-,fill=orange] (3.0374396157697507, 0.7297215589949635) circle[radius=0.04];
            \fill[-,fill=orange] (2.7940219060685765, 1.0940219060685765) circle[radius=0.04];
            \fill[-,fill=orange] (2.4297215589949634, 1.337439615769751) circle[radius=0.04];
            \fill[-,fill=orange] (2.0, 1.4229165483835167) circle[radius=0.04];
            \fill[-,fill=orange] (1.5702784410050366, 1.337439615769751) circle[radius=0.04];
            \fill[-,fill=orange] (1.2059780939314235, 1.0940219060685765) circle[radius=0.04];
            \fill[-,fill=orange] (0.962560384230249, 0.7297215589949637) circle[radius=0.04];
            \fill[-,fill=orange] (0.8770834516164834, 0.3000000000000001) circle[radius=0.04];
            \fill[-,fill=orange] (0.9625603842302488, -0.12972155899496346) circle[radius=0.04];
            \fill[-,fill=orange] (1.2059780939314233, -0.4940219060685765) circle[radius=0.04];
            \fill[-,fill=orange] (1.5702784410050359, -0.7374396157697507) circle[radius=0.04];
            \fill[-,fill=orange] (1.9999999999999998, -0.8229165483835166) circle[radius=0.04];
            \fill[-,fill=orange] (2.429721558994964, -0.7374396157697509) circle[radius=0.04];
            \fill[-,fill=orange] (2.7940219060685765, -0.4940219060685767) circle[radius=0.04];
            \fill[-,fill=orange] (3.0374396157697507, -0.12972155899496424) circle[radius=0.04];
            \fill[-,fill=orange] (3.22518116031096, 0.3) circle[radius=0.04];
            \fill[-,fill=orange] (3.1319197976297257, 0.7688565316968415) circle[radius=0.04];
            \fill[-,fill=orange] (2.8663339066378826, 1.1663339066378826) circle[radius=0.04];
            \fill[-,fill=orange] (2.468856531696842, 1.4319197976297258) circle[radius=0.04];
            \fill[-,fill=orange] (2.0, 1.5251811603109602) circle[radius=0.04];
            \fill[-,fill=orange] (1.5311434683031586, 1.4319197976297258) circle[radius=0.04];
            \fill[-,fill=orange] (1.1336660933621174, 1.1663339066378826) circle[radius=0.04];
            \fill[-,fill=orange] (0.8680802023702743, 0.7688565316968416) circle[radius=0.04];
            \fill[-,fill=orange] (0.7748188396890399, 0.30000000000000016) circle[radius=0.04];
            \fill[-,fill=orange] (0.8680802023702741, -0.1688565316968414) circle[radius=0.04];
            \fill[-,fill=orange] (1.1336660933621174, -0.5663339066378825) circle[radius=0.04];
            \fill[-,fill=orange] (1.5311434683031577, -0.8319197976297255) circle[radius=0.04];
            \fill[-,fill=orange] (1.9999999999999998, -0.9251811603109601) circle[radius=0.04];
            \fill[-,fill=orange] (2.468856531696842, -0.8319197976297255) circle[radius=0.04];
            \fill[-,fill=orange] (2.866333906637882, -0.5663339066378827) circle[radius=0.04];
            \fill[-,fill=orange] (3.1319197976297257, -0.16885653169684228) circle[radius=0.04];
		\end{axis}

	\end{tikzpicture}

%% file: implementation.tex
\section{Implementation}\label{sec:implementation}
The proposed scheme was implemented within PyFR \citep{Witherden2014}, a high-order flux reconstruction solver that can target massively-parallel CPU and GPU computing architectures, and simulations were performed on up to 40 NVIDIA V100 GPUs. For the spatial discretization, the solution nodes were placed on the Gauss--Lobatto quadrature nodes for tensor-product elements and $\alpha$-optimized points \citep{Hesthaven2008DG} for simplex elements, but any closed nodal set would essentially be identical due to the linearity of the flux. Temporal integration was performed via an explicit, fourth-order, four-stage Runge--Kutta scheme. Although the positivity property of the element-wise mean required by the positivity-preserving limiter theoretically requires strong stability preserving temporal integration, this is almost never an issue in practice, and it is therefore advantageous to utilize more efficient temporal schemes. The time step was chosen as the minima of the collision time scale $\tau$ and the maximum allowable time step by the Courant–-Friedrichs–-Lewy (CFL) condition $\Delta t_{CFL}$, i.e.,
\begin{equation}
    \Delta t = \min \left (\tau, \Delta t_{CFL} \right),
\end{equation}
where $\Delta t_{CFL}$ was computed using the estimate of \citet{Cockburn2001} as 
\begin{equation}
    \Delta t_{CFL} = \frac{CFL}{2p+1}\frac{h_{\min}}{c_{\max}}.
\end{equation}
Here, $CFL = 0.5$ is the chosen CFL number, $p$ is the order of the scheme, $h_{\min}$ is the minimum element edge length in the domain, and $c_{\max} = r_{\max}$ is the maximum particle velocity in the domain. For certain problems, such as ones in the zero Mach number limit or ones that converge to a steady state, it may be advantageous to utilize implicit time stepping due to the stiffness of the source term. However, as this work is primarily focused on problems tending towards direct numerical simulation, where the time step limit imposed by the source term usually does not drastically differ from the CFL-based time step, explicit time stepping was deemed preferable due to its significant benefits in terms of computational efficiency on massively-parallel GPU computing architectures. For the numerical experiments in this work, the time step was generally limited by the CFL condition except for cases in the continuum limit (i.e., low Mach number, high Reynolds number flows).

To simplify the presentation, the initial conditions are given in terms of the primitive variables $\mathbf q = [\rho, \mathbf{U}, P]^T$. The initial distribution function was set as the \textit{modified} Maxwellian corresponding to the initial macroscopic variables through 5 iterations of the DVM, after which the DVM residual was generally on the order of machine precision. Throughout the simulations, a fixed number of iterations (2) was used for the DVM, which will be later shown to be sufficient to ensure conservation of the macroscopic flow variables. The extent of the velocity domain was computed as presented in \cref{ssec:velocity} and the velocity tolerance was set as $\epsilon_{\mathbf{u}} = 10^{-15}$ as the estimated extent did not vary strongly with the tolerance (e.g., a velocity tolerance of $\epsilon_{\mathbf{u}} = 10^{-6}$ would only decrease the extent by 37\%). However, for cases with internal degrees of freedom, since the extent of the internal energy domain varied more strongly with the internal energy tolerance (as shown in \cref{tab:zeta_max}), the internal energy tolerance was instead set as $\epsilon_{\zeta} = 10^{-6}$. The effects of this truncation will be explored in \cref{sec:results}. The velocity discretization was generally chosen with the goal of minimizing the necessary resolution, particularly for larger-scale problems, usually through qualitative convergence studies performed by incrementally increasing the resolution and comparing with established results. 

Unless otherwise stated, the Knudsen number was computed with respect to a unit reference length and the discrete maximum of the initial wavespeed in the domain, and the velocity dimension $m$ was set equal to the spatial dimension $d$. Furthermore, the collision time $\tau$ was set as constant based on the Knudsen number in almost all scenarios, except for problems where a variable collision time model is explicitly defined. For some problems, comparisons were made between the Boltzmann--BGK approach and a standard Navier--Stokes approach implemented within an identical codebase. The Navier--Stokes results were computed using the HLLC \citep{Toro1994} Riemann solver for the inviscid fluxes and the BR2 approach of \citet{Bassi2000} for the viscous fluxes, and a unit Prandtl number and identical specific heat ratio was used for consistency with the Boltzmann--BGK approach.

%% file: results.tex
\section{Results}\label{sec:results}

\subsection{Smooth Pulse Propagation}
\begin{figure}[tbhp]
\centering
    \adjustbox{width=0.4\linewidth, valign=b}{\input{figs/pulse}}
    \caption{\label{fig:pulse} Density profile of the smooth pulse propagation problem at $t=1$ with varying Knudsen numbers computed using a $\mathbb P_5$ scheme with 20 elements, the discrete velocity model (DVM) with $N_v = 32$, and $\delta = 0$ ($\gamma = 3$).}
\end{figure}
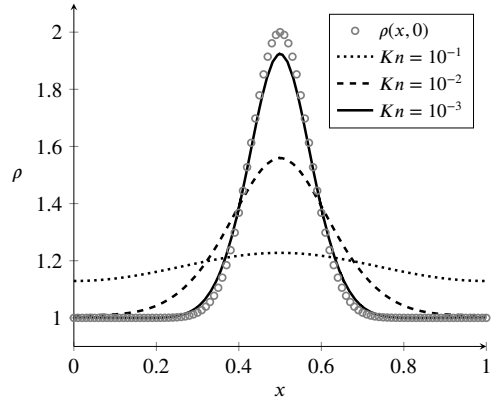

As an initial evaluation of the accuracy of the proposed scheme, the propagation of a smooth density pulse is considered. The problem is defined on a one-dimensional periodic domain $\Omega^x = [0,1]$ and the initial conditions are given as
\begin{equation}
\mathbf{q}(x, 0) = 
    \begin{bmatrix}
           \rho \\
           U \\
           P
    \end{bmatrix} 
    =
    \begin{bmatrix}
           1 + \exp\left(-\beta \left(x - 0.5\right)^2 \right)\\
           1 \\
           1
    \end{bmatrix},
\end{equation}
where $\beta = 100$ is a parameter controlling the pulse width. The problem was investigated with $Kn = 10^{-1}, 10^{-2},$ and $10^{-3}$ based on a unit characteristic length and the maximum wavespeed in the domain, with an example of the resulting density profiles shown in \cref{fig:pulse}. Additionally, the effects of including internal degrees of freedom were explored. Given a single velocity dimension ($m = 1$), a comparison was made between the monatomic case ($\delta = 0$, $\gamma = 3$) and the polyatomic case mimicking a diatomic molecule in three dimensions ($\delta = 4$, $\gamma = 1.4$).

\subsubsection{Spatial Convergence}
To verify the high-order spatial accuracy of the flux reconstruction approach, the convergence of the density error was evaluated. After one flow-through of the domain ($t=1$), the $L^\infty$ norm of the density error was calculated as
\begin{equation}
    \epsilon_{\rho, \infty} = \| \rho(x, t) - \rho_{\mathrm{ref}}(x, t) \|_{\infty,\Omega^{x}},
\end{equation}
where the extremum was computed on the discrete solution nodes and $\rho_{\mathrm{ref}}$ is the reference density computed with a highly-resolved numerical scheme. For this case, the reference simulation was performed using a $\mathbb P_5$ scheme with 100 elements. To isolate the effects of the spatial discretization error from the velocity discretization error, a very high resolution, $N_v = 128$, was used for the velocity space for all tests with $N_\zeta = N_v$ for the polyatomic case. The convergence of the density error for $Kn = 10^{-1}, 10^{-2},$ and $10^{-3}$ is shown in \cref{tab:spatial_errorkn-1}, \cref{tab:spatial_errorkn-2}, and \cref{tab:spatial_errorkn-3}, respectively. For almost all cases, the expected theoretical $p+1$ convergence rate was recovered or exceeded, although some variation is expected due to the lack of an analytic solution. These results confirm the high-order spatial accuracy afforded by the flux reconstruction scheme. 

\begin{figure}[tbhp]
    \centering
    \begin{tabularx}{\textwidth}{r | c c c c | c c c c}
         & \multicolumn{4}{c}{$\delta = 0$ ($\gamma = 3$)} & \multicolumn{4}{|c}{$\delta = 4$ ($\gamma = 1.4$)} \\ \midrule
        $\Delta x$ & $\mathbb{P}_2$ &$\mathbb{P}_3$ &$\mathbb{P}_4$ &$\mathbb{P}_5$ & $\mathbb{P}_2$ &$\mathbb{P}_3$ &$\mathbb{P}_4$ &$\mathbb{P}_5$ \\ \midrule
        $1/4$ & \num{3.37e-02} & \num{1.97e-03} & \num{1.20e-04} & \num{6.61e-05} & \num{3.76e-03} & \num{4.91e-04} & \num{1.98e-04} & \num{2.96e-05} \\
        $1/8$ & \num{7.03e-04} & \num{5.42e-05} & \num{2.61e-06} & \num{7.24e-08} & \num{7.46e-04} & \num{5.55e-05} & \num{2.80e-06} & \num{1.59e-07} \\
        $1/12$ & \num{4.52e-05} & \num{1.49e-06} & \num{5.51e-08} & \num{1.69e-09} & \num{2.22e-04} & \num{1.19e-05} & \num{4.67e-07} & \num{1.88e-08} \\
        $1/16$ & \num{1.84e-05} & \num{4.75e-07} & \num{1.10e-08} & \num{2.81e-10} & \num{9.68e-05} & \num{3.91e-06} & \num{1.05e-07} & \num{3.61e-09} \\
        $1/20$ & \num{9.04e-06} & \num{1.91e-07} & \num{3.61e-09} & \num{7.27e-11} & \num{5.00e-05} & \num{1.63e-06} & \num{3.64e-08} & \num{9.58e-10} \\
        \midrule
        RoC & $5.26$ & $5.98$ & $6.71$ & $8.63$ &$2.70$ & $3.55$ & $5.31$ & $6.35$ \\
    \end{tabularx}
    \captionof{table}{\label{tab:spatial_errorkn-1}  Convergence of the $L^\infty$ norm of the density error with respect to mesh resolution $\Delta x$ and varying approximations orders for the smooth pulse propagation problem with $Kn = 10^{-1}$. Rate of convergence shown on bottom.}
\end{figure}
\begin{figure}[tbhp]
    \centering
    \begin{tabularx}{\textwidth}{r | c c c c | c c c c}
         & \multicolumn{4}{c}{$\delta = 0$ ($\gamma = 3$)} & \multicolumn{4}{|c}{$\delta = 4$ ($\gamma = 1.4$)} \\ \midrule
        $\Delta x$ & $\mathbb{P}_2$ &$\mathbb{P}_3$ &$\mathbb{P}_4$ &$\mathbb{P}_5$ & $\mathbb{P}_2$ &$\mathbb{P}_3$ &$\mathbb{P}_4$ &$\mathbb{P}_5$ \\ \midrule
        $1/4$ & \num{1.43e-01} & \num{1.63e-02} & \num{3.33e-03} & \num{1.20e-03} & \num{3.06e-02} & \num{5.36e-03} & \num{5.93e-03} & \num{5.28e-04} \\
        $1/8$ & \num{1.59e-02} & \num{9.56e-04} & \num{1.30e-04} & \num{1.15e-05} & \num{5.92e-03} & \num{4.87e-04} & \num{5.59e-05} & \num{4.10e-06} \\
        $1/12$ & \num{2.85e-03} & \num{1.82e-04} & \num{1.54e-05} & \num{1.34e-06} & \num{1.88e-03} & \num{9.30e-05} & \num{6.53e-06} & \num{4.08e-07} \\
        $1/16$ & \num{9.67e-04} & \num{5.76e-05} & \num{3.90e-06} & \num{2.88e-07} & \num{8.19e-04} & \num{2.94e-05} & \num{1.68e-06} & \num{9.55e-08} \\
        $1/20$ & \num{4.55e-04} & \num{2.36e-05} & \num{1.21e-06} & \num{8.12e-08} & \num{4.33e-04} & \num{1.19e-05} & \num{5.80e-07} & \num{2.93e-08} \\
        \midrule
        RoC & $3.63$ & $4.06$ & $4.93$ & $5.92$ & $2.65$ & $3.81$ & $5.72$ & $6.07$  \\
    \end{tabularx}
    \captionof{table}{\label{tab:spatial_errorkn-2}  Convergence of the $L^\infty$ norm of the density error with respect to mesh resolution $\Delta x$ and varying approximations orders for the smooth pulse propagation problem with $Kn = 10^{-2}$. Rate of convergence shown on bottom.}
\end{figure}
\begin{figure}[tbhp]
    \centering
    \begin{tabularx}{\textwidth}{r | c c c c | c c c c}
         & \multicolumn{4}{c}{$\delta = 0$ ($\gamma = 3$)} & \multicolumn{4}{|c}{$\delta = 4$ ($\gamma = 1.4$)} \\ \midrule
        $\Delta x$ & $\mathbb{P}_2$ &$\mathbb{P}_3$ &$\mathbb{P}_4$ &$\mathbb{P}_5$ & $\mathbb{P}_2$ &$\mathbb{P}_3$ &$\mathbb{P}_4$ &$\mathbb{P}_5$ \\ \midrule
        $1/4$ & \num{3.54e-01} & \num{1.32e-01} & \num{4.69e-02} & \num{1.57e-02} & \num{5.07e-02} & \num{9.73e-03} & \num{5.09e-03} & \num{1.52e-03} \\
        $1/8$ & \num{1.03e-01} & \num{1.24e-02} & \num{2.16e-03} & \num{2.59e-04} & \num{6.34e-03} & \num{3.75e-04} & \num{2.97e-05} & \num{6.76e-06} \\
        $1/12$ & \num{3.25e-02} & \num{1.98e-03} & \num{2.25e-04} & \num{1.83e-05} & \num{2.10e-03} & \num{7.56e-05} & \num{4.11e-06} & \num{6.39e-07} \\
        $1/16$ & \num{1.08e-02} & \num{6.20e-04} & \num{4.76e-05} & \num{4.53e-06} & \num{9.09e-04} & \num{1.95e-05} & \num{1.16e-06} & \num{2.68e-07} \\
        $1/20$ & \num{4.37e-03} & \num{2.32e-04} & \num{1.76e-05} & \num{1.39e-06} & \num{4.62e-04} & \num{8.86e-06} & \num{5.24e-07} & \num{1.24e-07} \\
        \midrule
        RoC & $2.71$ & $3.96$ & $4.98$ & $5.85$ & $2.90$ & $4.37$ & $5.71$ & $5.89$\\
    \end{tabularx}
    \captionof{table}{\label{tab:spatial_errorkn-3} Convergence of the $L^\infty$ norm of the density error with respect to mesh resolution $\Delta x$ and varying approximations orders for the smooth pulse propagation problem with $Kn = 10^{-3}$. Rate of convergence shown on bottom.}
\end{figure}

\subsubsection{Velocity Convergence}
\begin{figure}[tbhp]
        \subfloat[$Kn = 10^{-1}$]{\adjustbox{width=0.33\linewidth, valign=b}{\input{figs/velocity_convergence_kn0p1}}}
        \subfloat[$Kn = 10^{-2}$]{\adjustbox{width=0.33\linewidth, valign=b}{\input{figs/velocity_convergence_kn0p01}}}
        \subfloat[$Kn = 10^{-3}$]{\adjustbox{width=0.33\linewidth, valign=b}{\input{figs/velocity_convergence_kn0p001}}}
        \newline
        \caption{\label{fig:velocity_error} Convergence of the $L^\infty$ norm of the density error with respect to velocity space resolution $N_v$ computed using a $\mathbb P_5$ scheme with 20 elements and $Kn = 10^{-1}$ (left), $Kn = 10^{-2}$ (middle), and $Kn = 10^{-3}$ (right) with and without the discrete velocity model (DVM). Results for $\delta = 0$ $(\gamma = 3)$ and $\delta = 4$ $(\gamma = 1.4)$ shown in black and red, respectively. Diverged solutions not plotted.}
    \end{figure}
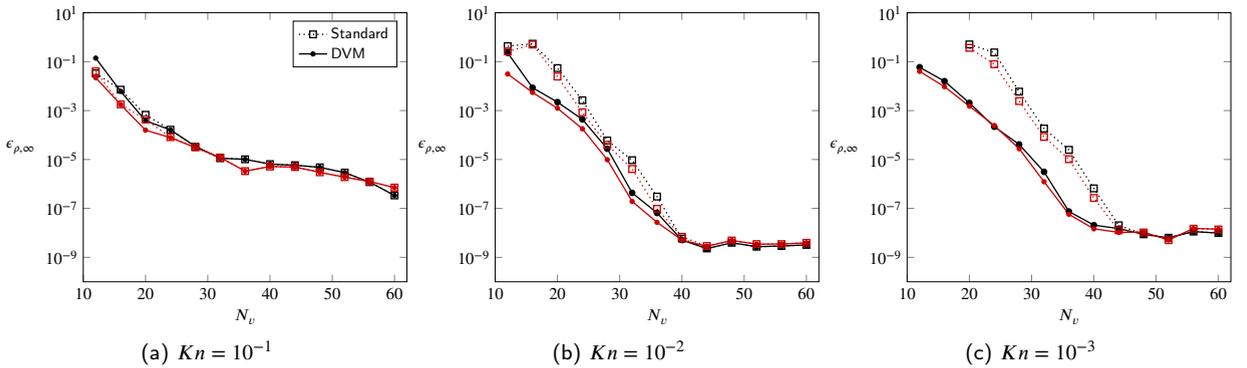

To evaluate the effects of the velocity/internal energy discretization as well as the discrete velocity model on the solution, the density error was similarly analyzed while varying the velocity/internal energy resolution. To isolate the effects of the velocity/internal energy discretization error from the spatial discretization error, a highly-resolved $\mathbb P_5$ scheme with 20 elements was used for all tests. The $L^\infty$ norm of the density error with respect to the resolution $N_v$, $N_\zeta$ after one flow-through of the domain is shown in \cref{fig:velocity_error} for $Kn = 10^{-1}, 10^{-2},$ and $10^{-3}$ with and without the DVM and internal degrees of freedom, respectively. For $Kn = 10^{-1}$, the results of the DVM and the standard approach were essentially identical over the resolution range. However, for smaller $Kn$, the difference between the standard approach and the DVM became more evident. For $Kn = 10^{-2}$ and $Kn = 10^{-3}$, the use of the DVM resulted in 1-3 orders of magnitude lower error than the standard approach over a large range of resolution, with this disparity becoming more evident with lower $Kn$. Furthermore, the use of the DVM sometimes resulted in a lower minimal stable resolution requirement for the velocity/internal energy spaces, with an example being the case of $Kn = 10^{-3}$ where the DVM was stable even with lowest resolution ($N_v, N_\zeta = 12$) while the standard approach diverged for $N_v, N_\zeta < 20$.

\begin{figure}[tbhp]
        \subfloat[$Kn = 10^{-1}$]{\adjustbox{width=0.33\linewidth, valign=b}{\input{figs/conservation_kn0p1}}}
        \subfloat[$Kn = 10^{-2}$]{\adjustbox{width=0.33\linewidth, valign=b}{\input{figs/conservation_kn0p01}}}
        \subfloat[$Kn = 10^{-3}$]{\adjustbox{width=0.33\linewidth, valign=b}{\input{figs/conservation_kn0p001}}}
        \newline
        \caption{\label{fig:velocity_masscon} Convergence of the mass conservation error with respect to velocity space resolution $N_v$ computed using a $\mathbb P_5$ scheme with 20 elements and $Kn = 10^{-1}$ (left), $Kn = 10^{-2}$ (middle), and $Kn = 10^{-3}$ (right) with and without the discrete velocity model (DVM). Results for $\delta = 0$ $(\gamma = 3)$ and $\delta = 4$ $(\gamma = 1.4)$ shown in black and red, respectively. Diverged solutions not plotted. }
    \end{figure}
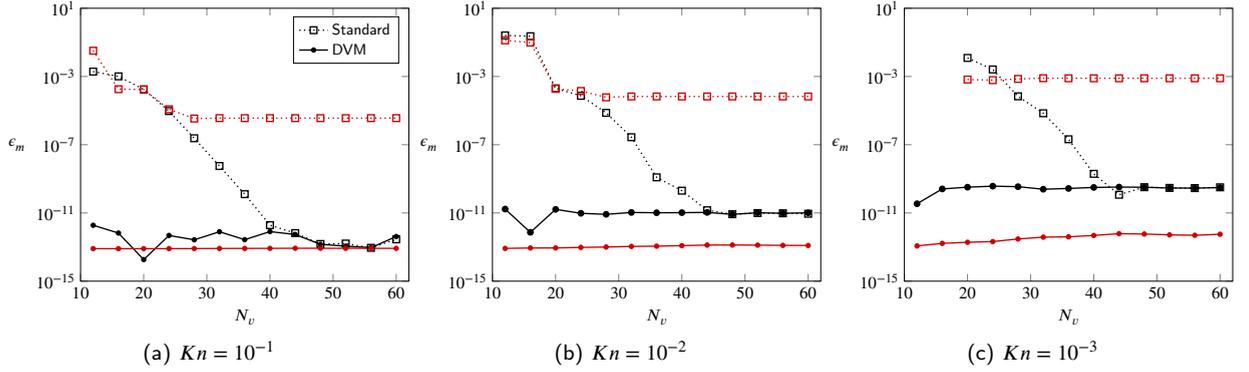

To verify the conservation properties of the DVM in comparison to the standard approach, the mass conservation error, defined as 
\begin{equation}
    \epsilon_{m} = \left | \int_{\Omega^x} \rho(x, t)\ \mathrm{d}x - \int_{\Omega^x} \rho\\(x, 0)\ \mathrm{d}x\right |,
\end{equation}
was similarly analyzed after one flow-through of the domain. The convergence of this mass conservation error with respect to the velocity/internal energy resolution is shown in \cref{fig:velocity_masscon} for $Kn = 10^{-1}, 10^{-2},$ and $10^{-3}$. In all cases, with just two iterations of the DVM, the mass conservation error was essentially negligible regardless of the resolution, varying from $10^{-13}$ to $10^{-10}$ across the various problem setups. However, without the DVM, the conservation error was significant, particularly at lower resolutions where errors of $\mathcal O(10^{-2})$ and larger were observed. For the monatomic case, where the domain bounds were chosen such that the values of the distribution function outside of the domain were essentially negligible ($\epsilon_{\mathbf{u}}$), the conservation error without the DVM did eventually converge to the DVM results with increasing resolution, but the resolution required to match the DVM results was significantly larger, such that the velocity/internal energy discretization required a factor of $3{-}4$ times as many nodes \textit{per dimension}. In contrast, for the polyatomic case, the bounds of the internal energy domain were chosen such that value of the distribution function outside of it was small but non-negligible ($\epsilon_{\zeta}$) to decrease the resolution requirements. In this case, the mass conservation error with the DVM was still negligible, but without the DVM, the error was significant and did \textit{not} converge with increasing resolution. This effect in the standard approach can be attributed to the error introduced by truncating the infinite domain. As a result, the DVM offers even more advantages in the overall computational cost as it allows for the use of smaller domain sizes which can be better resolved for the same number of discrete nodes.  

\subsection{Double Expansion Wave}
A double expansion wave was subsequently considered as a verification of the entropy-satisfying and positivity-preserving properties of the scheme. This problem, known as the \textit{123 problem} (or Test 2) in \citet{Toro1997}, is solved on the domain $\Omega^x = [0,1]$ with the initial conditions
\begin{equation}\label{eq:expansion}
    \mathbf{q}(x, 0) = 
    \begin{cases}
    \mathbf{q}_L, \quad \mathrm{if } x \leq 0.5,\\
    \mathbf{q}_R, \quad \mathrm{else},
    \end{cases}
    \quad \mathrm{where} \quad 
    \mathbf{q}_L = 
    \begin{bmatrix}
           1 \\
           -2 \\
           0.4
    \end{bmatrix},
    \quad
    \mathbf{q}_R = 
    \begin{bmatrix}
           1 \\
           2 \\
           0.4
    \end{bmatrix}.
\end{equation}
The initial discontinuous state develops into two outrunning expansion waves with a low density/pressure region in the center. The resulting symmetric wave structure is difficult to resolve, particularly in the center region, with even robust first-order numerical schemes for gas dynamics frequently failing \citep{Toro1997}. 

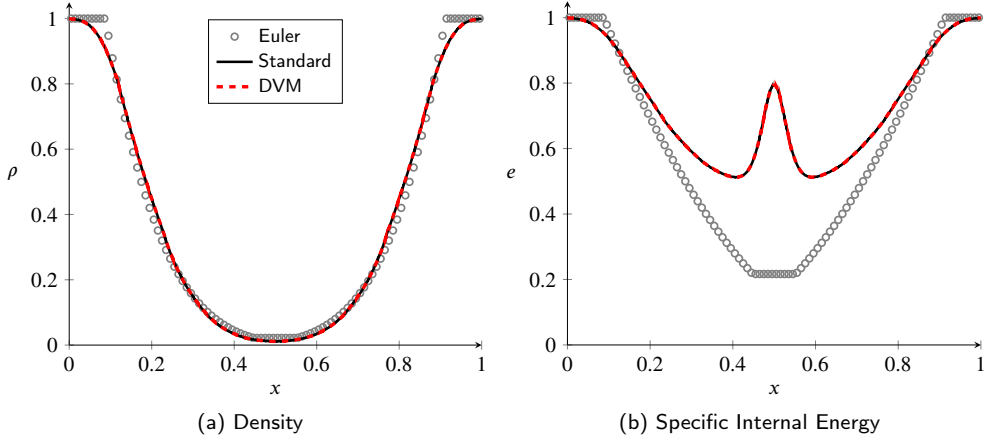
\begin{figure}[tbhp]
        \subfloat[Density]{\adjustbox{width=0.4\linewidth, valign=b}{\input{figs/expansion_density_Kn0p01}}}
        \subfloat[Specific Internal Energy]{\adjustbox{width=0.4\linewidth, valign=b}{\input{figs/expansion_energy_Kn0p01}}}
        \newline
        \caption{\label{fig:expansion_kn0p01} Density (left) and specific internal energy (right) profile for the double expansion wave problem at $t = 0.15$ with $Kn = 10^{-2}$ computed using a $\mathbb P_3$ scheme with 100 elements, $N_v = N_\zeta = 32$, and $\delta = 4$ ($\gamma = 1.4$). Results shown in comparison to the exact Euler profiles. }
    \end{figure}

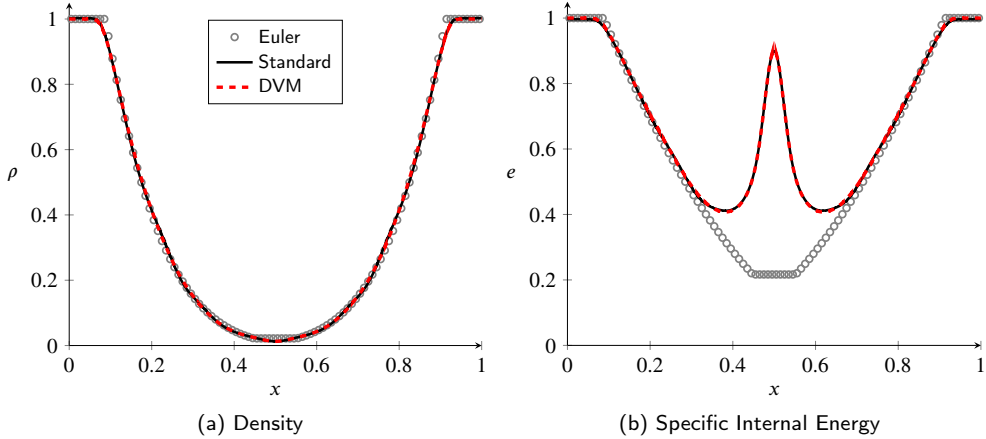
\begin{figure}[tbhp]
        \subfloat[Density]{\adjustbox{width=0.4\linewidth, valign=b}{\input{figs/expansion_density_Kn0p001}}}
        \subfloat[Specific Internal Energy]{\adjustbox{width=0.4\linewidth, valign=b}{\input{figs/expansion_energy_Kn0p001}}}
        \newline
        \caption{\label{fig:expansion_kn0p001}  Density (left) and specific internal energy (right) profile for the double expansion wave problem at $t = 0.15$ with $Kn = 10^{-3}$ computed using a $\mathbb P_3$ scheme with 100 elements, $N_v = N_\zeta = 32$, and $\delta = 4$ ($\gamma = 1.4$). Results shown in comparison to the exact Euler profiles. }
    \end{figure}

The problem was simulated using a $\mathbb P_3$ scheme with $N_e = 100$ elements and $N_v = N_\zeta = 32$ with Dirichlet boundary conditions. For consistency with the original problem setup, the polyatomic case was chosen, with $\delta = 4$ such that $\gamma = 1.4$. The predicted density and specific internal energy $e = \theta/(\gamma -1)$ at $t = 0.15$ is shown in \cref{fig:expansion_kn0p01} and \cref{fig:expansion_kn0p001} for $Kn = 10^{-2}$ and $Kn = 10^{-3}$, respectively, based on the initial macroscopic state. For both Knudsen numbers, the density profile was well-resolved even in the low-density center region, with convergence to the exact Euler results qualitatively observed with decreasing $Kn$ and near identical behavior with and without the DVM. For the specific internal energy profile, which is generally much more difficult to accurately resolve in this region as errors are amplified by the low density and pressure, the predicted results showed good agreement in the decreasing $Kn$ limit on the exterior of the expansion waves. However, the center region shows a spike in the internal energy profile, indicative of spurious physical entropy generation in the center. This effect is quite common for gas dynamics solvers, with almost all schemes producing this energy spike (see \citet{Toro1997}, Section 6.4). Interestingly enough, we also observe this effect when solving the Boltzmann equation, which should, in theory, be entropy-satisfying. 

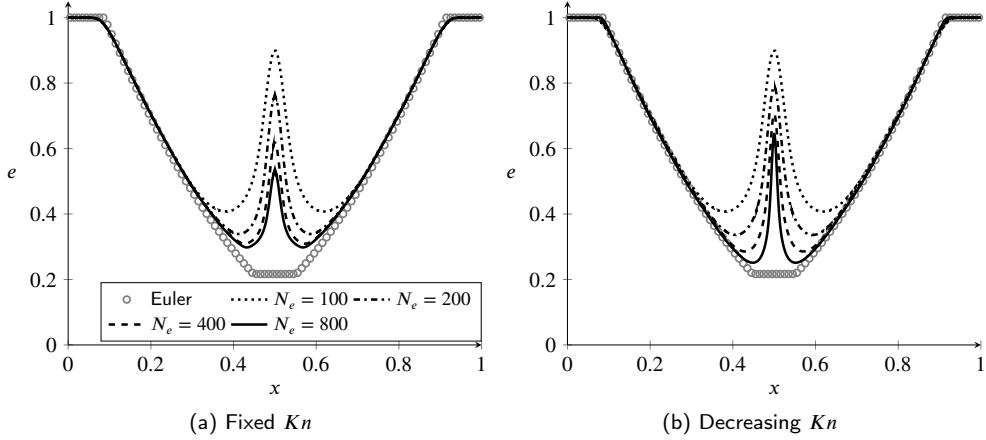
\begin{figure}[tbhp]
        \subfloat[Fixed $Kn$]{\adjustbox{width=0.4\linewidth, valign=b}{\input{figs/expansion_energy_conv_fixed}}}
        \subfloat[Decreasing $Kn$]{\adjustbox{width=0.4\linewidth, valign=b}{\input{figs/expansion_energy_conv_decreasing}}}
        \newline
        \caption{\label{fig:expansion_conv_mesh}  
        Convergence of the specific internal energy profile with respect to mesh resolution $N_e$ for the double expansion wave problem at $t = 0.15$ computed using a $\mathbb P_3$ scheme, $N_v = N_\zeta = 32$ with the discrete velocity model, and $\delta = 4$ ($\gamma = 1.4$). Left: Discontinuous initial conditions ($\Delta = 0$). Right: Continuous initial conditions ($\Delta \sim h$). Results shown in comparison to the exact Euler profiles.
        }
    \end{figure}

To discern the cause of this spurious entropy generation, we explore several possibilities. One may attribute this problem to the effects of a non-zero Knudsen number, but \cref{fig:expansion_kn0p01} and \cref{fig:expansion_kn0p001} show that although the profiles in the expansion region tend toward the Euler results with decreasing $Kn$, the magnitude of the spike actually \textit{increases} with decreasing $Kn$. Another explanation may be that this is an artifact of the mesh resolution, which is generally the case for gas dynamics solvers but would certainly raise questions about the entropy-satisfying properties of the proposed scheme. To explore this possibility, the problem was solved on a progression of meshes with increasing resolution, both with a fixed $Kn = 10^{-3}$ and with a proportionally decreasing $Kn = h/10$, where $h = 1/N_e$ is the mesh size. The predicted specific internal energy profiles for these cases are shown in \cref{fig:expansion_conv_mesh}. It can be seen that for the decreasing $Kn$ case (i.e., constant mesh Knudsen number), the profiles converge to the Euler results, but the spike does not proportionally decrease. In contrast, the fixed $Kn$ case (e.g., decreasing mesh Knudsen number) shows a significant decrease in the energy spike with increasing resolution, indicating that the problem is more associated with the mesh Knudsen number than the resolution itself. 

\begin{figure}[tbhp]
        \subfloat[Discontinuous ($\Delta = 0$)]{\adjustbox{width=0.4\linewidth, valign=b}{\input{figs/expansion_energy_conv_sharp}}}
        \subfloat[Continuous ($\Delta \sim h$)]{\adjustbox{width=0.4\linewidth, valign=b}{\input{figs/expansion_energy_conv_smooth}}}
        \newline
        \caption{\label{fig:expansion_conv_smooth}  
        Convergence of the specific internal energy profile with respect to mesh resolution $N_e$ and $Kn = h/10$ for the double expansion wave problem at $t = 0.15$ computed using a $\mathbb P_3$ scheme, $N_v = N_\zeta = 32$ with the discrete velocity model, and $\delta = 4$ ($\gamma = 1.4$). Left: Discontinuous initial conditions ($\Delta = 0$). Right: Continuous initial conditions ($\Delta \sim h$). Results shown in comparison to the exact Euler profiles.
        }
    \end{figure}
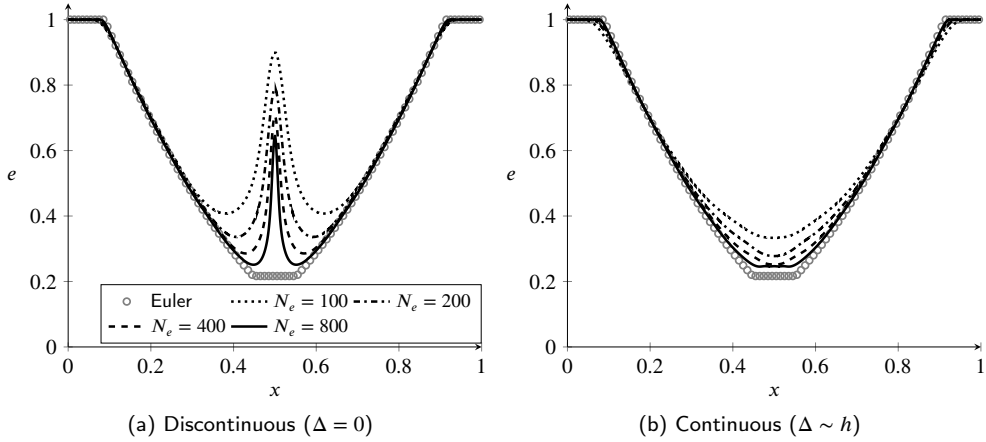

We posit that this effect is related to a violation of the resolution requirements proposed in \cref{ssec:disc}. Although it may seem that for the given parameters, the resolution condition in \cref{eq:resreq} is easily satisfied, the assumptions posed in \cref{ssec:disc} are not valid for the initial conditions --- regardless of the Knudsen number of the problem, the initial discontinuity thickness given by \cref{eq:expansion} is exactly zero. Therefore, the actual mesh Knudsen number for the initial state is zero, such that the scales are not properly resolved and \cref{eq:resreq} is not actually satisfied. This issue is exacerbated by the stationary nature of the center region in symmetric Riemann problems, which does not benefit from ``correcting'' effects such as compression in shocks or numerical dissipation in traveling contact discontinuities. To verify this claim, the previous mesh convergence study was run again utilizing the case of a decreasing Knudsen ($Kn = h/10$), but the discontinuous original initial condition was modified to a \textit{continuous} initial condition with a finite discontinuity thickness $\Delta$. The modified initial velocity profile was chosen as
\begin{equation}
    U(x, 0) = \frac{U_R - U_L}{2} \tanh\left[\frac{1}{h} (x - 0.5) \right] + \frac{U_R + U_L}{2} = 2\tanh\left[\frac{1}{h} (x - 0.5) \right],
\end{equation}
such that the discontinuity thickness $\Delta$ is of $\mathcal O(h)$ which satisfies the resolution requirements at $t = 0$ and recovers the Euler equations in the $Kn\to 0$, $h \to 0$ limit. A comparison of the predicted specific internal energy profile for the discontinuous and continuous initial conditions is shown in \cref{fig:expansion_conv_smooth}. As it can clearly be seen, modifying the initial conditions to be consistent with the system in question completely removes the spurious spike in specific internal energy (and consequently, physical entropy) while converging towards the Euler results in the $Kn\to 0$, $h \to 0$ limit. Note that the modification was only necessary due to a pathological test case known to be sensitive to this issue --- for practical purposes (and for the remaining numerical experiments in this work), the modification to the initial conditions was neglected as its effect was generally negligible in comparison to numerical approximation errors.
    
\subsection{Normal Shock Structure}
As a verification of the proposed method for more complex flow regimes and a comparison of the Boltzmann-BGK method to continuum approaches, the structure of a normal shock in argon was computed over a range of Mach numbers. For this problem, the domain is set as $\Omega^x = [-25, 25]$ and the initial conditions are given as 
\begin{equation}
    \mathbf{q}(x, 0) = 
    \begin{cases}
    \mathbf{q}_L, \quad \mathrm{if } x \leq 0,\\
    \mathbf{q}_R, \quad \mathrm{else},
    \end{cases}
    \quad \mathrm{where} \quad 
    \mathbf{q}_L = 
    \begin{bmatrix}
           1 \\
           M \sqrt{\gamma} \\
           1
    \end{bmatrix}.
\end{equation}
The left-hand state $\mathbf{q}_L$ is given in terms of the incoming (upstream) Mach number $M$, and the right-hand state $\mathbf{q}_R$ can be computed through the Rankine--Hugoniot conditions as 
\begin{equation}
    \frac{\rho_R}{\rho_L} = \frac{(\gamma + 1)M^2}{(\gamma - 1)M^2 + 2}, \quad \quad
    \frac{U_R}{U_L} = \frac{(\gamma - 1)M^2 + 2}{(\gamma + 1)M^2}, \quad \quad \mathrm{and} \quad \quad 
    \frac{P_R}{P_L} = \frac{2\gamma M^2 - (\gamma - 1)}{\gamma +1}.
\end{equation}
For comparison with experimental data and other numerical approaches, the specific heat ratio $\gamma$ was set as $5/3$ (i.e., $\delta = 2$) to mimic a monatomic molecule (argon) in three dimensions. The quantity of interest for this problem is the inverse ratio of the shock thickness $\Delta$, defined as 
\begin{equation}
    \Delta = \frac{\rho_R - \rho_L}{\max \left[\partial_x \rho (x) \right]},
\end{equation}
relative to the upstream mean free path $\lambda_L$.

Since the shock structure is highly sensitive to the transport coefficients, the use of a constant collision time $\tau$, for which the viscosity would depend on pressure, can result in inaccurate predictions. One can instead adapt the collision time based on the solution to recover more physically consistent transport coefficients similarly to the approach of \citet{Mieussens2000} as 
\begin{equation}
    \tau(\mathbf{Q}) = \tau_{ref} \frac{\rho_{ref} \theta_{ref}^{1 - \omega}}{\rho \theta^{1 - \omega}}.
\end{equation}
Here, $\omega$ refers to the viscosity law exponent of the gas, taken as $\omega = 0.81$ for argon, and the reference quantities are taken from the incoming macroscopic state $\mathbf{Q}_L$. The reference collision time $\tau_{ref}$ was computed by setting the Knudsen number as $Kn = 1$ based on a unit reference length and the incoming macroscopic state $\mathbf{Q}_L$, such that $\lambda_L = Kn$.

\begin{figure}[tbhp]
        \tikzexternaldisable
        \adjustbox{width=0.48\linewidth, valign=b}{\input{figs/normalshock_itr}}
        \newline
        \caption{\label{fig:normalshock_itr} Inverse thickness ratio with respect to inflow Mach number for a stationary normal shock in argon ($\delta = 2$, $\gamma = 5/3$) computed using a $\mathbb P_3$ scheme with 100 elements and $N_v = N_\zeta = 32$. Navier--Stokes results of \citet{Mieussens2000} and experimental results of \citet{Alsmeyer1976} (squares), \citet{Linzer1963} (triangles), and \citet{Camac1964} (diamonds) shown for reference.}
        \tikzexternalenable
    \end{figure}
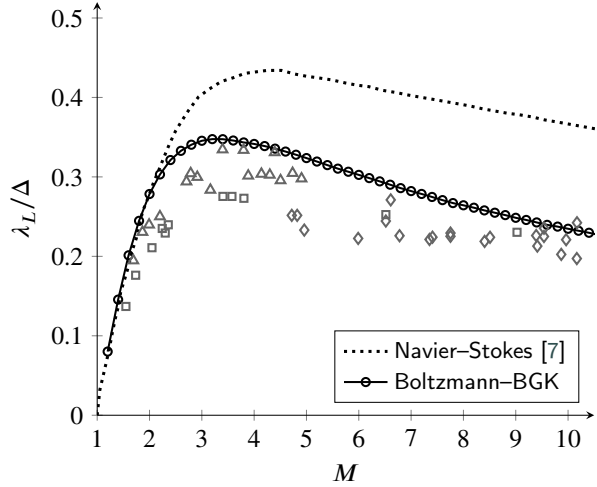

The shock structure was computed using a $\mathbb P_3$ scheme with 100 elements with $N_v = N_\zeta = 32$ which allows for the structure to be fully resolved based on the Knudsen number and the spatial resolution. The boundary conditions were set to Dirichlet on both ends. To compute the shock thickness, the corrected gradients were used, with the common interface solution taken as the mean of the interior and exterior values similarly to the BR1 approach of \citet{Bassi1997}. The solution was advanced to a final time of $t=100$, after which the $L^\infty$ norm of the temporal residual was generally converged to values less than $10^{-5}$. The computed inverse thickness ratio $\lambda_L/\Delta$ is shown in \cref{fig:normalshock_itr} over a range of Mach numbers in comparison to the Navier--Stokes results of \citet{Mieussens2000} and experimental results of \citet{Alsmeyer1976}, \citet{Linzer1963}, and \citet{Camac1964}. The predictions of the present Boltzmann--BGK approach were in much better agreement with the experimental data than the Navier--Stokes results, properly predicting both the critical Mach number where this inverse thickness ratio is at its maxima and the decay of the ratio with increasing Mach number. The predicted inverse thickness ratio was generally on the higher end of the experimental data but still within the range of uncertainty between the various experiments. 

\begin{figure}[tbhp]
        \subfloat[$M = 3.8$]{\adjustbox{width=0.48\linewidth, valign=b}{\input{figs/normalshock_profile_m3p8}}}
        \subfloat[$M = 9.0$]{\adjustbox{width=0.48\linewidth, valign=b}{\input{figs/normalshock_profile_m9p0}}}
        \newline
        \caption{\label{fig:normalshock_profiles} Normalized density, velocity, and pressure profiles for a stationary normal shock in argon ($\delta = 2$, $\gamma = 5/3$) at $M = 3.8$ (left) and $M = 9.0$ (right) computed using a $\mathbb P_3$ scheme with 100 elements and $N_v = N_\zeta = 32$. Experimental density profiles of \citet{Alsmeyer1976} (square markers) and DSMC density profiles of \citet{Bentley2009} (triangle markers) shown for reference. }
    \end{figure}
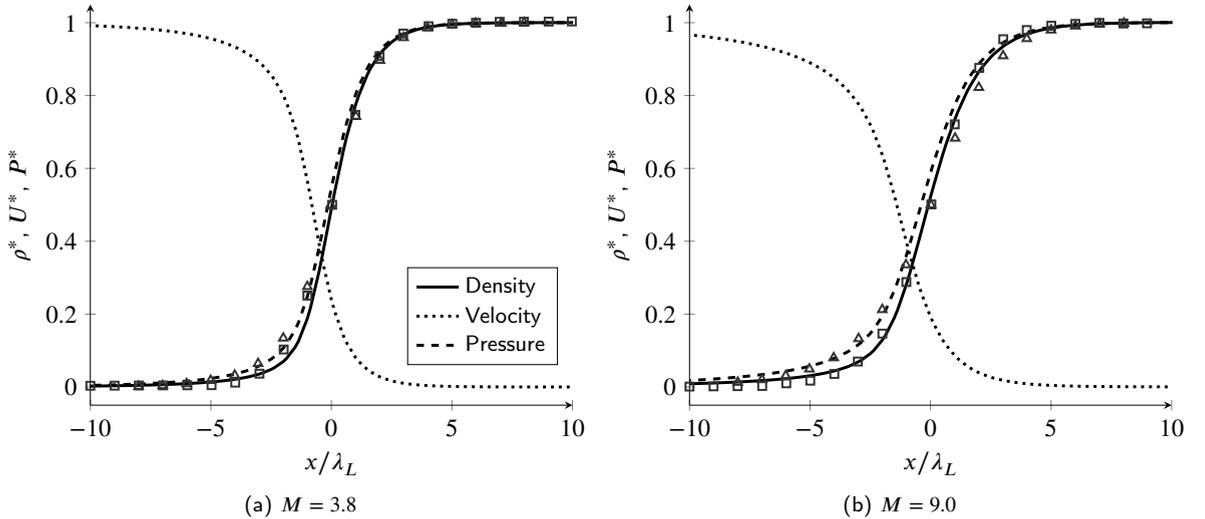

For a more direct comparison of the shock structures, the predicted normalized density profiles ($\rho^*$) as well as normalized velocity ($U^*$) and pressure ($P^*$) profiles at two Mach numbers, $M = 3.8$ and $M = 9.0$, were compared to the experimental density profiles of \citet{Alsmeyer1976} and the direct simulation Monte Carlo (DSMC) density profiles of \citet{Bentley2009} (computed using the method of \citet{Bird1994}) in \cref{fig:normalshock_profiles}. At $M = 3.8$, the predicted density profile was generally in good agreement with the experimental and DSMC results with respect to both the thickness and the structure, although some minor undershoot in the upstream section of the shock was observed. At $M = 9.0$, even better agreement with the experiment was observed, with the predicted density profiles matching up almost exactly and with notably better predictions than the DSMC results. The benefits of the kinetic approach to predicting shock structures are clearly evident in these results. 
    
To observe the effects of the various flow conditions on the distribution function, the profiles of the distribution function at various locations in the flow were visualized for $M = 3.8$ and $M = 9.0$. Since the velocity/internal energy domain is two-dimensional in this case, the velocity component $f_u(u)$ was extracted by taking the discrete maxima over the internal energy domain, i.e.,
\begin{equation}
    f_u(u) = \underset{\zeta}{\mathrm{max}}\ f(u, \zeta).
\end{equation}

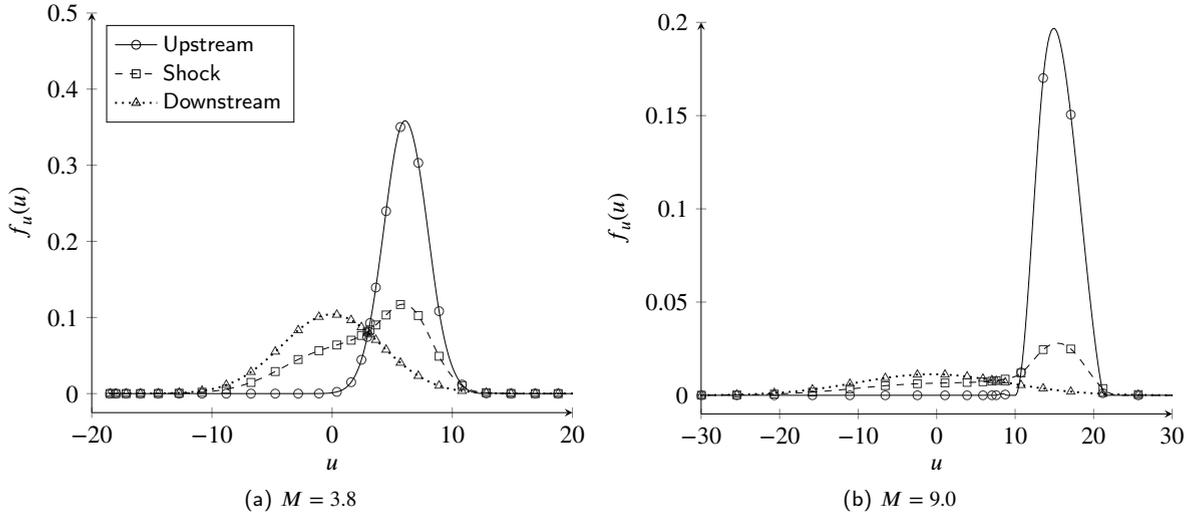
\begin{figure}[tbhp]
        \subfloat[$M = 3.8$]{\adjustbox{width=0.48\linewidth, valign=b}{\input{figs/normalshock_pdf_m3p8}}}
        \subfloat[$M = 9.0$]{\adjustbox{width=0.48\linewidth, valign=b}{\input{figs/normalshock_pdf_m9p0}}}
        \newline
        \caption{\label{fig:normalshock_pdf}
        Velocity component of the distribution function at various spatial locations for a stationary normal shock in argon ($\delta = 2$, $\gamma = 5/3$) at $M = 3.8$ (left) and $M = 9.0$ (right) computed using a $\mathbb P_3$ scheme with 100 elements and $N_v = N_\zeta = 32$. Nodal values shown as markers, spline interpolation shown as lines.         }
    \end{figure}

The predicted distribution function profiles for $M = 3.8$ and $M = 9.0$ are shown in \cref{fig:normalshock_pdf} at the upstream ($x = -25$), shock ($x = 0$), and downstream ($x = 25$) locations. The discrete values of the distribution function are shown as markers, but for presentation purposes, a spline interpolation of these values is also shown. The profiles indicate that the \textit{a priori} approach for choosing the extent of the velocity domain is sufficient as all of the statistically significant variation in the profiles was completely captured within the domain. For both Mach numbers, the upstream and downstream distribution functions appeared Maxwellian, indicating the flow is in thermodynamic equilibrium, whereas at the shock showed clear non-equilibrium behavior. The distribution functions both also showed the increase of the thermal velocity across the shock. These results appear to agree qualitatively with the results of \citet{Mieussens2000} obtained with a three-dimensional velocity space at a slightly different Mach number. 

\subsection{Sod Shock Tube}
Extensions to unsteady gas dynamics with discontinuities was performed through the Sod shock tube \citep{Sod1978}. This one-dimensional test case exhibits the three main features of the Riemann problem, namely a shock wave, a contact discontinuity, and a rarefaction wave, which makes it a suitable test for the discontinuity-resolving properties of the proposed scheme. The problem is computed on the domain $\Omega^x = [0, 1]$ and the initial conditions are given as
\begin{equation}
    \mathbf{q}(x, 0) = 
    \begin{cases}
    \mathbf{q}_L, \quad \mathrm{if } x \leq 0.5,\\
    \mathbf{q}_R, \quad \mathrm{else},
    \end{cases}
    \quad \mathrm{where} \quad 
    \mathbf{q}_L = 
    \begin{bmatrix}
           1 \\
           0 \\
           1
    \end{bmatrix},
    \quad
    \mathbf{q}_R = 
    \begin{bmatrix}
           0.125 \\
           0 \\
           0.1
    \end{bmatrix},
\end{equation}
with Dirichlet boundary conditions on both ends. The specific heat ratio $\gamma$ was set to 1.4 ($\delta = 4$) for consistency with the original problem setup.

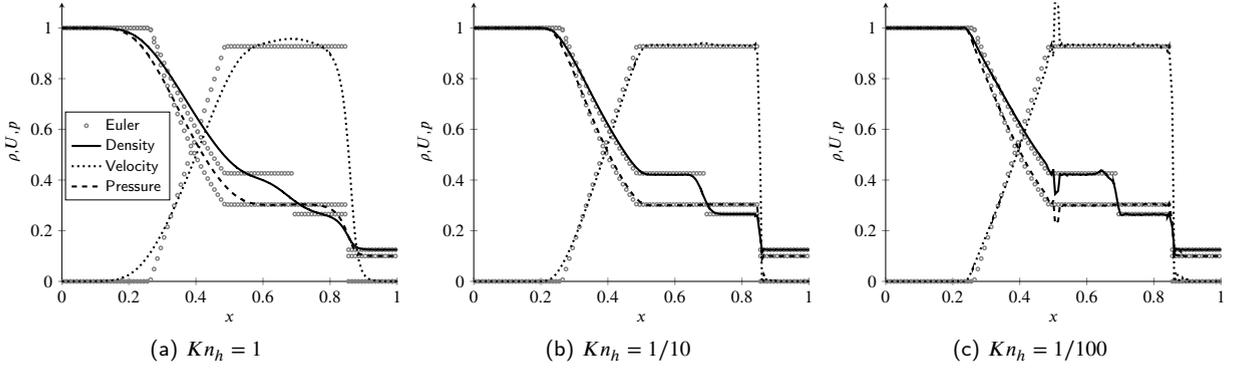
\begin{figure}[tbhp]
        \subfloat[$Kn_h = 1$]{\adjustbox{width=0.33\linewidth, valign=b}{\input{figs/sod_density_kn1}}}
        \subfloat[$Kn_h = 1/10$]{\adjustbox{width=0.33\linewidth, valign=b}{\input{figs/sod_density_kn0p1}}}
        \subfloat[$Kn_h = 1/100$]{\adjustbox{width=0.33\linewidth, valign=b}{\input{figs/sod_density_kn0p01}}}
        \newline
        \caption{\label{fig:sod_kn} Density, velocity, and pressure profiles for the Sod shock tube problem at $t = 0.2$ with varying mesh Knudsen number $Kn_h$ computed using a $\mathbb P_3$ scheme with 50 elements, $N_v = N_\zeta = 16$, and $\delta = 4$ ($\gamma = 1.4$). Results shown in comparison to the exact Euler profiles. }
    \end{figure}

To evaluate the validity of the resolution requirements posed in \cref{ssec:disc}, the parameters of the scheme were fixed while modulating the mesh Knudsen number. A $\mathbb P_3$ scheme with 50 elements was used with $N_v = N_\zeta = 16$ and mesh Knudsen numbers of $Kn_h = 1$, $1/10$, and $1/100$ based on the initial macroscopic state. The predicted macroscopic solution profiles are shown in \cref{fig:sod_kn} at $t = 0.2$ for these respective mesh Knudsen numbers. For $Kn_h = 1$, the solution was excessively diffused, although the locations and magnitudes of most discontinuities were reasonably predicted. With $Kn_h = 1/10$ (the proposed resolution limit), the solution was in good agreement with the exact Euler results, with excellent resolution of the location and magnitude of the discontinuities and minimal observable spurious oscillations. The contact discontinuity was more diffused compared to the shock, but this behavior is expected due to the compressive nature of shock waves. When the mesh Knudsen number was lowered past the proposed limit to $Kn_h = 1/100$, spurious oscillations began to appear in the predicted profiles, severely degrading the accuracy of the solution. However, the location and magnitude of the discontinuities were still reasonably well predicted, indicating that the prediction of the governing physics is not severely degraded even with a numerically ill-behaved solution. 

The proposed resolution limit was further explored with a mesh convergence study. The mesh Knudsen number was fixed at $Kn_h = 1/10$ while the mesh resolution was varied with $N_e = 50$, $100$, and $200$. The predicted solution profiles showed excellent convergence to the exact Euler results with increasing resolution (i.e., decreasing Knudsen number). Sub-element resolution of the shock was generally observed, whereas the contact was generally resolved over 1-2 elements. As mentioned in \cref{ssec:disc}, the resolution limit does not guarantee monotonicity of the solution in the vicinity of a discontinuity. This was observed in the predicted velocity profiles with very minor overshoots in the vicinity of the shock. However, the predicted density and pressure profiles were (visually) monotonic in the vicinity of the discontinuities. These results indicate that the mesh Knudsen number offers a suitable indicator metric for the resolution of discontinuous structures. 
   \begin{figure}[tbhp]
        \subfloat[$N_e = 50$]{\adjustbox{width=0.33\linewidth, valign=b}{\input{figs/sod_density_50}}}
        \subfloat[$N_e = 100$]{\adjustbox{width=0.33\linewidth, valign=b}{\input{figs/sod_density_100}}}
        \subfloat[$N_e = 200$]{\adjustbox{width=0.33\linewidth, valign=b}{\input{figs/sod_density_200}}}
        \newline
        \caption{\label{fig:sod_ne} Density, velocity, and pressure profiles for the Sod shock tube problem at $t = 0.2$ with a fixed mesh Knudsen number $Kn_h = 1/10$ computed using a $\mathbb P_3$ scheme with varying mesh resolution, $N_v = N_\zeta = 16$, and $\delta = 4$ ($\gamma = 1.4$). Results shown in comparison to the exact Euler profiles.}
    \end{figure}
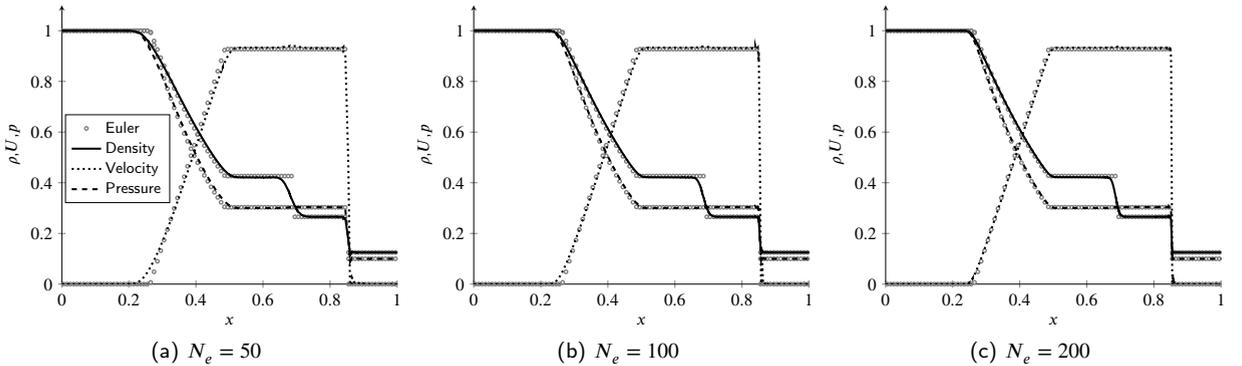

\subsection{Mach 3 Forward Facing Step}
The proposed approach was then extended to large-scale problems on two-dimensional unstructured grids with complex flow physics including strong shocks and boundary interactions through the forward facing step problem of \citet{Woodward1984}. This problem originally consists of an inviscid Mach 3 flow in a wind tunnel with a step perturbation that forms strong shocks, rarefaction fans, and contact discontinuities. The domain is set as $\Omega^{\mathbf{x}} = [0, 3] \times [0, 1] \setminus [0.6, 3] \times [0, 0.2]$ with the initial conditions set as $\mathbf{q}(\mathbf{x}, 0) = [1, 1, 0, 1/(\gamma M^2)]^T$ for a specific heat ratio $\gamma = 1.4$ ($\delta = 3)$ and Mach number $M = 3$. Dirichlet boundary conditions were applied on the inlet ($x = 0$) and Neumann boundary conditions were applied on the outlet ($x = 3$). The remaining boundaries were set to specular wall boundary conditions. To vastly simplify the implementation of the specular wall boundary conditions, the velocity offsets were set as $\mathbf{U}_0 = \mathbf{0}$ and the polar resolution $N_{\phi}$ was set as a multiple of 4, such that the velocity space was symmetric about the wall normal directions and the reflection operator in \cref{eq:bcspecular} could be easily and efficiently implemented without requiring interpolation in velocity space. 

The problem was solved using a $\mathbb P_4$ scheme on a series of meshes of increasing resolution with $N_r = N_{\phi} = 16$ and $N_{\zeta} = 8$. The coarse, medium, and fine meshes were generated using approximately homogeneous unstructured triangles with characteristic mesh lengths of $h = 1/50$, $1/100$, and $1/200$, respectively, resulting in $2.3{\cdot}10^{4}$, $7.6{\cdot}10^{4}$, and $3.1{\cdot}10^{5}$ elements. For the given resolution in the spatial, velocity, and internal energy domains, this resulted in approximately 703 million total degrees of freedom for the coarse mesh, 2.38 billion total degrees of freedom for the medium mesh, and 9.65 billion total degrees of freedom for the fine mesh. The Knudsen number was set based on the initial conditions as $Kn = 1{\cdot}10^{-2}$, $5{\cdot} 10^{-3}$, and $2{\cdot}10^{-3}$ for the coarse, medium, and fine meshes, respectively, such that the mesh Knudsen number was approximately constant and the flow structures could be fully resolved. The approximate computational cost of the finest case was 480 GPU-hours (15 hours on 32 GPUs), which, given such a large number of degrees of freedom, shows the efficiency of the proposed approach.

    \begin{figure}[htbp!]
        \centering
        \subfloat[Coarse mesh, $Kn = 1{\cdot}10^{-2}$] {\adjustbox{width=0.8\linewidth,valign=b}{
            \includegraphics[width=\textwidth]{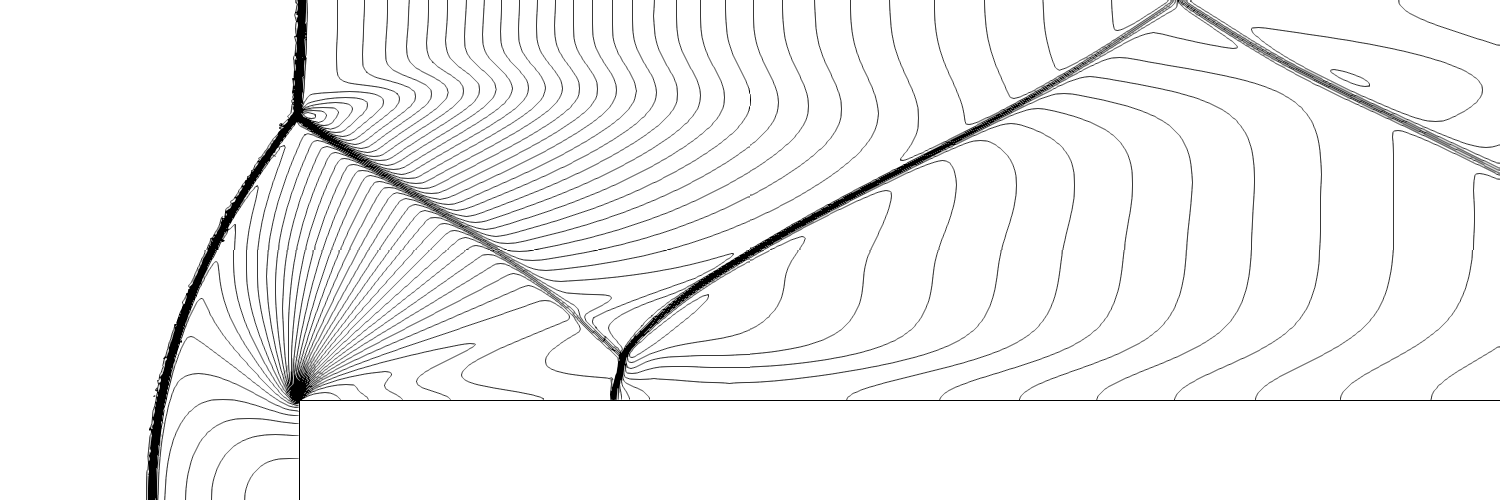}
        }}
        \newline
        \subfloat[Medium mesh, $Kn = 5{\cdot}10^{-3}$] {\adjustbox{width=0.8\linewidth,valign=b}{
            \includegraphics[width=\textwidth]{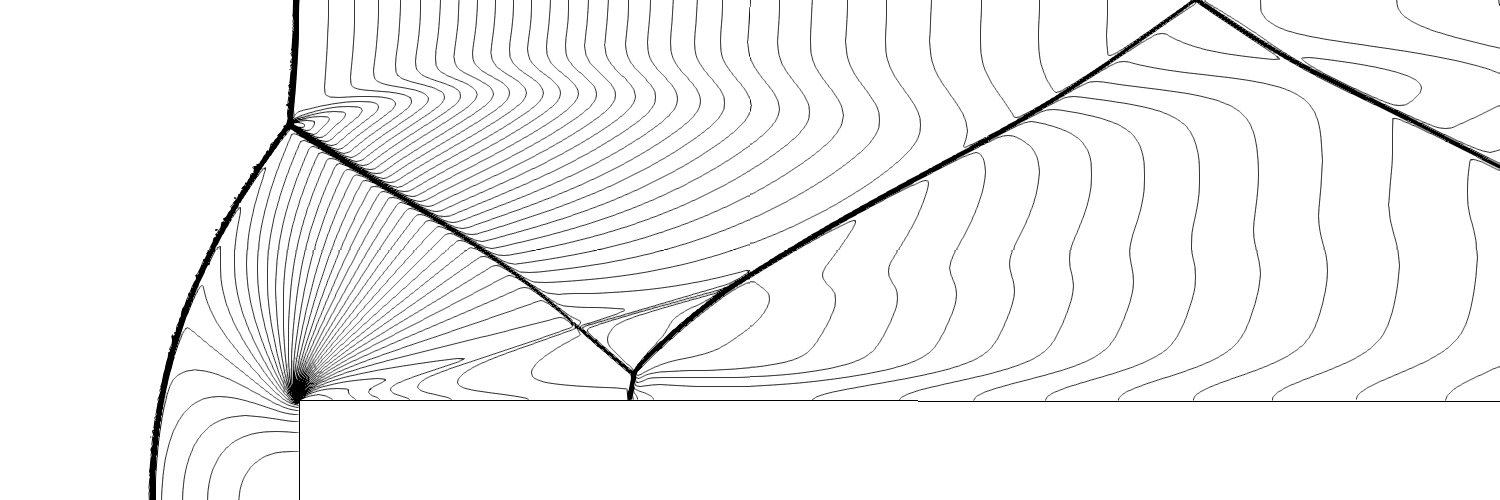}
        }}
        \newline
        \subfloat[Fine mesh, $Kn = 2{\cdot}10^{-3}$] {\adjustbox{width=0.8\linewidth,valign=b}{
            \includegraphics[width=\textwidth]{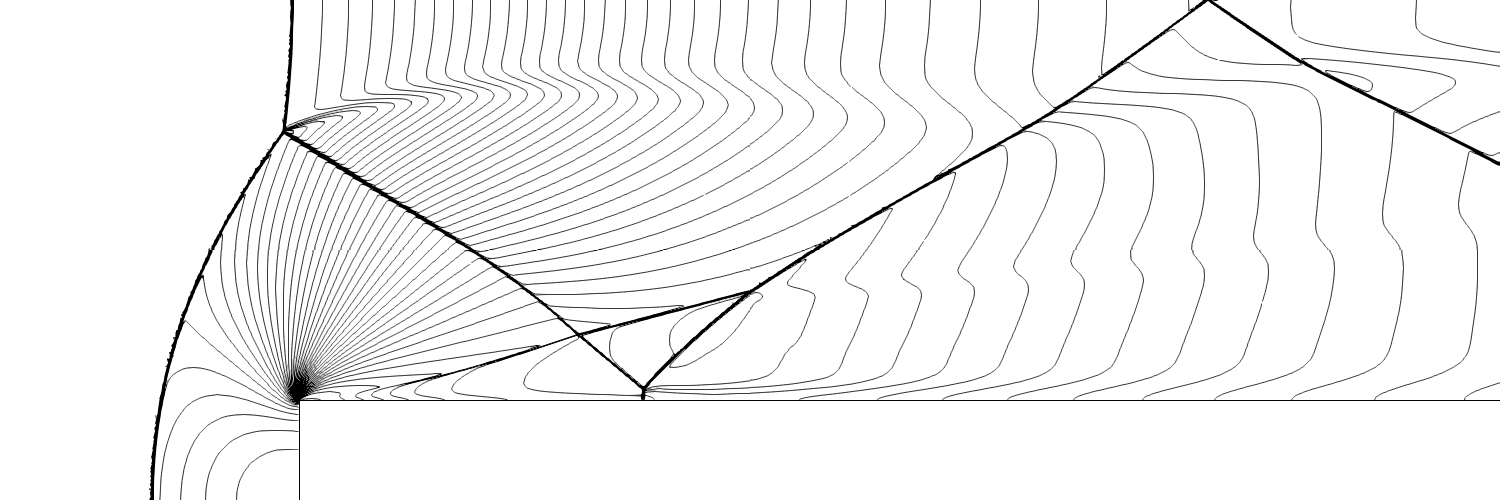}
        }}
        \newline
        \caption{\label{fig:step} Isocontours of density (equispaced on the range $[0,5]$) for the Mach 3 forward facing step problem at $t = 12$ computed using a $\mathbb P_4$ scheme with $N_r = N_{\phi} = 16$, $N_{\zeta} = 8$, and $\delta = 3$ ($\gamma = 1.4$). Top row: Coarse mesh ($2.3{\cdot}10^{4}$ elements) with $Kn = 1{\cdot}10^{-2}$. Middle row: Medium mesh ($7.6{\cdot}10^{4}$ elements) with $Kn = 5{\cdot}10^{-3}$. Bottom row: Fine mesh ($3.1{\cdot}10^{5}$ elements) with $Kn = 2{\cdot}10^{-3}$.}
    \end{figure}
    
The predicted density profiles at $t = 12$, represented through isocontours equispaced on the range $[0,5]$, are shown in \cref{fig:step} for the various meshes and problem setups. The results show the canonical structure of the forward-facing step, with a strong shock reflected from the step, a rarefaction fan centered on the corner, and a Mach stem emanating from the top wall. For each case, the discontinuities in the flow were very well resolved and showed no spurious oscillations in the vicinity of the shocks. Furthermore, the thickness of the discontinuities decreased accordingly with decreasing mesh size and Knudsen number. With decreasing Knudsen number, it can be seen that the results are converging to predictions in the Euler limit (see \citet{Dumbser2014}, Fig. 8). These effects were most noticeable in the ``liftoff'' of the Mach stem emanating from the bottom wall, which should decrease to zero in the absence of physical and numerical viscosity, and the behavior of the contact line emerging from the triple point. In the Euler limit, the onset of Kelvin--Helmholtz instabilities should become evident along this contact line, but with a finite Knudsen number, these instabilities can be stabilized due to viscous effects. 
    
    \begin{figure}[htbp!]
        \centering
        \subfloat[Coarse mesh, $Kn = 1{\cdot}10^{-2}$] {\adjustbox{width=0.3\linewidth,valign=b}{
            \includegraphics[width=\textwidth]{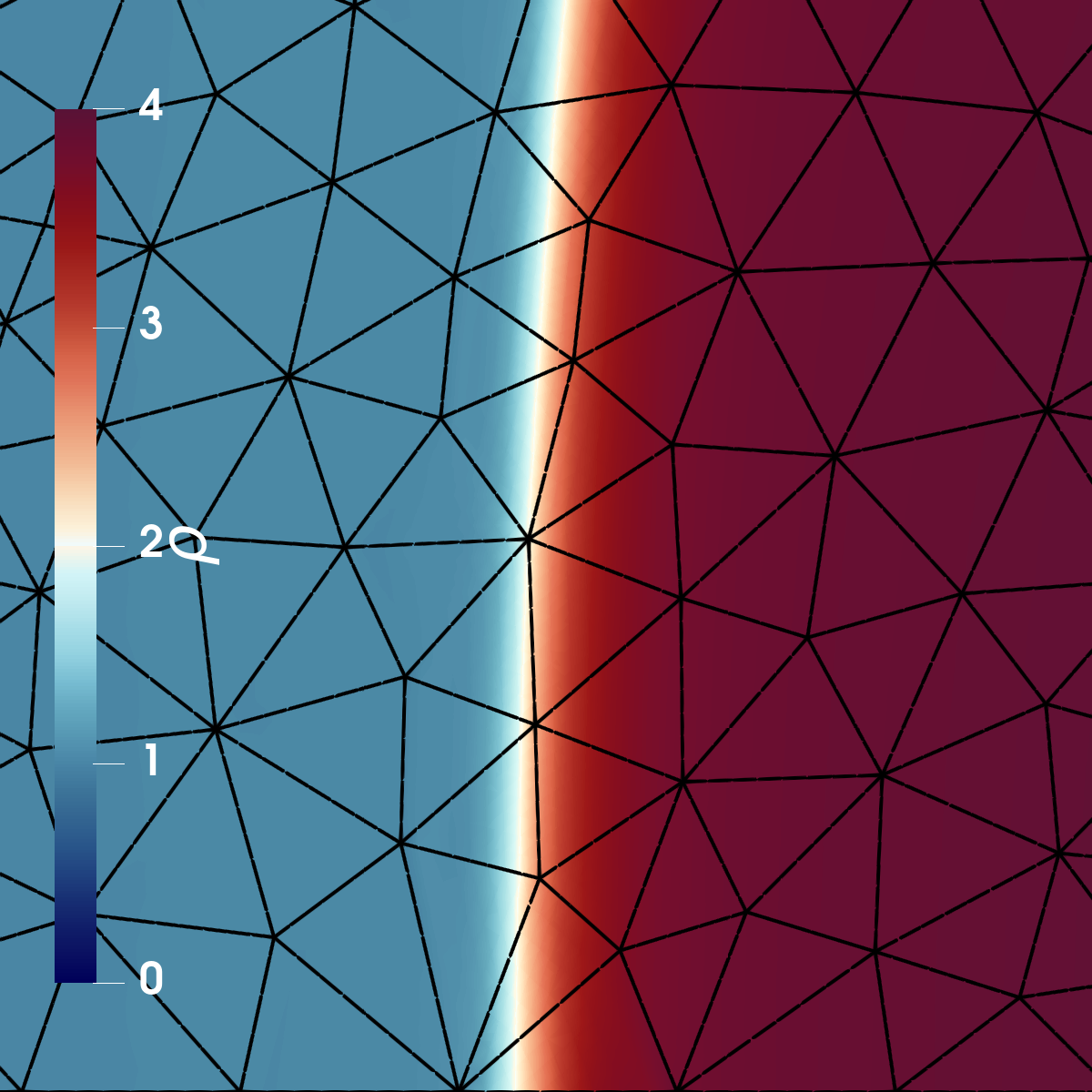}
        }}
        \subfloat[Medium mesh, $Kn = 5{\cdot}10^{-3}$] {\adjustbox{width=0.3\linewidth,valign=b}{
            \includegraphics[width=\textwidth]{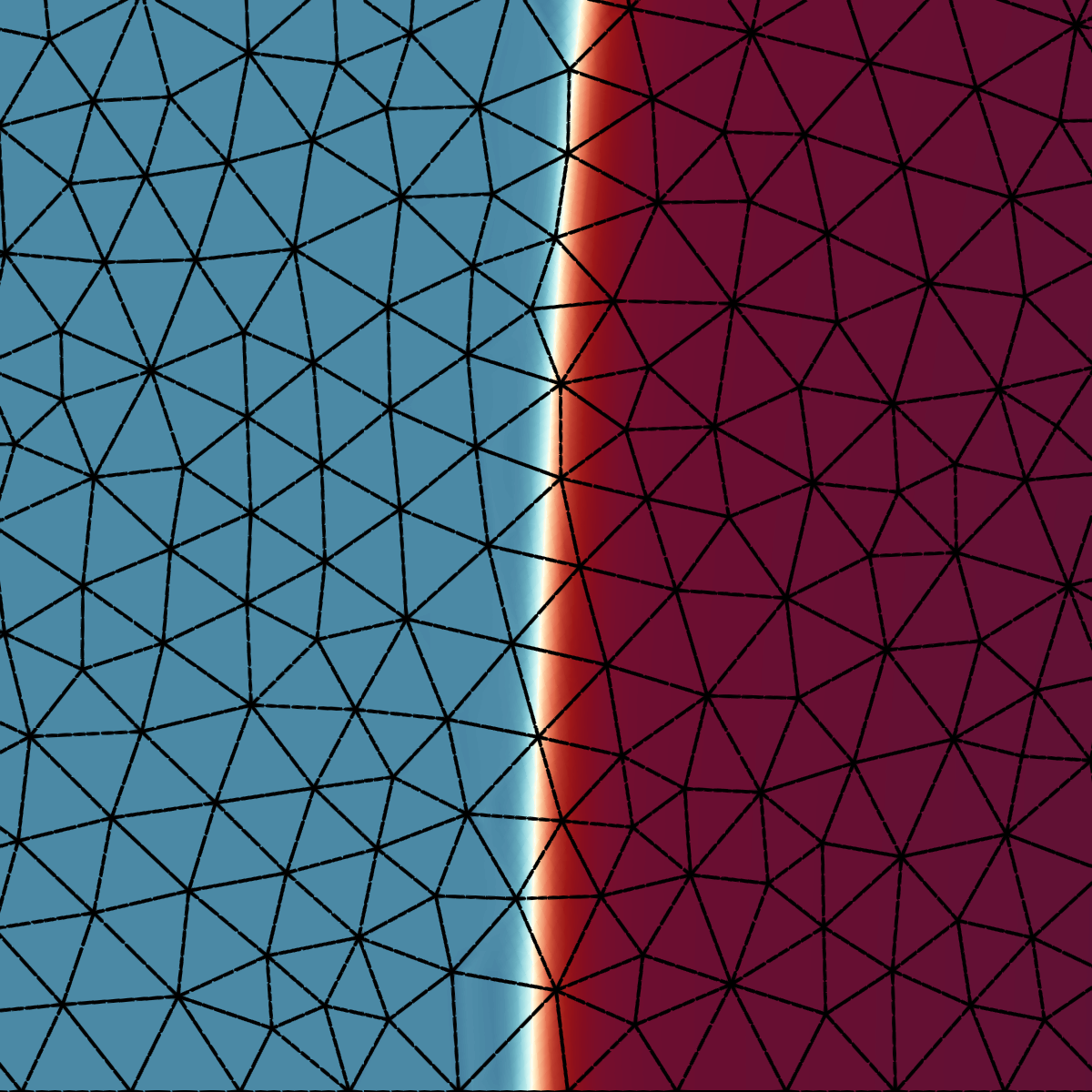}
        }}
        \subfloat[Fine mesh, $Kn = 2{\cdot}10^{-3}$] {\adjustbox{width=0.3\linewidth,valign=b}{
            \includegraphics[width=\textwidth]{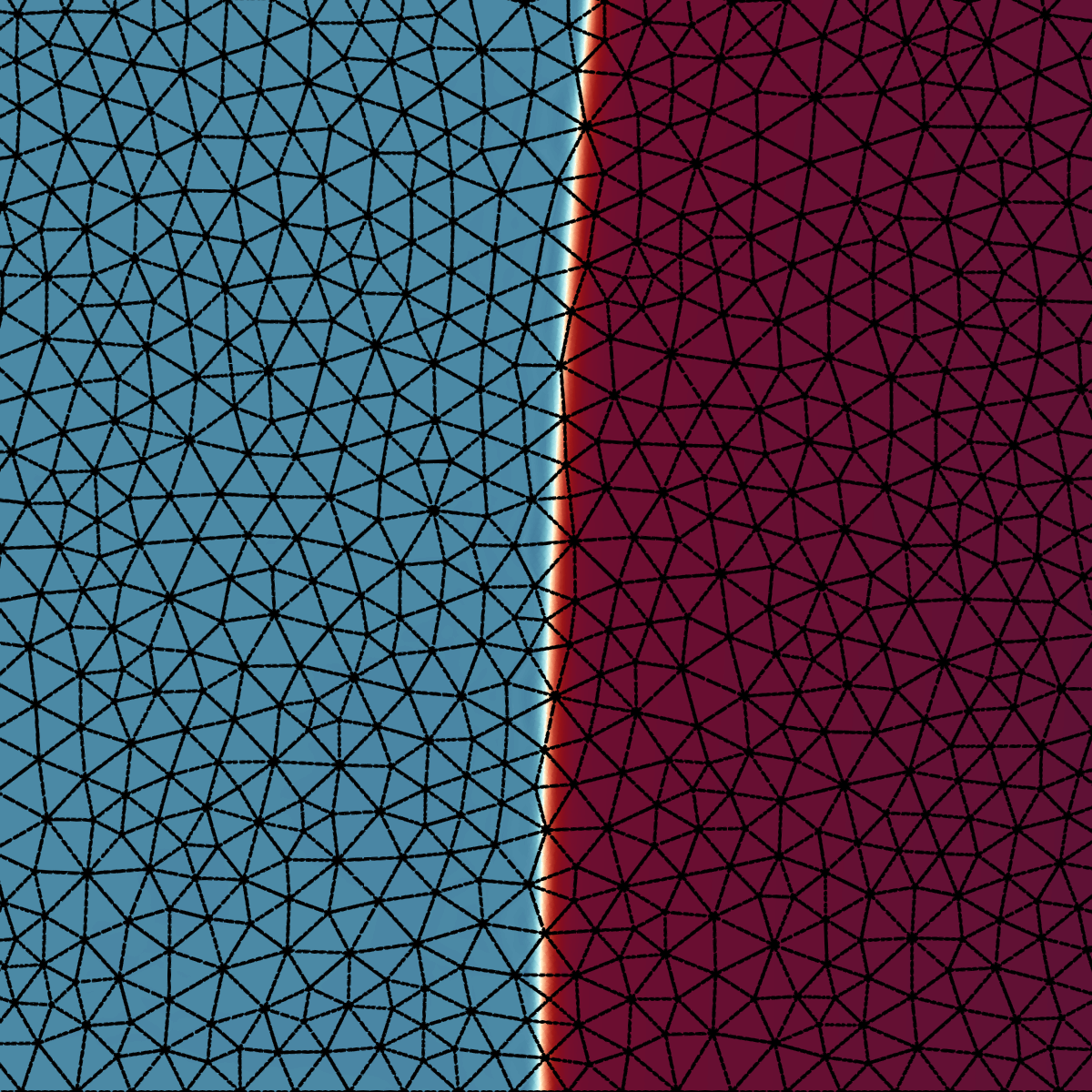}
        }}
        \newline
        \caption{\label{fig:stepmesh} 
        Enlarged view of the contours of density with mesh overlay near the leading shock for the Mach 3 forward facing step problem at $t = 12$ computed using a $\mathbb P_4$ scheme on a coarse (left), medium (middle), and fine (right) mesh.
        }
    \end{figure}

To evaluate the ability of the proposed approach to resolve shock structures, the contours of density in the region of the leading shock are shown in \cref{fig:stepmesh} in comparison to the computational mesh. It can be seen that the flow fields show quite impressive sub-element resolution of the shock, effectively predicting a monotonic shock structure over a length of only 1-2 spatial solution points. Furthermore, this behavior was maintained with decreasing mesh scale by fixing the mesh Knudsen number. The level of resolution of shocks relative to the mesh scale afforded by the present Boltzmann--BGK approach \textit{without any numerical shock capturing scheme} rival that of a high-order discontinuous spectral element method augmented with highly-resolved \textit{a posteriori} WENO-type reconstruction \citep{Dumbser2014}. These results indicate that the proposed approach can be a very effective tool for directly resolving shock structures in more complex fluid flows.

\subsection{Kelvin--Helmholtz Instability}
As an initial evaluation of the Boltzmann--BGK approach in predicting turbulent flow phenomena, the onset and development of a two-dimensional Kelvin--Helmholtz instability was computed, where two counter-flowing fluids of varying density exhibit a fluid instability which transitions into complex vortical flow. The domain is set as periodic on $\Omega^{\mathbf{x}} = [-0.5, 0.5]^2$ and the initial conditions are given as
\begin{equation}
    \mathbf{q}(x, y, 0) = 
    \begin{cases}
    \mathbf{q}_L, \quad \mathrm{if}\ |y| \leq 0.25,\\
    \mathbf{q}_R, \quad \mathrm{else},
    \end{cases}
    \quad \mathrm{where} \quad 
    \mathbf{q}_L = 
    \begin{bmatrix}
          2 \\
          0.5 \\
          V(x, y) \\
          6.25
    \end{bmatrix},
    \quad
    \mathbf{q}_R = 
    \begin{bmatrix}
          1 \\
          -0.5 \\
          V(x, y) \\
          6.25
    \end{bmatrix},
\end{equation}
which yields an Atwood number of $1/3$ and a maximum Mach number of $0.2$. To seed instabilities in the flow, a single mode is excited similarly to \citet{Springel2010} by setting the initial vertical velocity component as
\begin{equation}
    V(x, y) = \kappa \sin (\omega \pi x)\left [ \exp \left (-\frac{(y - 0.25)^2}{2\sigma^2} \right) +  \exp \left (-\frac{(y + 0.25)^2}{2\sigma^2}\right) \right],
\end{equation}
where the parameters $\kappa = 0.02$, $\omega = 128$, and $\sigma = 0.0035$ dictate the perturbation strength, frequency, and vertical decay, respectively. The collision time $\tau$ was set constant according to a reference dynamic viscosity $\mu = \tau P = 10^{-5}$ based on the initial static pressure. Due to the low Mach number, the variations in pressure, and therefore viscosity, were relatively low.

To verify the ability of the proposed scheme to predict the onset of flow instabilities consistent with the hydrodynamic equations, a comparison between the solution computed with the Boltzmann--BGK equation and the solution computed with the Navier--Stokes equations was performed. For both cases, the solution was computed on a uniform $400^2$ quadrilateral mesh with a $\mathbb P_3$ scheme, and to reduce the overall computational cost, the monatomic case was chosen ($\delta = 0$, $\gamma = 2$). For the velocity space discretization of the Boltzmann--BGK approach, the resolution was set as $N_r = N_{\phi} = 32$, yielding approximately $2.62$ billion degrees of freedom. To maintain consistency between the two methods, the Navier--Stokes approach was computed with $Pr = 1$, $\gamma = 2$, and $\mu = 10^{-5}$, and the entropy filtering method of \citet{Dzanic2022} was used to stabilize the high-order scheme in the vicinity of discontinuities. 

\begin{figure}[htbp!]
        \centering
        \subfloat[Boltzmann--BGK]{\adjustbox{width=\linewidth,valign=b}{
            \includegraphics[width=\textwidth]{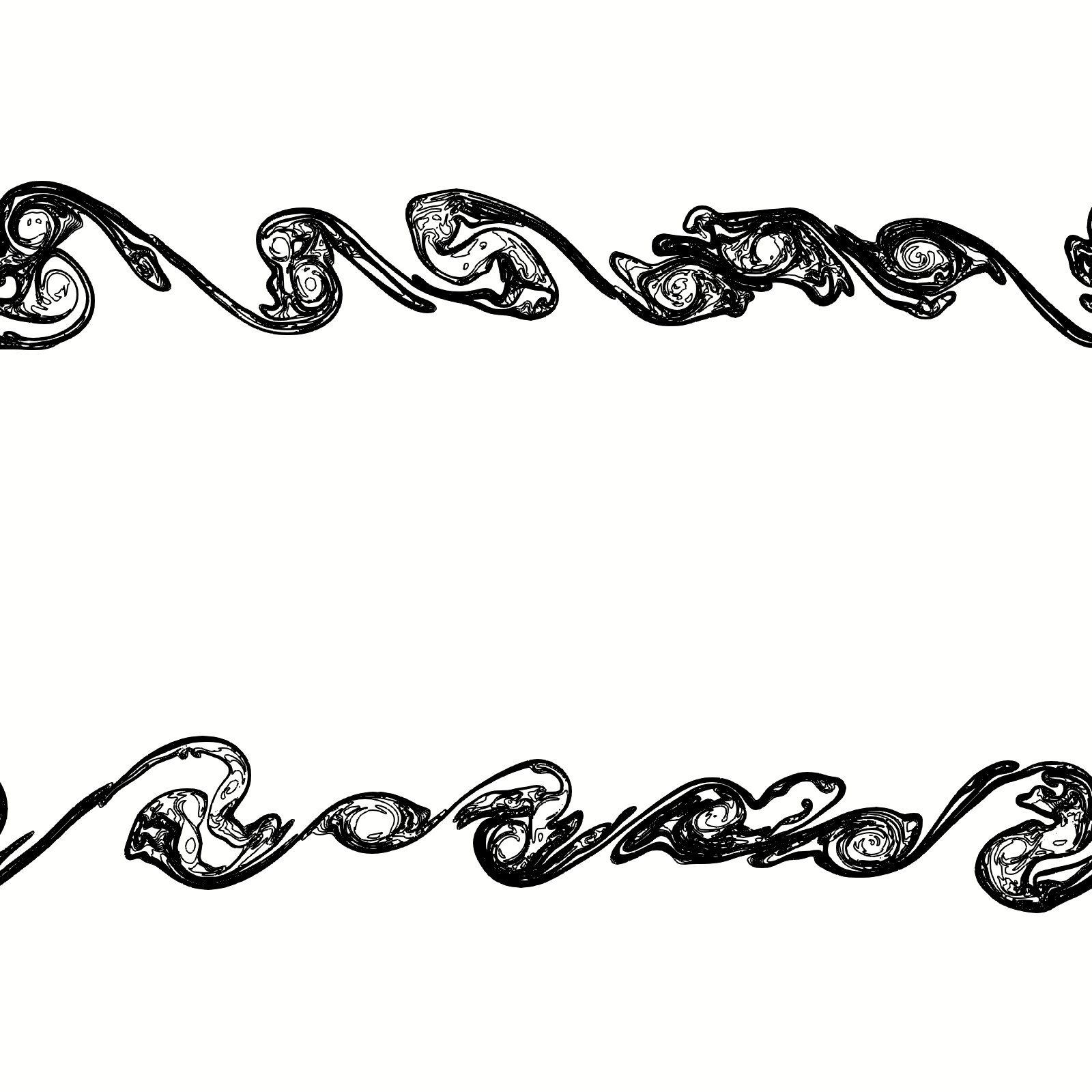}
            \quad \includegraphics[width=\textwidth]{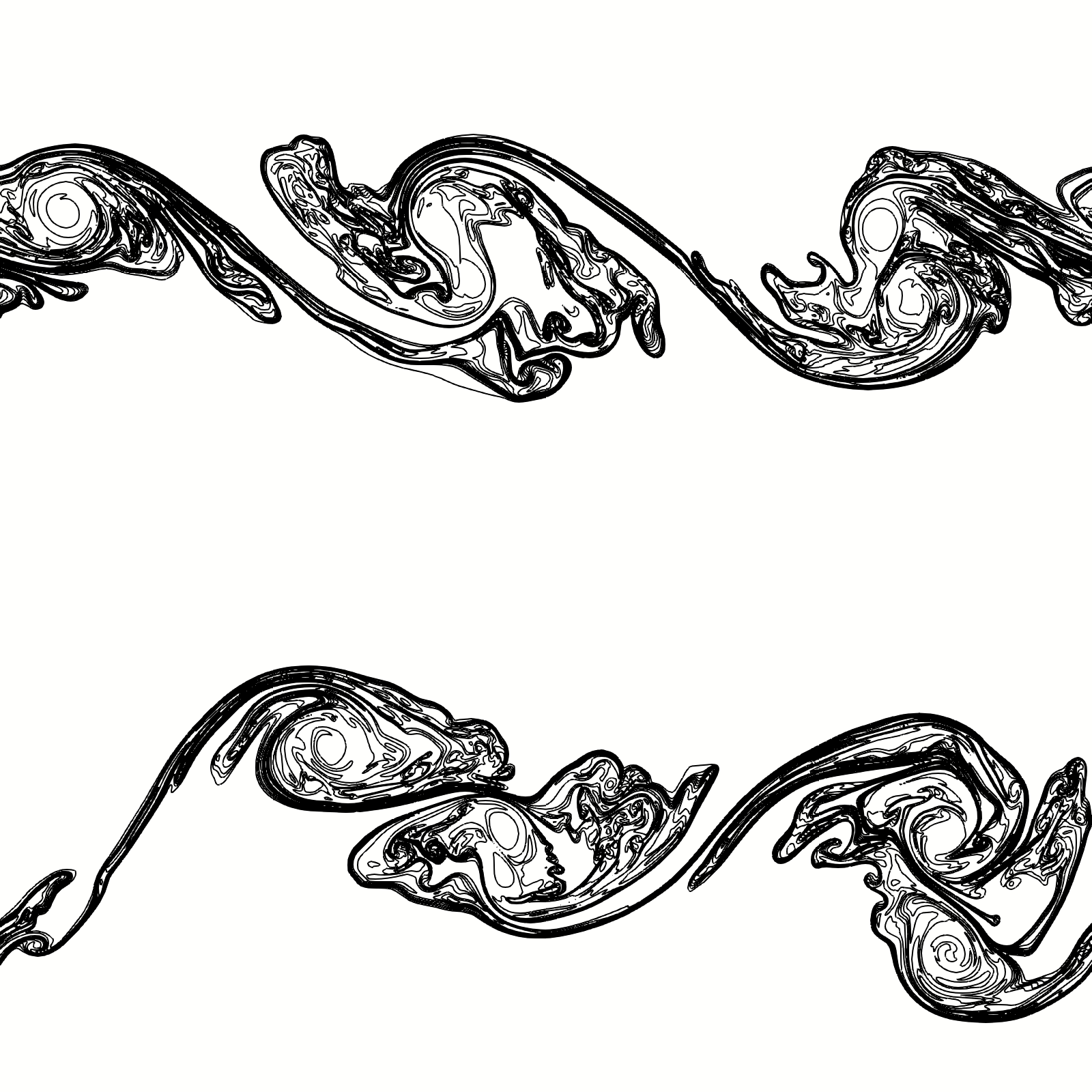}
            \quad \includegraphics[width=\textwidth]{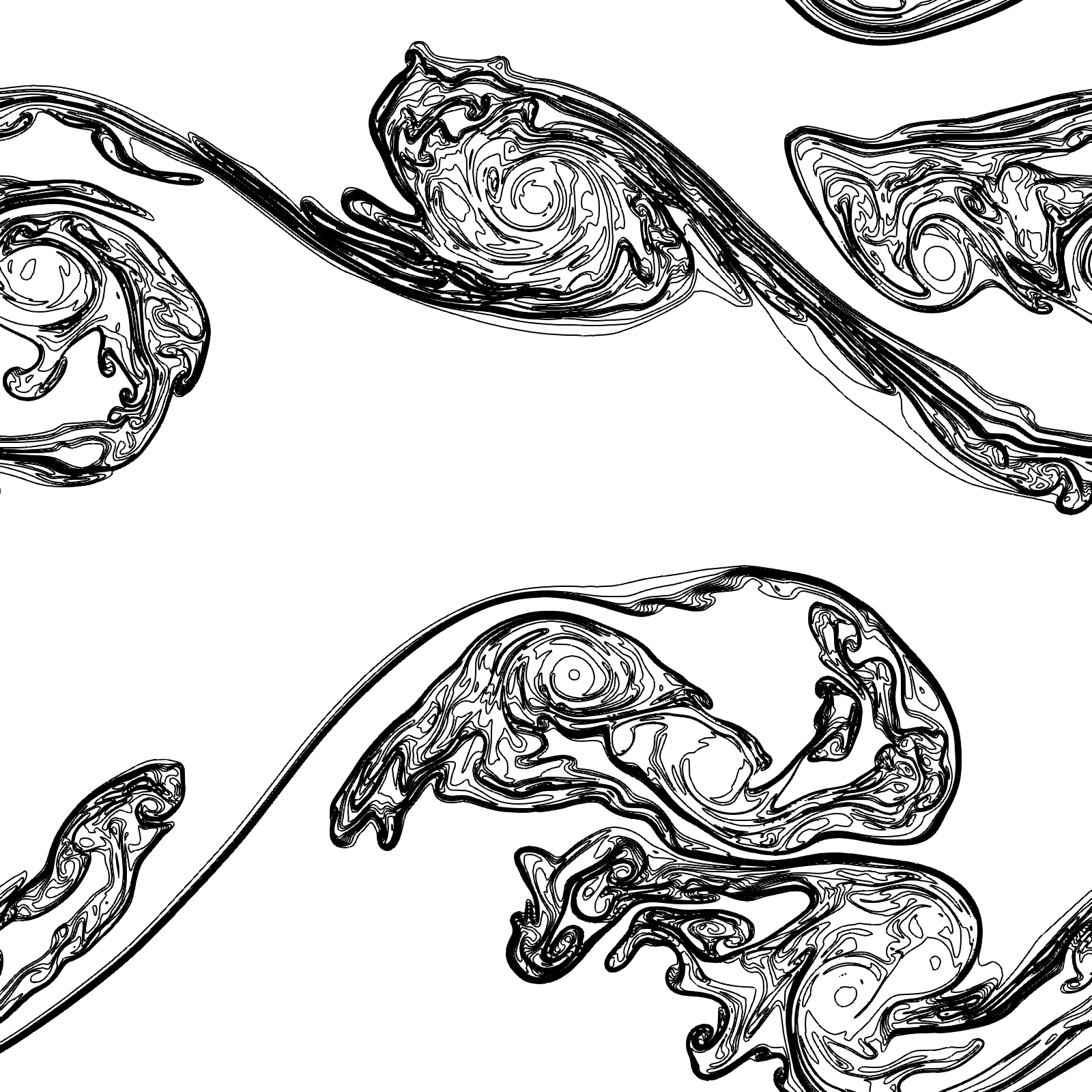}
            \quad \includegraphics[width=\textwidth]{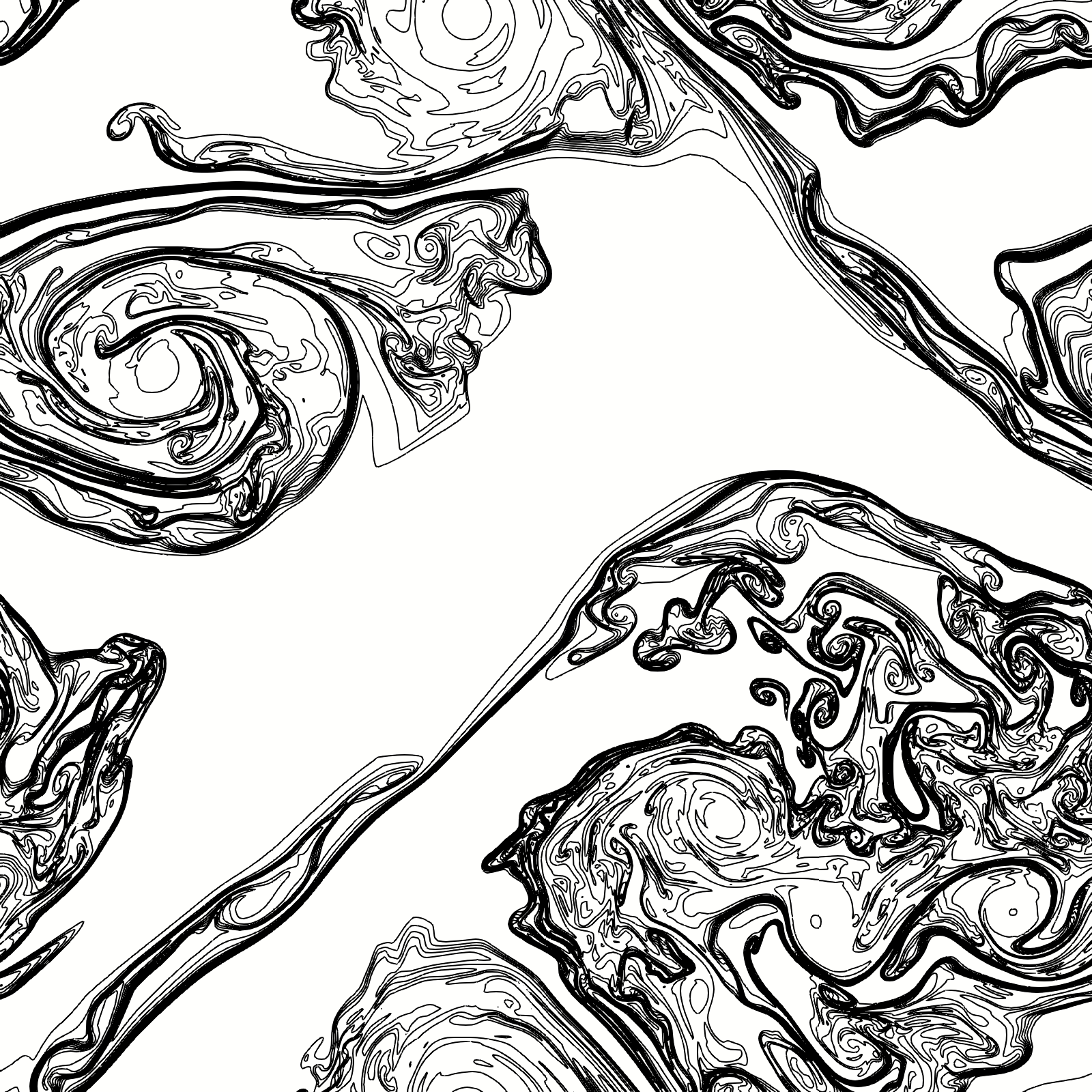}
        }}
        \newline
        \subfloat[Navier--Stokes]{\adjustbox{width=\linewidth,valign=b}{
            \includegraphics[width=\textwidth]{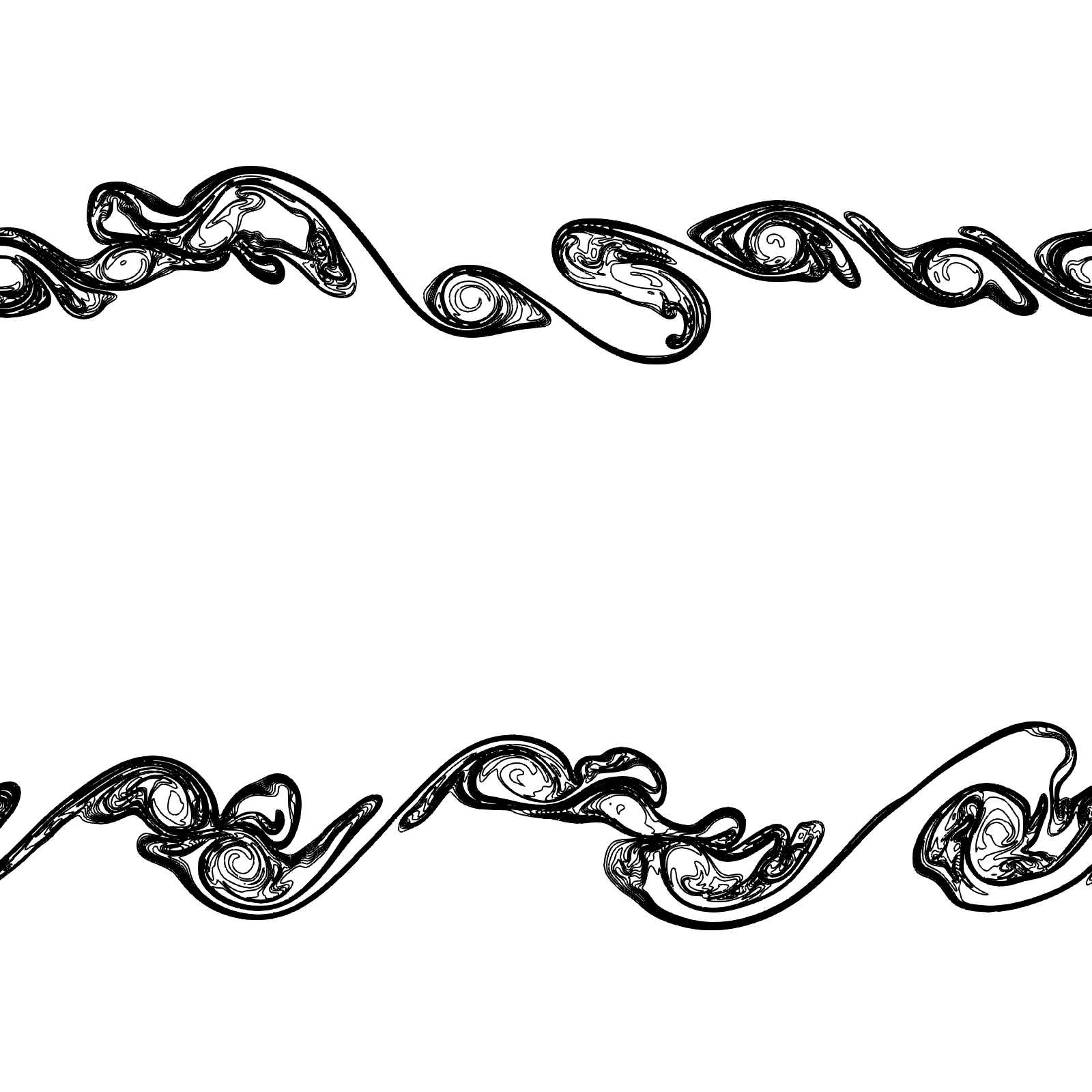}
            \quad \includegraphics[width=\textwidth]{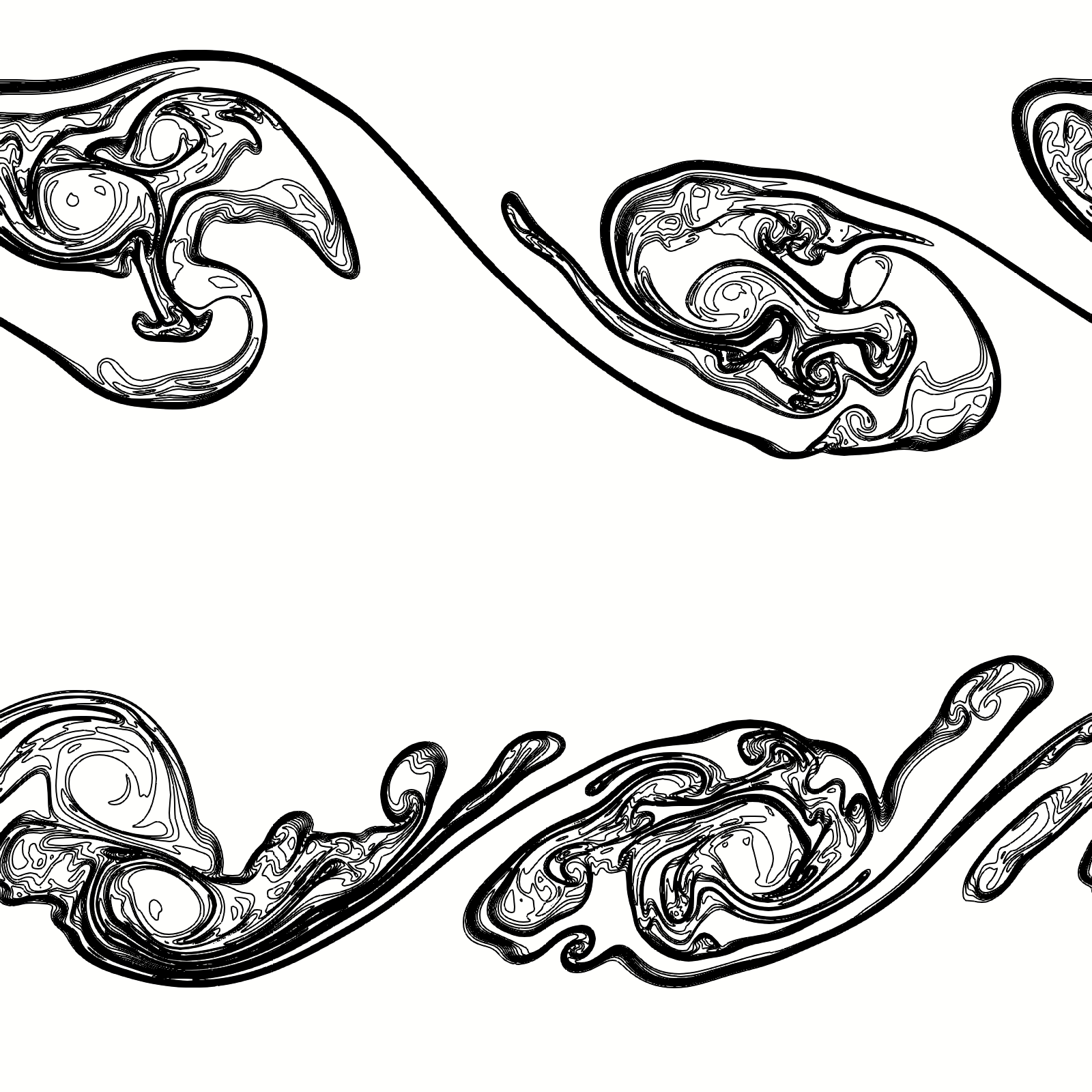}
            \quad \includegraphics[width=\textwidth]{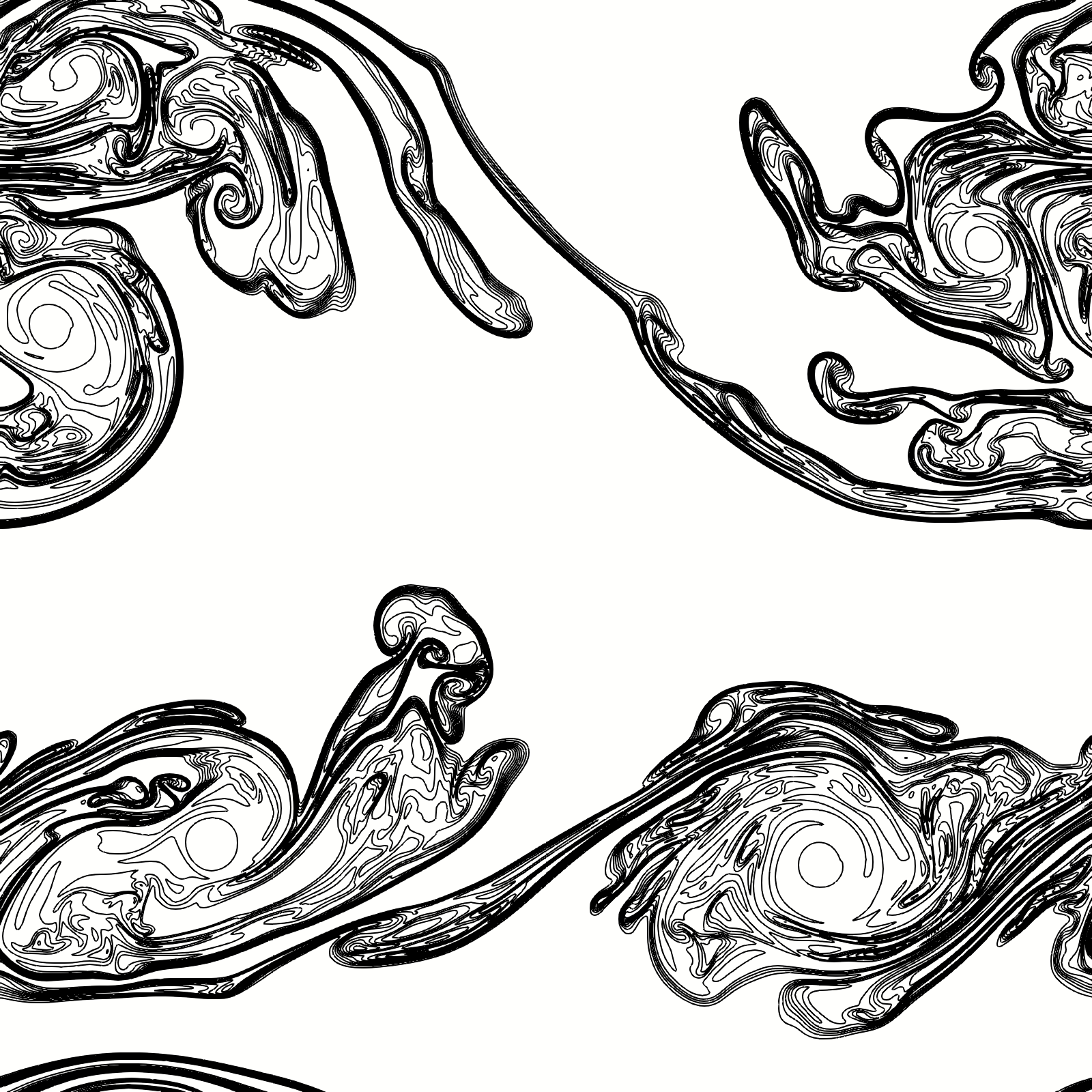}
            \quad \includegraphics[width=\textwidth]{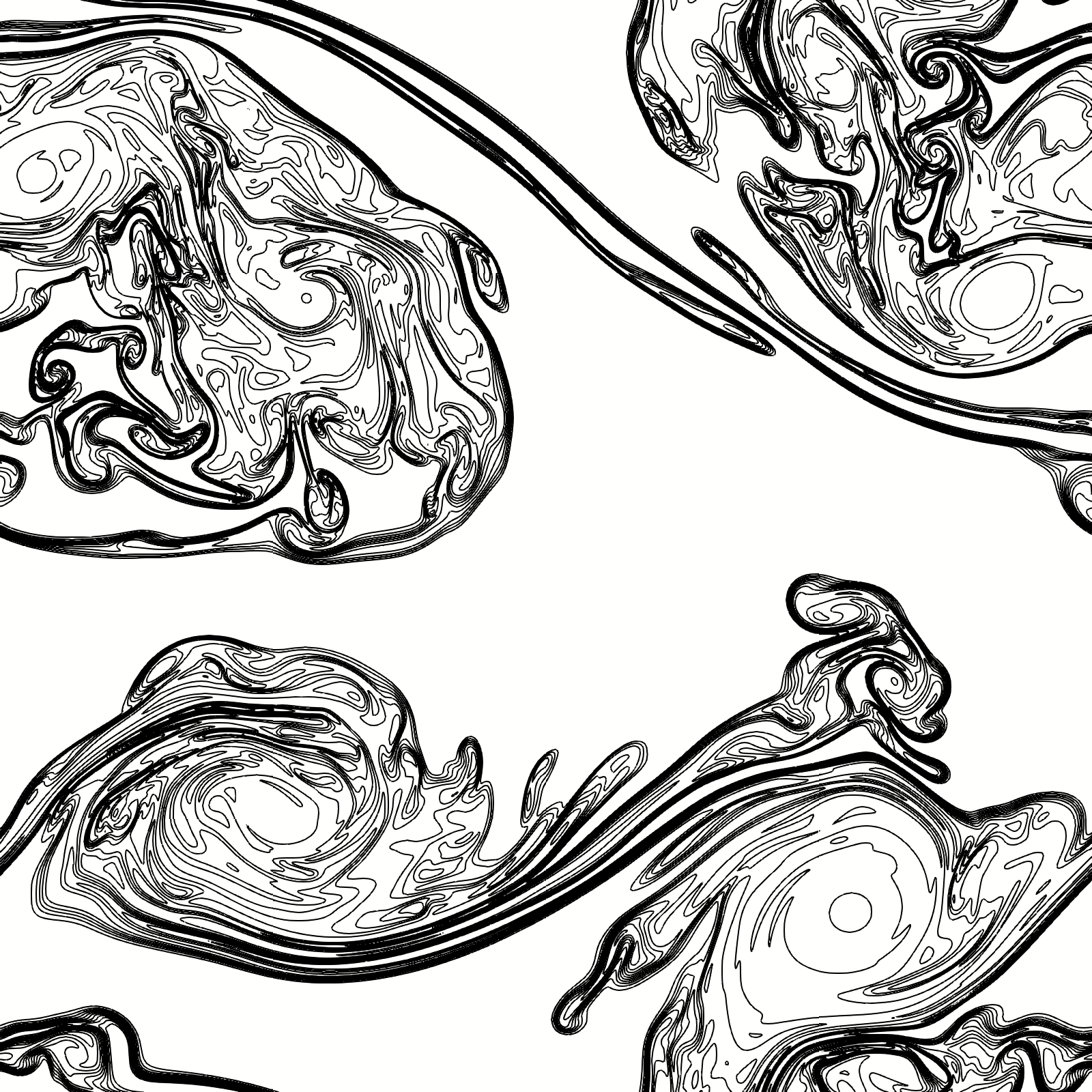}
        }}
        \newline
        \caption{\label{fig:kh} Isocontours of density (equispaced on the range $[1,2]$) for the Kelvin--Helmholtz instability problem at $t = 1$ (left), $t = 2$ (middle-left), $t = 3$ (middle-right), and $t=4$ (right) computed using a $\mathbb P_3$ scheme with $400^2$ elements and $\mu = 10^{-5}$. Top row: Boltzmann--BGK with $N_r = N_{\phi} = 32$ and $\delta = 0$ ($\gamma = 2$). Bottom row: Navier--Stokes with $Pr = 1$ and $\gamma = 2$. }
    \end{figure}

The predicted density profiles, represented as equispaced isocontours on the range $[1, 2]$, at unit time intervals are shown in \cref{fig:kh} for both the Boltzmann--BGK and Navier--Stokes approaches. It can be seen that even though transport in the Boltzmann--BGK equation is linear, the predicted results evidently show the development of the instability into complex vortical flow. In fact, in the earlier time intervals, the large-scale coherent structures of the flow were almost identically computed by the Navier--Stokes and the Boltzmann--BGK approaches. The two flow fields did start to diverge at later times, but this is expected due to the chaotic nature of the Kelvin--Helmholtz instability. Furthermore, the range of scales of the flow features were generally similar between the two approaches, indicating that the turbulent length scales are properly resolved by the Boltzmann--BGK approach. Some discrepancies in the smallest scales did exist between the two approaches, with the Boltzmann--BGK approach showing marginally smaller flow features. This may be attributed to the numerical dissipation introduced by the discontinuity capturing approach for the Navier--Stokes method which artificially inflates the physical viscosity and the minor fluctuations in viscosity of the Boltzmann--BGK approach related to low pressure regions in the vortex cores. However, the marked similarities between the predicted flow fields still indicate that the Boltzmann--BGK approach can accurately resolve multi-scale nonlinear flow phenomena in a manner that is consistent with predictions in the continuum limit. 

\begin{figure}[tbhp]
        \adjustbox{width=0.4\linewidth,valign=b}{\input{figs/kh_enst.tex}}
        \newline
        \caption{\label{fig:kh_enst} 
        Temporal evolution of the volume-integrated enstrophy for the Kelvin--Helmholtz instability problem computed using a $\mathbb P_3$ scheme with $400^2$ elements and $\mu = 10^{-5}$. Boltzmann--BGK results computed with $N_r = N_{\phi} = 32$ and $\delta = 0$ ($\gamma = 2$). Navier--Stokes results computed with $Pr = 1$ and $\gamma = 2$.
        }
    \end{figure}
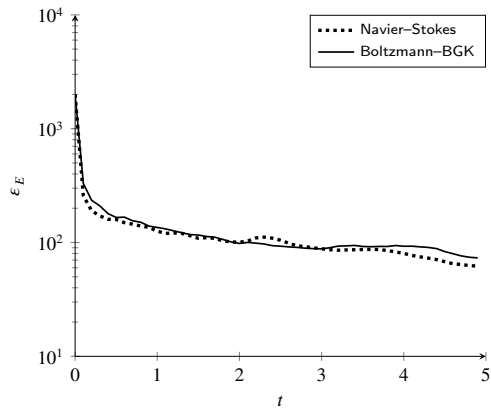

For a more quantitative comparison, the evolution of the volume-integrated enstrophy, defined as 
\begin{equation}
    \varepsilon_E = \frac{1}{V}\int_{\Omega^{\mathbf{x}}} \frac{1}{2}\rho \boldsymbol{\omega}{\cdot} \boldsymbol{\omega} \ \mathrm{d}{\mathbf{x}},
\end{equation}
where $\boldsymbol{\omega}$ is the vorticity and $V = 1$ is the domain volume, was calculated for both approaches. The predicted enstrophy profiles are shown in \cref{fig:kh_enst} over the range $0 \leq t \leq 5$. Very good agreement between the Navier--Stokes results and the Boltzmann--BGK results was observed, both with the early rapid decay of enstrophy associated with the initial velocity perturbations as well as with the long-term evolution of the enstrophy profile driven by turbulent mixing and viscous dissipation. These results give an initial quantitative verification of the ability of the Boltzmann--BGK approach to resolve complex fluid flows. 

\subsection{Compressible Taylor--Green Vortex}
As the final and most rigorous test case, the transition to turbulence of a three-dimensional compressible Taylor--Green vortex was studied. To the authors' knowledge, this is the first direct numerical simulation of three-dimensional turbulent flows performed by directly solving the Boltzmann equation, made possible by the efficiency of the proposed approach. This canonical test case consists of an initially laminar flow that transitions as time evolves, forming large coherent vortical structures that eventually break down to small-scale isotropic turbulence. The problem is solved the periodic domain $\Omega^{\mathbf{x}} = [-\pi, \pi]^3$, and the initial conditions are given as 
    \begin{equation}
         \mathbf{q}(x,y,z,0) =  \begin{bmatrix}
            1 \\ \hphantom{-}\sin(x)\cos(y)\cos(z)\\ -\cos(x)\sin(y)\cos(z) \\ 0\\ P_0 + \frac{1}{16}\left (\cos(2x) + \cos(2y) \right)\left (\cos(2z + 2) \right)
        \end{bmatrix}.
    \end{equation}
The reference pressure was set as $P_0 = 1/(\gamma M^2)$, where $\gamma = 5/3$ and $M = 0.2$. The monatomic case ($\delta = 0$) was chosen to make the computational cost tractable and the Mach number was slightly increased over its typical value of $M = 0.08$ to decrease the stiffness of the source term. The collision time was set constant as $\tau = \mu/P_0$, where the reference dynamic viscosity was set as $\mu_{ref} = 1/1600$ corresponding to a Reynolds number of 1600. Through \textit{a posteriori} analysis, the variation in the viscosity due to the fixed collision time was deemed to be negligible due to the low Mach number. 

The problem was solved on a series of increasingly-resolved meshes and compared to a standard Navier--Stokes approach. As a point of reference, direct numerical simulation was performed by solving the flow field with the same Navier--Stokes approach ($\gamma = 5/3$, $Pr = 1$) on a $N_e = 64^3$ mesh with a $\mathbb P_3$ scheme, a suitable level of resolution for capturing all statistically significant physical scales \citep{Cox2021}. Comparison was then performed on a series of meshes with $N_e = 16^3$, $32^3$, and $48^3$, also discretized with a $\mathbb P_3$ scheme. For the Boltzmann--BGK approach, the velocity domain was discretized with $N_r = 16$,  $N_{\phi} = 16$, and $N_{\psi} = 8$, which will later be shown to be sufficiently resolved. For the finest case, where the resolution level corresponded to approximately 14.5 billion degrees of freedom, the computational cost was on the order of 3000 GPU-hours (72 hours on 40 GPUs). In comparison, the cost of the Navier--Stokes approach for the same mesh was approximately 16 GPU hours.

   \begin{figure}[tbhp]
        \subfloat[$N_e = 16^3$]{\adjustbox{width=0.33\linewidth, valign=b}{\input{figs/tgv_64_enst}}}
        \subfloat[$N_e = 32^3$]{\adjustbox{width=0.33\linewidth, valign=b}{\input{figs/tgv_128_enst}}}
        \subfloat[$N_e = 48^3$]{\adjustbox{width=0.33\linewidth, valign=b}{\input{figs/tgv_192_enst}}}
        \newline
        \caption{\label{fig:tgv_enst}  Dissipation measured by enstrophy for the compressible Taylor-Green vortex at $Re = 1600$ computed using a $\mathbb P_3$ scheme with $16^3$ (left), $32^3$ (middle), and $48^3$ (right) elements. Boltzmann--BGK results computed with $N_r = 16$,  $N_{\phi} = 16$, $N_{\psi} = 8$, and $\delta = 0$ ($\gamma = 5/3$). Navier--Stokes results computed with $Pr = 1$ and $\gamma = 5/3$. 
        }
    \end{figure}
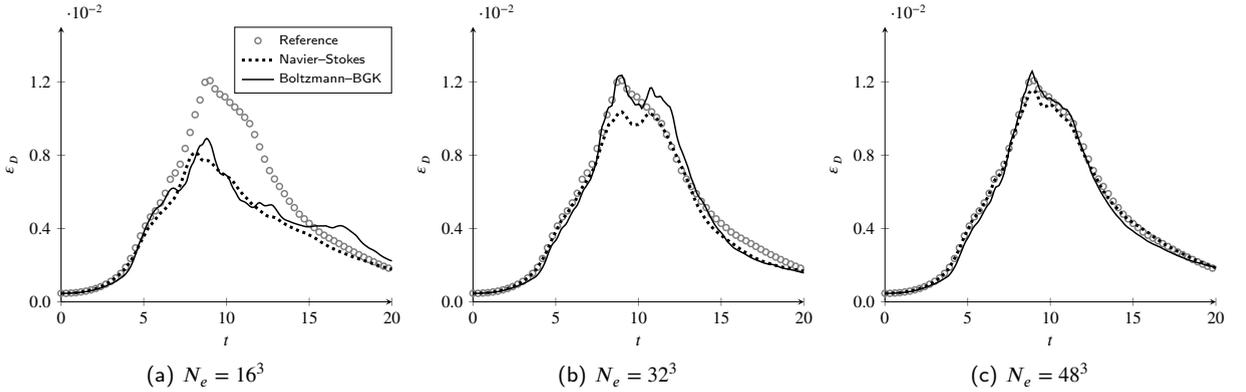
    
The standard metric of comparison for the Taylor--Green vortex is the dissipation $\varepsilon_D$ measured by a scaled form of the non-dimensional volume-integrated enstrophy, given as 
\begin{equation}
    \varepsilon_D = \frac{\beta}{V}\int_{\Omega^{\mathbf{x}}} \frac{1}{2}\rho \boldsymbol{\omega}{\cdot} \boldsymbol{\omega} \ \mathrm{d}{\mathbf{x}} = \beta \varepsilon_E,
\end{equation}
where $\beta = 2 \mu$ is the scaling factor and $V = 8 \pi^3$ is the domain volume. The viscosity was taken as a constant corresponding to the reference value $\mu_{ref}$, but the differences to the enstrophy values computed using the pressure-dependent viscosity $\mu = \tau P$ were confirmed to be negligible. The profiles of the enstrophy-based dissipation rate are shown in \cref{fig:tgv_enst} for the various meshes. For the coarsest mesh, the enstrophy was underpredicted for both the Boltzmann--BGK and the Navier--Stokes approaches in comparison to the reference results due to the low levels of resolution, However, the profiles showed relatively good agreement between the two approaches, with the Boltzmann--BGK results showing some more oscillations in the profile and an overprediction of the enstrophy at later times. When the resolution was increased to the medium mesh, the enstrophy profile of the Boltzmann--BGK approach showed much better agreement with the reference results, with better agreement at the enstrophy peak than the Navier--Stokes results at the same resolution but with a spurious secondary peak at a later time. For the finest mesh, the results of the Boltzmann--BGK approach converged excellently to the reference results with nearly an identical enstrophy profile. Furthermore, these results were consistent with the work of \citet{Cao2022} which utilized a gas kinetic scheme derived from the Chapman–-Enskog expansion of the Boltzmann--BGK equation that directly solves for the conserved flow variables. 

   \begin{figure}[tbhp]
        \subfloat[]{\adjustbox{width=0.4\linewidth, valign=b}{\input{figs/tgv_128_velconv}}}
        \subfloat[]{\adjustbox{width=0.4\linewidth, valign=b}{\input{figs/tgv_128_conservation}}}
        \newline
        \caption{\label{fig:tgv_cons}  Convergence of dissipation measured by enstrophy (left) and mass conservation error (right) with respect to velocity space resolution $N_v = N_r \times N_{\phi} \times N_{\psi}$ for the Boltzmann-BGK approximation of the compressible Taylor--Green vortex at $Re = 1600$ computed using a $\mathbb P_3$ scheme with the DVM and $N_e = 32^3$. Default velocity space resolution shown in red. 
        }
    \end{figure}
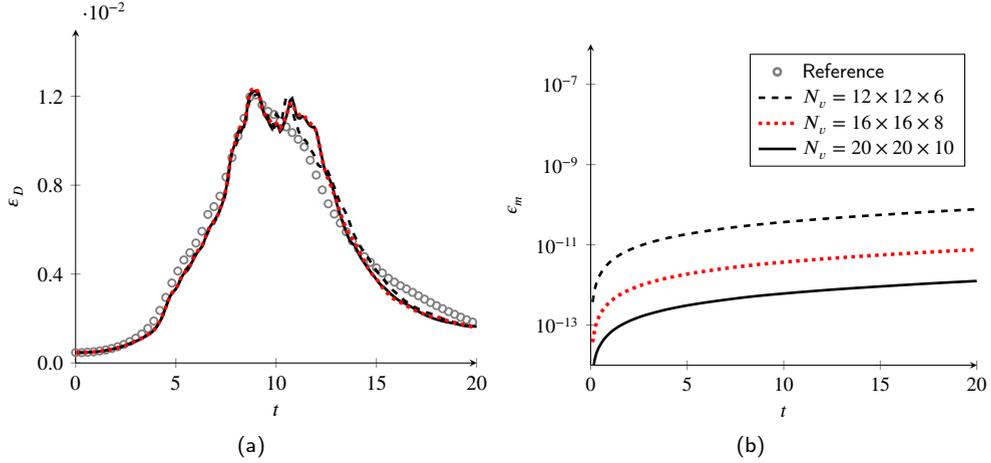

To verify that the results were converged in the velocity domain and the scheme remained conservative, the spatial resolution was fixed at $N_e = 32^3$ and the velocity resolution was modulated. The predicted enstrophy profiles and the mass conservation error for the various velocity resolution levels are shown in \cref{fig:tgv_cons}. In comparison to the standard resolution used in this study ($N_v = 16 \times 16 \times 8$), lowering the resolution to $N_v = 12 \times 12 \times 6$ resulted in a very similar enstrophy profiles, with just a slightly deviation in the second spurious enstrophy peak. When the resolution was increased to $N_v = 20 \times 20 \times 10$, the results were essentially indistinguishable from the standard resolution, indicating that the results were well-converged in the velocity domain for the choice of velocity resolution used for the experiments. Furthermore, regardless of the choice of the resolution, the scheme remained conservative, with a mass conservation error of less than $10^{-10}$ purely as a result of the linear accumulation of mass conservation error roughly on the order of machine precision per time step. 
    
    \begin{figure}[htbp!]
        \centering
        
        \subfloat[Boltzmann-BGK, $N_e = 16^3$] {\adjustbox{width=0.33\linewidth,valign=b}{
            \includegraphics[width=\textwidth]{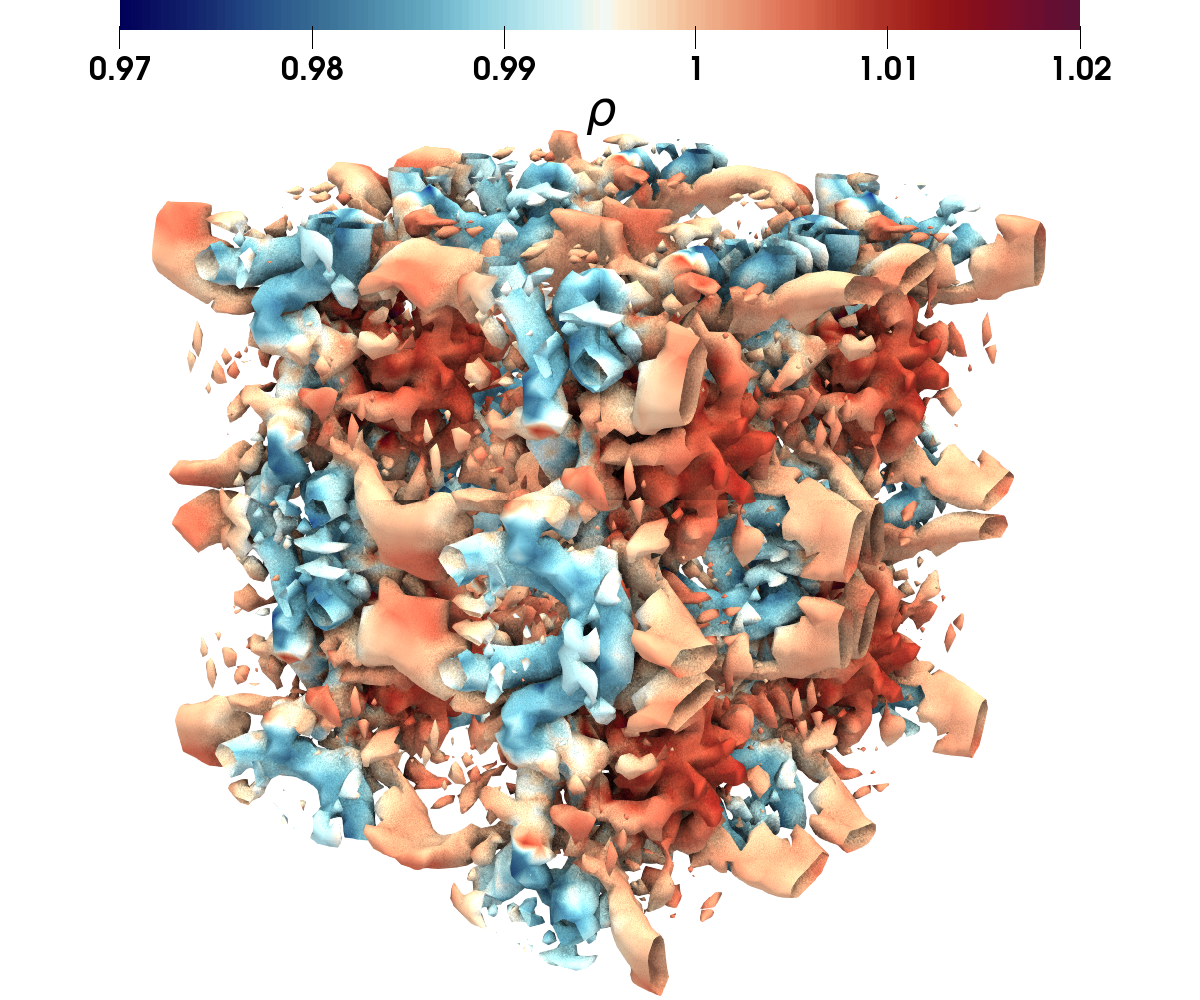}
        }}
        \subfloat[Navier--Stokes, $N_e = 16^3$] {\adjustbox{width=0.33\linewidth,valign=b}{
            \includegraphics[width=\textwidth]{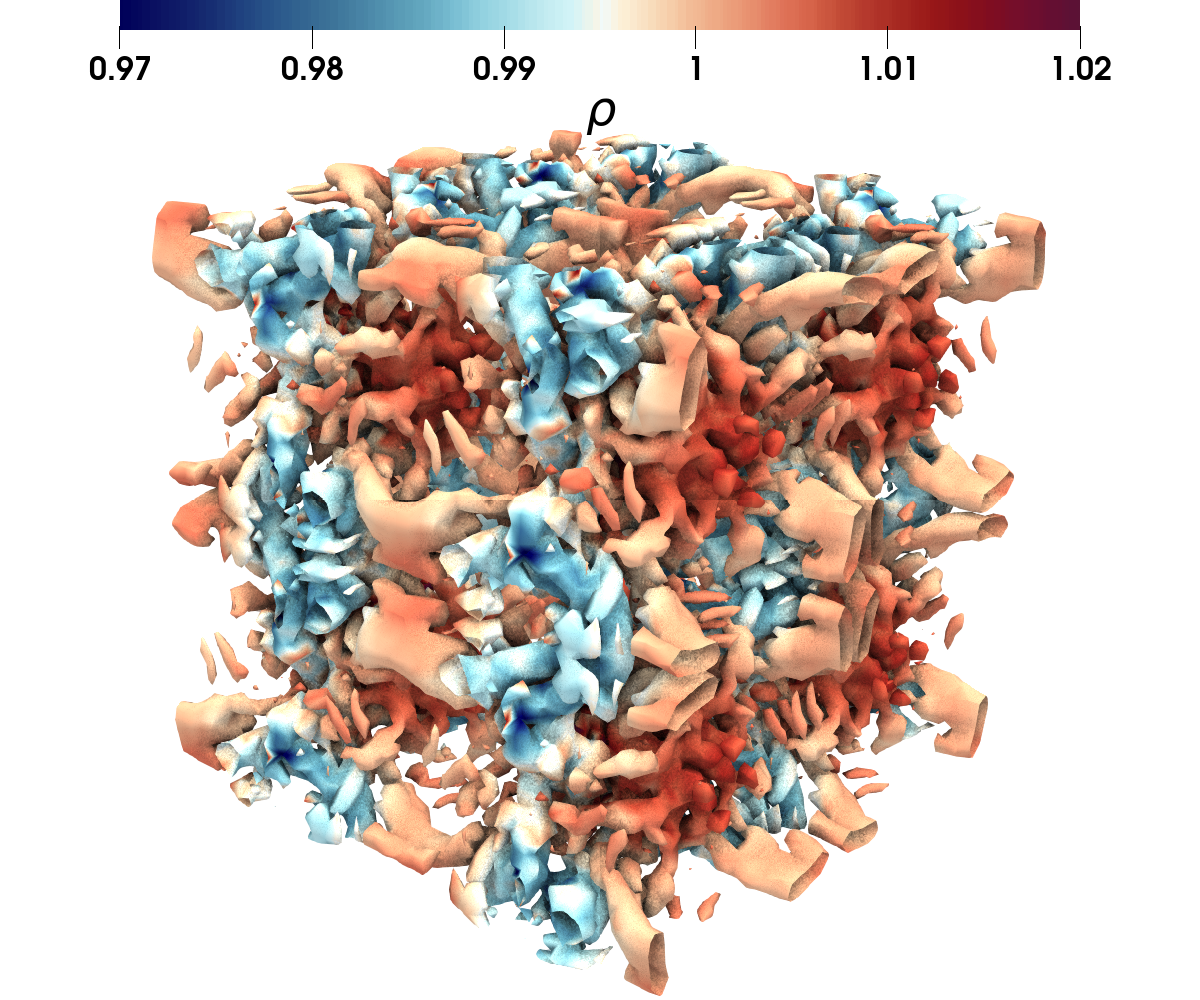}
        }}
        \newline
        \subfloat[Boltzmann-BGK, $N_e = 32^3$] {\adjustbox{width=0.33\linewidth,valign=b}{
            \includegraphics[width=\textwidth]{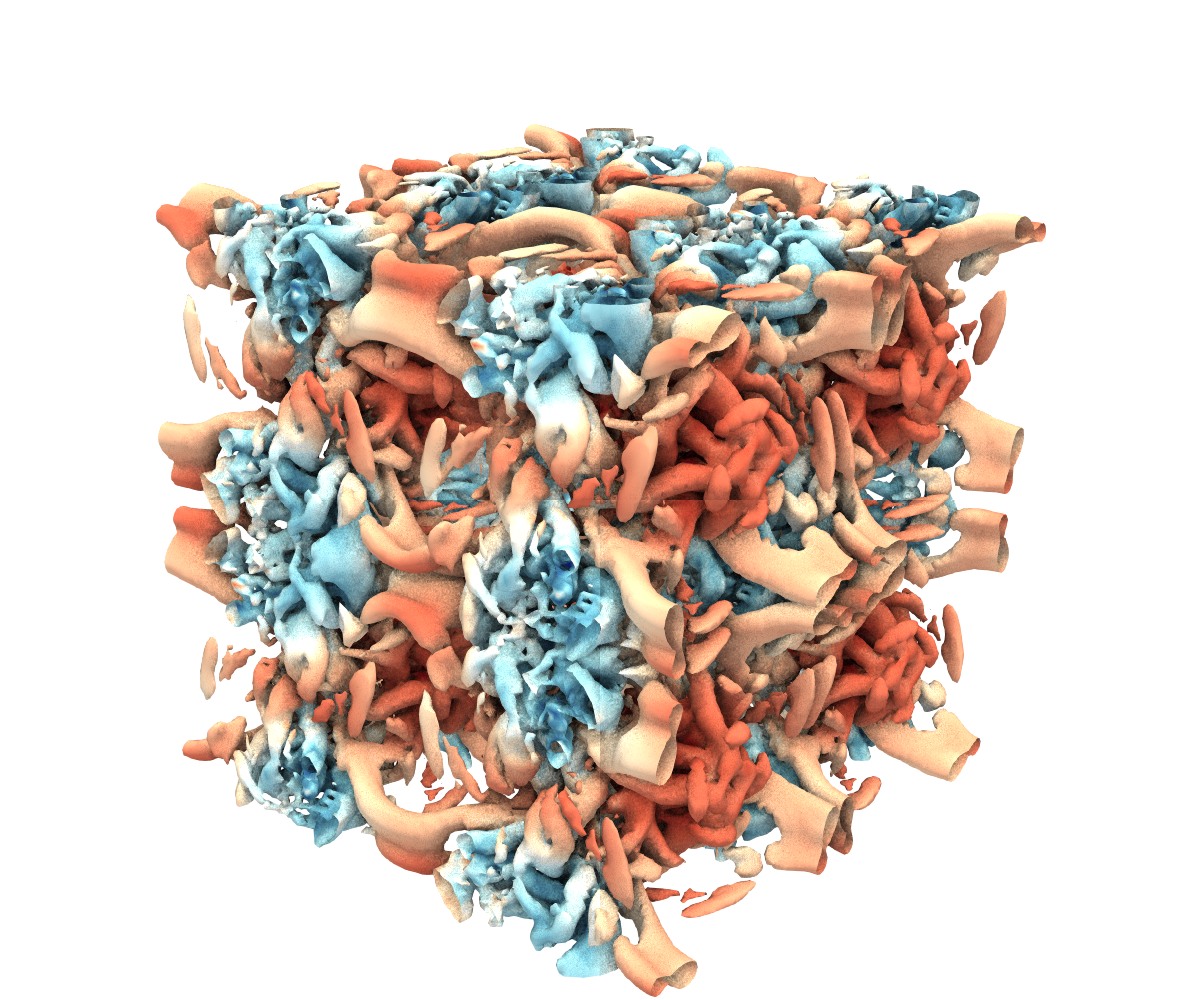}
        }}
        \subfloat[Navier--Stokes, $N_e = 32^3$] {\adjustbox{width=0.33\linewidth,valign=b}{
            \includegraphics[width=\textwidth]{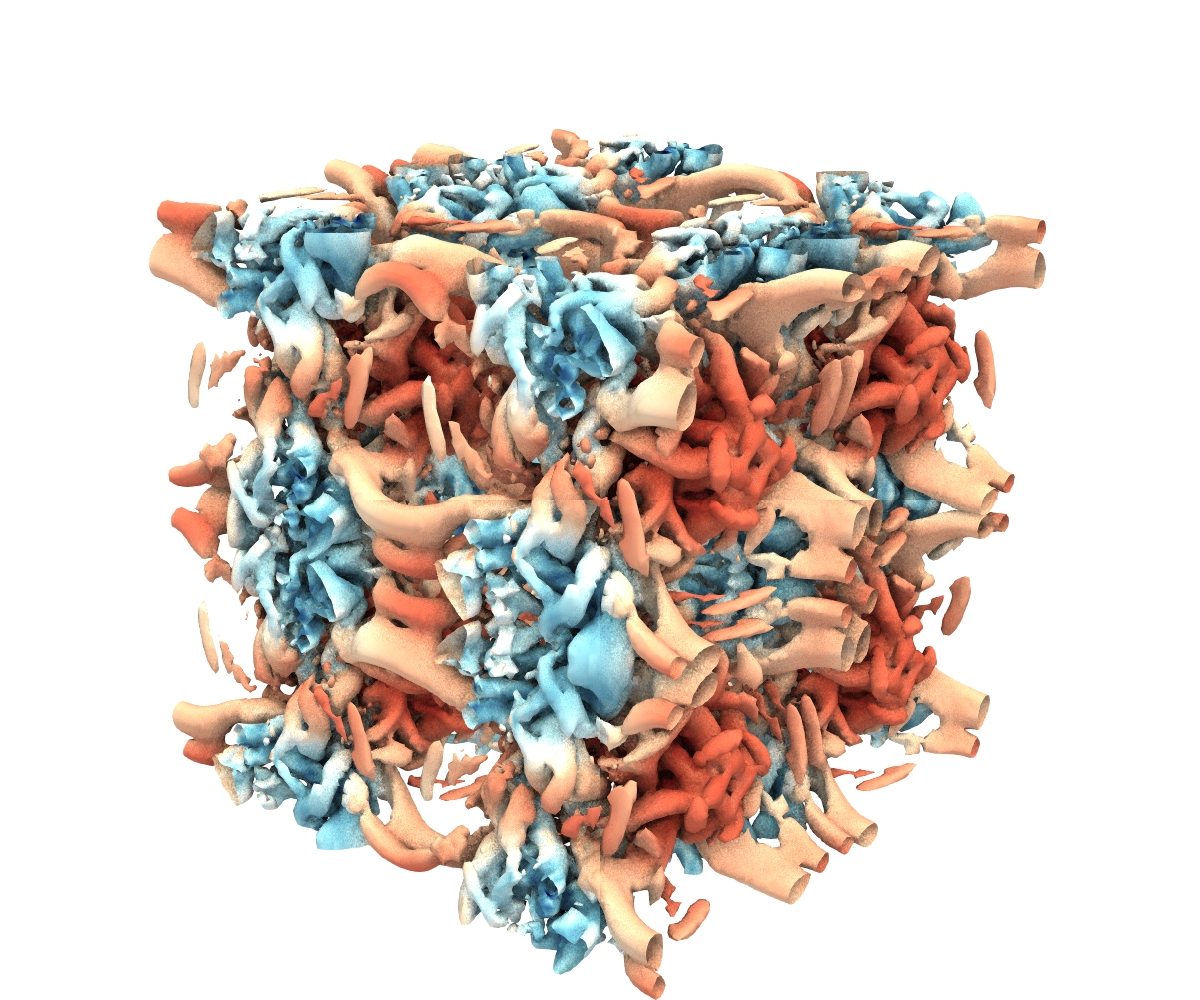}
        }}
        \newline
        \subfloat[Boltzmann-BGK, $N_e = 48^3$] {\adjustbox{width=0.33\linewidth,valign=b}{
            \includegraphics[width=\textwidth]{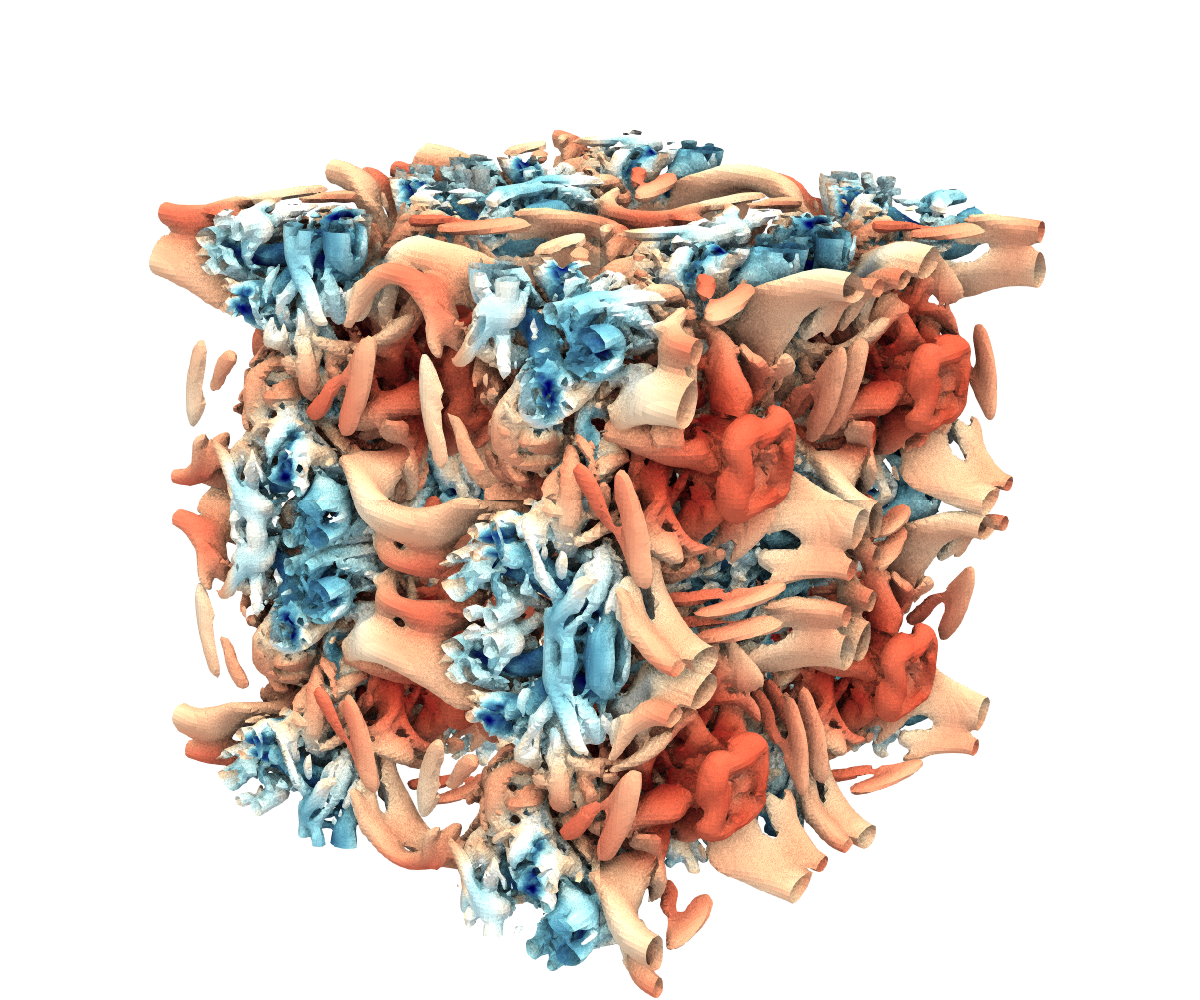}
        }}
        \subfloat[Navier--Stokes, $N_e = 48^3$] {\adjustbox{width=0.33\linewidth,valign=b}{
            \includegraphics[width=\textwidth]{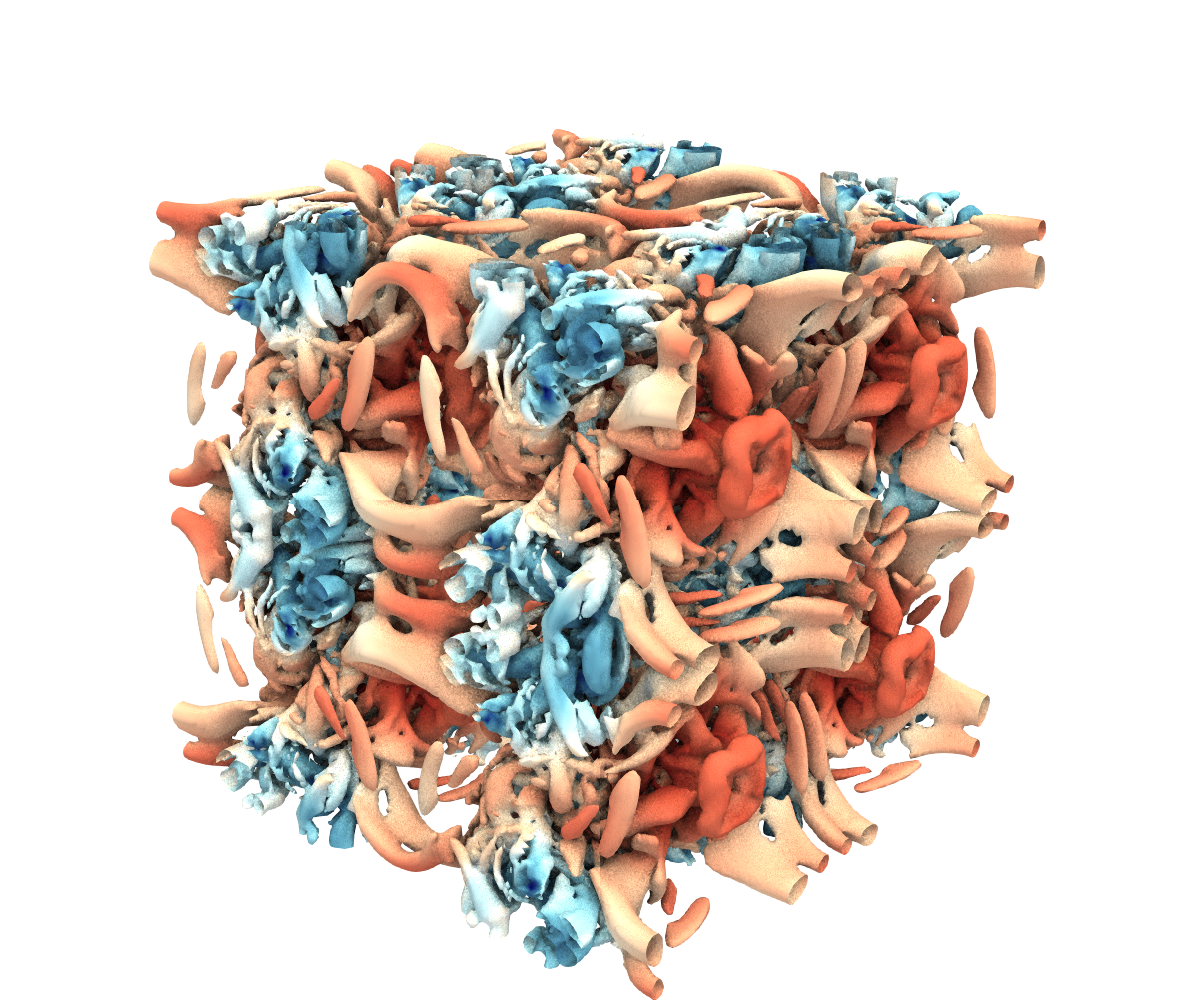}
        }}
        \newline
        \subfloat[Reference] {\adjustbox{width=0.33\linewidth,valign=b}{
            \includegraphics[width=\textwidth]{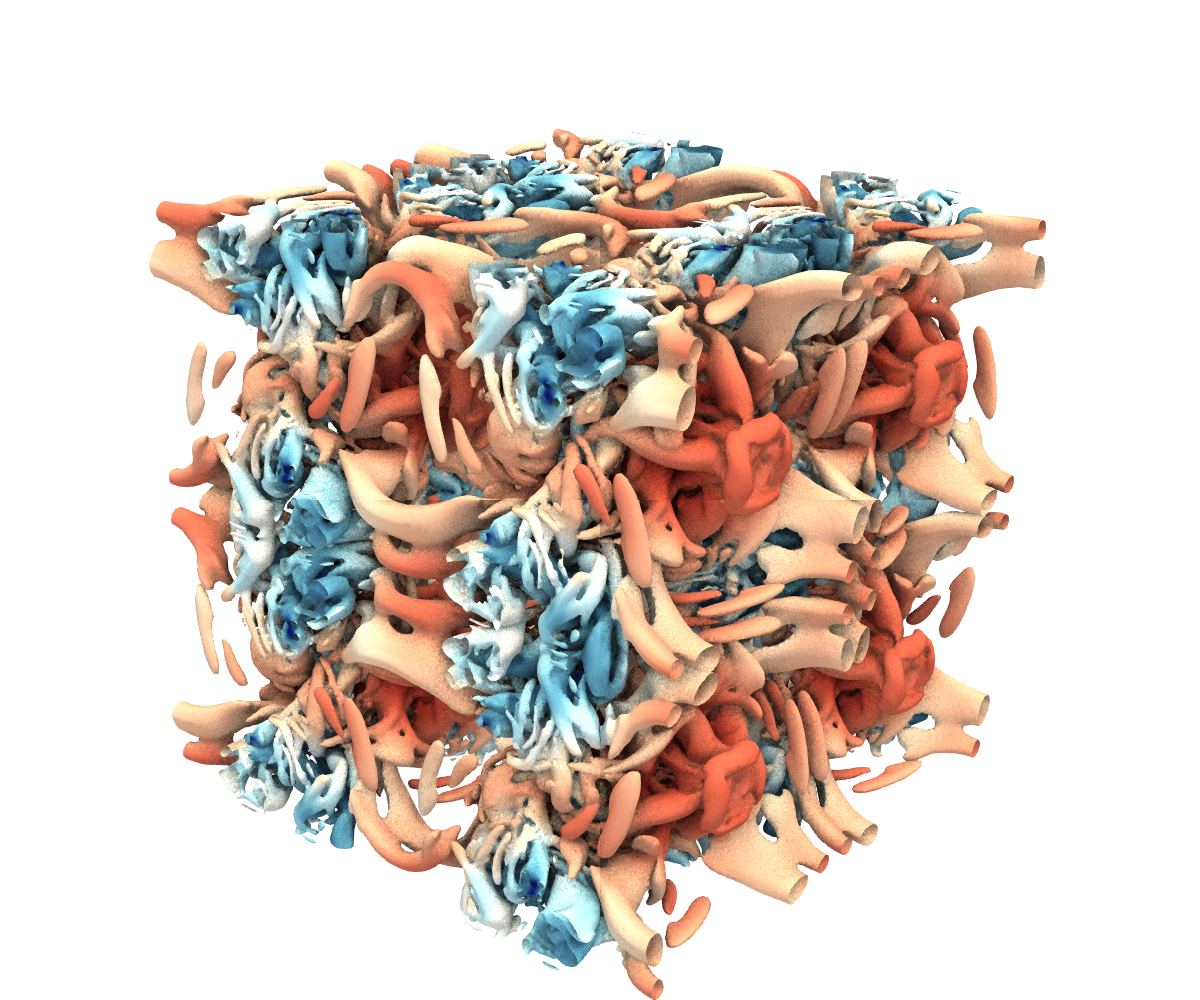}
        }}
        \hspace{50pt}
        
        \caption{\label{fig:tgv_qcrit} Isosurfaces of $Q$-criterion $= 0.2$ colored by density at $t = 10$ for the compressible Taylor--Green vortex at $Re = 1600$ computed using a $\mathbb P_3$ scheme with varying mesh resolution. Left column: Boltzmann--BGK results with $N_r = 16$,  $N_{\phi} = 16$, $N_{\psi} = 8$, and $\delta = 0$ ($\gamma = 5/3$). Right column: Navier--Stokes results with $Pr = 1$ and $\gamma = 5/3$. 
        }
    \end{figure}

A comparison between the flow fields was then performed to confirm that the turbulent structures in the flow were accurately predicted. The isosurfaces of Q-criterion (with $Q = 0.2$) colored by density are shown for the various meshes in \cref{fig:tgv_qcrit} at $t = 10$ where the enstrophy is near its peak. It can be seen that the dominant flow structures were accurately predicted by the Boltzmann--BGK approach in comparison to the reference results, even at the lowest mesh resolution, and these results were very similar to the Navier--Stokes predictions at the same resolution. When the resolution was increased, the predicted results converged excellently to the reference results, with effectively identical prediction of the flow structures and the density distribution within the domain.  

   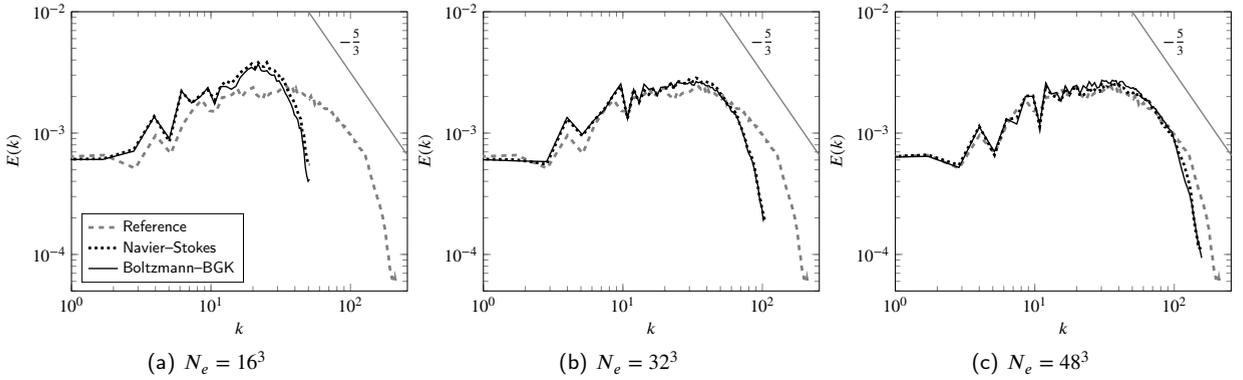
\begin{figure}[tbhp]
        \subfloat[$N_e = 16^3$]{\adjustbox{width=0.33\linewidth, valign=b}{\input{figs/tgv_64_spectra}}}
        \subfloat[$N_e = 32^3$]{\adjustbox{width=0.33\linewidth, valign=b}{\input{figs/tgv_128_spectra}}}
        \subfloat[$N_e = 48^3$]{\adjustbox{width=0.33\linewidth, valign=b}{\input{figs/tgv_192_spectra}}}
        \newline
        \caption{\label{fig:tgv_spectra} 
        Turbulent kinetic energy spectra at $t = 10$ for the compressible Taylor-Green vortex computed using a $\mathbb P_3$ scheme with $16^3$ (left), $32^3$ (middle), and $48^3$ (right) elements. Boltzmann--BGK results computed with $N_r = 16$,  $N_{\phi} = 16$, $N_{\psi} = 8$, and $\delta = 0$ ($\gamma = 5/3$). Navier--Stokes results computed with $Pr = 1$ and $\gamma = 5/3$. 
        }
    \end{figure}

As a final verification of the ability of the Boltzmann--BGK approach in predicting turbulent flow phenomena, the behavior of the turbulent kinetic energy spectra was examined. Similarly to the Q-criterion isosurfaces, this analysis was performed near the enstrophy peak at $t = 10$. The spectra $E(k)$ was computed by taking the three-dimensional fast Fourier transform of the square of the velocity magnitude ($\mathbf{U} \cdot \mathbf{U}$) interpolated onto a uniform grid of the same number of degrees of freedom and then performing a spherical averaging operation with the number of bins set identical to the maximum resolvable wavenumber in one-dimension ($N_e(p+1)/2$). The predicted energy spectra are shown in \cref{fig:tgv_spectra} for the various meshes. It can be seen that the predicted spectra converged to the reference results with increasing resolution. Furthermore, for a given resolution, the Boltzmann--BGK approach remarkably predicts the energy spectra almost identically to the Navier--Stokes approach. These results present a validation of the Boltzmann--BGK approach for three-dimensional turbulent flows and indicate that the approach can accurately resolve turbulent flow phenomena in a consistent manner with respect to the hydrodynamic equations.

%% file: figs/pulse.tex
\begin{tikzpicture}[spy using outlines={rectangle, height=3cm,width=2.3cm, magnification=3, connect spies}]
	\begin{axis}[name=plot1,
		axis line style={latex-latex},
	    axis x line=left,
        axis y line=left,
		xlabel={$x$},
    	xmin=0,xmax=1,
    	xtick={0, 0.2, 0.4, 0.6, 0.8, 1.0},
    	ylabel={$\rho$},
    	ymin=0.9,ymax=2.1,
    	ytick={1.0, 1.2, 1.4, 1.6, 1.8, 2.0},
    	ylabel style={rotate=-90},
    	legend style={at={(0.97, 0.97)},anchor=north east,font=\small},
    	legend cell align={left},
    	style={font=\normalsize}]

    \addplot [domain=0:1, samples=101,unbounded coords=jump, color=gray, style={thick}, only marks, mark=o, mark options={scale=0.8}]{1+exp(-100*(x - 0.5)^2)};
    \addlegendentry{$\rho(x, 0)$};
    
    \addplot[color=black, 
            style={very thick, dotted}]
    table[x=x, y=r, col sep=comma]{./figs/data/pulse_Kn0p1.csv};
    \addlegendentry{$Kn = 10^{-1}$};
    
    \addplot[color=black, 
            style={very thick, dashed}]
    table[x=x, y=r, col sep=comma]{./figs/data/pulse_Kn0p01.csv};
    \addlegendentry{$Kn = 10^{-2}$};
    
    \addplot[color=black, 
            style={very thick}]
    table[x=x, y=r, col sep=comma]{./figs/data/pulse_Kn0p001.csv};
    \addlegendentry{$Kn = 10^{-3}$};

	\end{axis}
\end{tikzpicture}

%% file: figs/velocity_convergence_kn0p1.tex
\begin{tikzpicture}[spy using outlines={rectangle, height=3cm,width=2.3cm, magnification=3, connect spies}]
	\begin{semilogyaxis}[name=plot1,
		axis line style={latex-latex},
		xlabel={$N_v$},
    	xmin=10,xmax=62,
    	ylabel={$\epsilon_{\rho, \infty}$},
		ylabel style={rotate=-90},
    	ymin=1e-10,ymax=10,
    	legend style={at={(0.97, 0.97)},anchor=north east,font=\small},
    	legend cell align={left},
    	style={font=\normalsize}]
    	
        \addplot[color=black, 
                style={thick, dotted},
                mark=square,
                mark options=solid]
        table[x=N, y=Linf, col sep=comma]{./figs/data/velocity_convergence_kn0.1_gamma3.csv};
        \addlegendentry{Standard};
        
        \addplot[color=black, 
                style={thick},
                mark=*,
                mark options=solid,
                mark options={scale=0.6}]
        table[x=N, y=Linf, col sep=comma]{./figs/data/velocity_convergence_kn0.1_gamma3_DVM.csv};
        \addlegendentry{DVM};
        
        \addplot[color=red!80!black, 
                style={thick, dotted},
                mark=square,
                mark options=solid]
        table[x=N, y=Linf, col sep=comma]{./figs/data/velocity_convergence_kn0.1_gamma14.csv};
        
        \addplot[color=red!80!black, 
                style={thick},
                mark=*,
                mark options=solid,
                mark options={scale=0.6} ]
        table[x=N, y=Linf, col sep=comma]{./figs/data/velocity_convergence_kn0.1_gamma14_DVM.csv};

	\end{semilogyaxis}
\end{tikzpicture}

%% file: figs/velocity_convergence_kn0p01.tex
\begin{tikzpicture}[spy using outlines={rectangle, height=3cm,width=2.3cm, magnification=3, connect spies}]
	\begin{semilogyaxis}[name=plot1,
		axis line style={latex-latex},
		xlabel={$N_v$},
    	xmin=10,xmax=62,
    	ylabel={$\epsilon_{\rho, \infty}$},
		ylabel style={rotate=-90},
    	ymin=1e-10,ymax=10,
    	legend style={at={(0.03, 0.03)},anchor=south west,font=\small},
    	legend cell align={left},
    	style={font=\normalsize}]
    	
        \addplot[color=black, 
                style={thick, dotted},
                mark=square,
                mark options=solid]
        table[x=N, y=Linf, col sep=comma]{./figs/data/velocity_convergence_kn0.01_gamma3.csv};
        
        \addplot[color=black, 
                style={thick},
                mark=*,
                mark options=solid,
                mark options={scale=0.8}]
        table[x=N, y=Linf, col sep=comma]{./figs/data/velocity_convergence_kn0.01_gamma3_DVM.csv};
        
        \addplot[color=red!80!black, 
                style={thick, dotted},
                mark=square,
                mark options=solid]
        table[x=N, y=Linf, col sep=comma]{./figs/data/velocity_convergence_kn0.01_gamma14.csv};
        
        \addplot[color=red!80!black, 
                style={thick},
                mark=*,
                mark options=solid,
                mark options={scale=0.6} ]
        table[x=N, y=Linf, col sep=comma]{./figs/data/velocity_convergence_kn0.01_gamma14_DVM.csv};

	\end{semilogyaxis}
\end{tikzpicture}

%% file: figs/velocity_convergence_kn0p001.tex
\begin{tikzpicture}[spy using outlines={rectangle, height=3cm,width=2.3cm, magnification=3, connect spies}]
	\begin{semilogyaxis}[name=plot1,
		axis line style={latex-latex},
		xlabel={$N_v$},
    	xmin=10,xmax=62,
    	ylabel={$\epsilon_{\rho, \infty}$},
		ylabel style={rotate=-90},
    	ymin=1e-10,ymax=10,
    	legend style={at={(0.03, 0.03)},anchor=south west,font=\small},
    	legend cell align={left},
    	style={font=\normalsize}]
    	
        \addplot[color=black, 
                style={thick, dotted},
                mark=square,
                mark options=solid]
        table[x=N, y=Linf, col sep=comma]{./figs/data/velocity_convergence_kn0.001_gamma3.csv};
        
        \addplot[color=black, 
                style={thick},
                mark=*,
                mark options=solid,
                mark options={scale=0.8}]
        table[x=N, y=Linf, col sep=comma]{./figs/data/velocity_convergence_kn0.001_gamma3_DVM.csv};
        
        \addplot[color=red!80!black, 
                style={thick, dotted},
                mark=square,
                mark options=solid]
        table[x=N, y=Linf, col sep=comma]{./figs/data/velocity_convergence_kn0.001_gamma14.csv};
        
        \addplot[color=red!80!black, 
                style={thick},
                mark=*,
                mark options=solid,
                mark options={scale=0.6}]
        table[x=N, y=Linf, col sep=comma]{./figs/data/velocity_convergence_kn0.001_gamma14_DVM.csv};

	\end{semilogyaxis}
\end{tikzpicture}

%% file: figs/conservation_kn0p1.tex
\begin{tikzpicture}[spy using outlines={rectangle, height=3cm,width=2.3cm, magnification=3, connect spies}]
	\begin{semilogyaxis}[name=plot1,
		axis line style={latex-latex},
		xlabel={$N_v$},
    	xmin=10,xmax=62,
    	ylabel={$\epsilon_{m}$},
		ylabel style={rotate=-90},
    	ymin=1e-15,ymax=10,
    	legend style={at={(0.97, 0.97)},anchor=north east,font=\small},
    	legend cell align={left},
    	style={font=\normalsize}]
    	
        \addplot[color=black, 
                style={thick, dotted},
                mark=square,
                mark options=solid]
        table[x=N, y=MassErr, col sep=comma]{./figs/data/velocity_convergence_kn0.1_gamma3.csv};
        \addlegendentry{Standard};
        
        \addplot[color=black, 
                style={thick},
                mark=*,
                mark options=solid,
                mark options={scale=0.6}]
        table[x=N, y=MassErr, col sep=comma]{./figs/data/velocity_convergence_kn0.1_gamma3_DVM.csv};
        \addlegendentry{DVM};
        
        \addplot[color=red!80!black, 
                style={thick, dotted},
                mark=square,
                mark options=solid]
        table[x=N, y=MassErr, col sep=comma]{./figs/data/velocity_convergence_kn0.1_gamma14.csv};
        
        \addplot[color=red!80!black, 
                style={thick},
                mark=*,
                mark options=solid,
                mark options={scale=0.6} ]
        table[x=N, y=MassErr, col sep=comma]{./figs/data/velocity_convergence_kn0.1_gamma14_DVM.csv};

	\end{semilogyaxis}
\end{tikzpicture}

%% file: figs/conservation_kn0p01.tex
\begin{tikzpicture}[spy using outlines={rectangle, height=3cm,width=2.3cm, magnification=3, connect spies}]
	\begin{semilogyaxis}[name=plot1,
		axis line style={latex-latex},
		xlabel={$N_v$},
    	xmin=10,xmax=62,
    	ylabel={$\epsilon_{m}$},
		ylabel style={rotate=-90},
    	ymin=1e-15,ymax=10,
    	legend style={at={(0.03, 0.03)},anchor=south west,font=\small},
    	legend cell align={left},
    	style={font=\normalsize}]
    	
        \addplot[color=black, 
                style={thick, dotted},
                mark=square,
                mark options=solid]
        table[x=N, y=MassErr, col sep=comma]{./figs/data/velocity_convergence_kn0.01_gamma3.csv};
        
        \addplot[color=black, 
                style={thick},
                mark=*,
                mark options=solid,
                mark options={scale=0.8}]
        table[x=N, y=MassErr, col sep=comma]{./figs/data/velocity_convergence_kn0.01_gamma3_DVM.csv};
        
        \addplot[color=red!80!black, 
                style={thick, dotted},
                mark=square,
                mark options=solid]
        table[x=N, y=MassErr, col sep=comma]{./figs/data/velocity_convergence_kn0.01_gamma14.csv};
        
        \addplot[color=red!80!black, 
                style={thick},
                mark=*,
                mark options=solid,
                mark options={scale=0.6} ]
        table[x=N, y=MassErr, col sep=comma]{./figs/data/velocity_convergence_kn0.01_gamma14_DVM.csv};

	\end{semilogyaxis}
\end{tikzpicture}

%% file: figs/conservation_kn0p001.tex
\begin{tikzpicture}[spy using outlines={rectangle, height=3cm,width=2.3cm, magnification=3, connect spies}]
	\begin{semilogyaxis}[name=plot1,
		axis line style={latex-latex},
		xlabel={$N_v$},
    	xmin=10,xmax=62,
    	ylabel={$\epsilon_{m}$},
		ylabel style={rotate=-90},
    	ymin=1e-15,ymax=10,
    	legend style={at={(0.03, 0.03)},anchor=south west,font=\small},
    	legend cell align={left},
    	style={font=\normalsize}]
    	
        \addplot[color=black, 
                style={thick, dotted},
                mark=square,
                mark options=solid]
        table[x=N, y=MassErr, col sep=comma]{./figs/data/velocity_convergence_kn0.001_gamma3.csv};
        
        \addplot[color=black, 
                style={thick},
                mark=*,
                mark options=solid,
                mark options={scale=0.8}]
        table[x=N, y=MassErr, col sep=comma]{./figs/data/velocity_convergence_kn0.001_gamma3_DVM.csv};
        
        \addplot[color=red!80!black, 
                style={thick, dotted},
                mark=square,
                mark options=solid]
        table[x=N, y=MassErr, col sep=comma]{./figs/data/velocity_convergence_kn0.001_gamma14.csv};
        
        \addplot[color=red!80!black, 
                style={thick},
                mark=*,
                mark options=solid,
                mark options={scale=0.6}]
        table[x=N, y=MassErr, col sep=comma]{./figs/data/velocity_convergence_kn0.001_gamma14_DVM.csv};

	\end{semilogyaxis}
\end{tikzpicture}

%% file: figs/expansion_density_Kn0p01.tex
\begin{tikzpicture}[spy using outlines={rectangle, height=3cm,width=2.3cm, magnification=3, connect spies}]
	\begin{axis}[name=plot1,
		axis line style={latex-latex},
	    axis x line=left,
        axis y line=left,
		xlabel={$x$},
    	xmin=0, xmax=1,
    	xtick={0, 0.2, 0.4, 0.6, 0.8, 1.0},
    	ylabel={$\rho$},
    	ymin=0,ymax=1.05,
    	ytick={0, 0.2, 0.4, 0.6, 0.8, 1.0},
        clip mode=individual,
    	ylabel style={rotate=-90},
    	legend style={at={(0.5, 0.95)},anchor=north,font=\small},
    	legend cell align={left},
    	style={font=\normalsize}]
    	
        \addplot[color=gray, style={thick}, only marks, mark=o, mark options={scale=0.8}]
        table[x=x, y=r, col sep=comma]{./figs/data/expansion_exact.csv};
        \addlegendentry{Euler};
        
        \addplot[color=black, 
                style={very thick}]
        table[x=x, y=r, col sep=comma]{./figs/data/expansion_p3_N100_Kn0p01.csv};
        \addlegendentry{Standard};
        
        \addplot[color=red, 
                style={ultra thick, dashed}]
        table[x=x, y=r, col sep=comma]{./figs/data/expansion_p3_N100_Kn0p01_DVM.csv};
        \addlegendentry{DVM};

	\end{axis}
\end{tikzpicture}

%% file: figs/expansion_energy_Kn0p01.tex
\begin{tikzpicture}[spy using outlines={rectangle, height=3cm,width=2.3cm, magnification=3, connect spies}]
	\begin{axis}[name=plot1,
		axis line style={latex-latex},
	    axis x line=left,
        axis y line=left,
		xlabel={$x$},
    	xmin=0,xmax=1,
    	xtick={0, 0.2, 0.4, 0.6, 0.8, 1.0},
    	ylabel={$e$},
    	ymin=0,ymax=1.05,
        clip mode=individual,
    	ylabel style={rotate=-90},
    	legend style={at={(1.0, 0.03)},anchor=south east,font=\small},
    	legend cell align={left},
    	style={font=\normalsize}]
    	
        \addplot[color=gray, style={thick}, only marks, mark=o, mark options={scale=0.8}]
        table[x=x, y=e, col sep=comma]{./figs/data/expansion_exact.csv};
        
        \addplot[color=black, 
                style={very thick}]
        table[x=x, y=e, col sep=comma]{./figs/data/expansion_p3_N100_Kn0p01.csv};
        
        \addplot[color=red, 
                style={ultra thick, dashed}]
        table[x=x, y=e, col sep=comma]{./figs/data/expansion_p3_N100_Kn0p01_DVM.csv};

	\end{axis}
\end{tikzpicture}

%% file: figs/expansion_density_Kn0p001.tex
\begin{tikzpicture}[spy using outlines={rectangle, height=3cm,width=2.3cm, magnification=3, connect spies}]
	\begin{axis}[name=plot1,
		axis line style={latex-latex},
	    axis x line=left,
        axis y line=left,
		xlabel={$x$},
    	xmin=0,xmax=1,
    	xtick={0, 0.2, 0.4, 0.6, 0.8, 1.0},
    	ylabel={$\rho$},
    	ymin=0,ymax=1.05,
    	ytick={0, 0.2, 0.4, 0.6, 0.8, 1.0, 1.2},
        clip mode=individual,
    	ylabel style={rotate=-90},
    	legend style={at={(0.5, 0.95)},anchor=north,font=\small},
    	legend cell align={left},
    	style={font=\normalsize}]
    	
        \addplot[color=gray, style={thick}, only marks, mark=o, mark options={scale=0.8}]
        table[x=x, y=r, col sep=comma]{./figs/data/expansion_exact.csv};
        \addlegendentry{Euler};
        
        \addplot[color=black, 
                style={very thick}]
        table[x=x, y=r, col sep=comma]{./figs/data/expansion_p3_N100_Kn0p001.csv};
        \addlegendentry{Standard};
        
        \addplot[color=red, 
                style={ultra thick, dashed}]
        table[x=x, y=r, col sep=comma]{./figs/data/expansion_p3_N100_Kn0p001_DVM.csv};
        \addlegendentry{DVM};

	\end{axis}
\end{tikzpicture}

%% file: figs/expansion_energy_Kn0p001.tex
\begin{tikzpicture}[spy using outlines={rectangle, height=3cm,width=2.3cm, magnification=3, connect spies}]
	\begin{axis}[name=plot1,
		axis line style={latex-latex},
	    axis x line=left,
        axis y line=left,
		xlabel={$x$},
    	xmin=0,xmax=1,
    	xtick={0, 0.2, 0.4, 0.6, 0.8, 1.0},
    	ylabel={$e$},
    	ymin=0,ymax=1.05,
        clip mode=individual,
    	ylabel style={rotate=-90},
    	legend style={at={(1.0, 0.03)},anchor=south east,font=\small},
    	legend cell align={left},
    	style={font=\normalsize}]
    	
        \addplot[color=gray, style={thick}, only marks, mark=o, mark options={scale=0.8}]
        table[x=x, y=e, col sep=comma]{./figs/data/expansion_exact.csv};
        
        \addplot[color=black, 
                style={very thick}]
        table[x=x, y=e, col sep=comma]{./figs/data/expansion_p3_N100_Kn0p001.csv};
        
        \addplot[color=red, 
                style={ultra thick, dashed}]
        table[x=x, y=e, col sep=comma]{./figs/data/expansion_p3_N100_Kn0p001_DVM.csv};

	\end{axis}
\end{tikzpicture}

%% file: figs/expansion_energy_conv_fixed.tex
\begin{tikzpicture}[spy using outlines={rectangle, height=3cm,width=2.3cm, magnification=3, connect spies}]
	\begin{axis}[name=plot1,
		axis line style={latex-latex},
	    axis x line=left,
        axis y line=left,
		xlabel={$x$},
    	xmin=0,xmax=1,
    	xtick={0, 0.2, 0.4, 0.6, 0.8, 1.0},
    	ylabel={$e$},
    	ymin=0,ymax=1.05,
        clip mode=individual,
    	ylabel style={rotate=-90},
    	legend style={at={(1.0, 0.01)},anchor=south east,font=\small},
    	legend cell align={left},
    	legend columns=3,
    	style={font=\normalsize}]
    	
        \addplot[color=gray, style={thick}, only marks, mark=o, mark options={scale=0.8}]
        table[x=x, y=e, col sep=comma]{./figs/data/expansion_exact.csv};
        \addlegendentry{Euler};
        
        \addplot[color=black, 
                style={very thick, dotted}]
        table[x=x, y=e, col sep=comma]{./figs/data/expansion_p3_N100_Kn0p001_DVM.csv};
        \addlegendentry{$N_e = 100$};
        
        \addplot[color=black, 
                style={very thick, dash dot}]
        table[x=x, y=e, col sep=comma]{./figs/data/expansion_p3_N200_Kn0p001_DVM.csv};
        \addlegendentry{$N_e = 200$};
        
        \addplot[color=black, 
                style={very thick, dashed}]
        table[x=x, y=e, col sep=comma]{./figs/data/expansion_p3_N400_Kn0p001_DVM.csv};
        \addlegendentry{$N_e = 400$};
        
        \addplot[color=black, 
                style={very thick}]
        table[x=x, y=e, col sep=comma]{./figs/data/expansion_p3_N800_Kn0p001_DVM.csv};
        \addlegendentry{$N_e = 800$};

	\end{axis}
\end{tikzpicture}

%% file: figs/expansion_energy_conv_decreasing.tex
\begin{tikzpicture}[spy using outlines={rectangle, height=3cm,width=2.3cm, magnification=3, connect spies}]
	\begin{axis}[name=plot1,
		axis line style={latex-latex},
	    axis x line=left,
        axis y line=left,
		xlabel={$x$},
    	xmin=0,xmax=1,
    	xtick={0, 0.2, 0.4, 0.6, 0.8, 1.0},
    	ylabel={$e$},
    	ymin=0,ymax=1.05,
        clip mode=individual,
    	ylabel style={rotate=-90},
    	legend style={at={(1.0, 0.01)},anchor=south east,font=\small},
    	legend cell align={left},
    	legend columns=3,
    	style={font=\normalsize}]
    	
        \addplot[color=gray, style={thick}, only marks, mark=o, mark options={scale=0.8}]
        table[x=x, y=e, col sep=comma]{./figs/data/expansion_exact.csv};
        
        \addplot[color=black, 
                style={very thick, dotted}]
        table[x=x, y=e, col sep=comma]{./figs/data/expansion_p3_N100_Kn0p001_DVM.csv};
        
        \addplot[color=black, 
                style={very thick, dash dot}]
        table[x=x, y=e, col sep=comma]{./figs/data/expansion_p3_N200_Kn0p0005_DVM.csv};
        
        \addplot[color=black, 
                style={very thick, dashed}]
        table[x=x, y=e, col sep=comma]{./figs/data/expansion_p3_N400_Kn0p00025_DVM.csv};
        
        \addplot[color=black, 
                style={very thick}]
        table[x=x, y=e, col sep=comma]{./figs/data/expansion_p3_N800_Kn0p000125_DVM.csv};

	\end{axis}
\end{tikzpicture}

%% file: figs/expansion_energy_conv_sharp.tex
\begin{tikzpicture}[spy using outlines={rectangle, height=3cm,width=2.3cm, magnification=3, connect spies}]
	\begin{axis}[name=plot1,
		axis line style={latex-latex},
	    axis x line=left,
        axis y line=left,
		xlabel={$x$},
    	xmin=0,xmax=1,
    	xtick={0, 0.2, 0.4, 0.6, 0.8, 1.0},
    	ylabel={$e$},
    	ymin=0,ymax=1.05,
        clip mode=individual,
    	ylabel style={rotate=-90},
    	legend style={at={(1.0, 0.01)},anchor=south east,font=\small},
    	legend cell align={left},
    	legend columns=3,
    	style={font=\normalsize}]
    	
        \addplot[color=gray, style={thick}, only marks, mark=o, mark options={scale=0.8}]
        table[x=x, y=e, col sep=comma]{./figs/data/expansion_exact.csv};
        \addlegendentry{Euler};
        
        \addplot[color=black, 
                style={very thick, dotted}]
        table[x=x, y=e, col sep=comma]{./figs/data/expansion_p3_N100_Kn0p001_DVM.csv};
        \addlegendentry{$N_e = 100$};
        
        \addplot[color=black, 
                style={very thick, dash dot}]
        table[x=x, y=e, col sep=comma]{./figs/data/expansion_p3_N200_Kn0p0005_DVM.csv};
        \addlegendentry{$N_e = 200$};
        
        \addplot[color=black, 
                style={very thick, dashed}]
        table[x=x, y=e, col sep=comma]{./figs/data/expansion_p3_N400_Kn0p00025_DVM.csv};
        \addlegendentry{$N_e = 400$};
        
        \addplot[color=black, 
                style={very thick}]
        table[x=x, y=e, col sep=comma]{./figs/data/expansion_p3_N800_Kn0p000125_DVM.csv};
        \addlegendentry{$N_e = 800$};

	\end{axis}
\end{tikzpicture}

%% file: figs/expansion_energy_conv_smooth.tex
\begin{tikzpicture}[spy using outlines={rectangle, height=3cm,width=2.3cm, magnification=3, connect spies}]
	\begin{axis}[name=plot1,
		axis line style={latex-latex},
	    axis x line=left,
        axis y line=left,
		xlabel={$x$},
    	xmin=0,xmax=1,
    	xtick={0, 0.2, 0.4, 0.6, 0.8, 1.0},
    	ylabel={$e$},
    	ymin=0,ymax=1.05,
        clip mode=individual,
    	ylabel style={rotate=-90},
    	legend style={at={(1.0, 0.01)},anchor=south east,font=\small},
    	legend cell align={left},
    	legend columns=3,
    	style={font=\normalsize}]
    	
        \addplot[color=gray, style={thick}, only marks, mark=o, mark options={scale=0.8}]
        table[x=x, y=e, col sep=comma]{./figs/data/expansion_exact.csv};
        
        \addplot[color=black, 
                style={very thick, dotted}]
        table[x=x, y=e, col sep=comma]{./figs/data/expansion_p3_N100_Kn0p001_DVM_smooth.csv};
        
        \addplot[color=black, 
                style={very thick, dash dot}]
        table[x=x, y=e, col sep=comma]{./figs/data/expansion_p3_N200_Kn0p0005_DVM_smooth.csv};
        
        \addplot[color=black, 
                style={very thick, dashed}]
        table[x=x, y=e, col sep=comma]{./figs/data/expansion_p3_N400_Kn0p00025_DVM_smooth.csv};
        
        \addplot[color=black, 
                style={very thick}]
        table[x=x, y=e, col sep=comma]{./figs/data/expansion_p3_N800_Kn0p000125_DVM_smooth.csv};

	\end{axis}
\end{tikzpicture}

%% file: figs/normalshock_itr.tex
\begin{tikzpicture}[spy using outlines={rectangle, height=3cm,width=2.3cm, magnification=3, connect spies}]
	\begin{axis}[name=plot1,
		axis line style={latex-latex},
	    axis x line=left,
        axis y line=left,
		xlabel={$M$},
    	xmin=1, xmax=10.5,
    	xtick={1,2,3,4,5,6,7,8,9,10,11},
    	ylabel={$\lambda_L/\Delta$},
    	ymin=0,ymax=0.52,
    	ytick={0, 0.1, 0.2, 0.3, 0.4, 0.5},
        clip mode=individual,
    	legend style={at={(0.97, 0.03)},anchor=south east,font=\small},
    	legend cell align={left},
    	style={font=\normalsize}    	]

        \addplot[color=black, style={very thick, dotted}]
        table[x=M, y=y, col sep=comma]{./figs/data/normalshock_itr_navierstokes.csv};
        \addlegendentry{Navier--Stokes \citep{Mieussens2000}};
        
         
        \addplot[color=black, 
                style={thick},
                mark=o,
                mark options=solid,
                mark options={scale=0.8}]
        table[x=M, y=y, col sep=comma]{./figs/data/normalshock_itr_bgk.csv} ;
        \addlegendentry{Boltzmann--BGK};

        \addplot[color=black!60, style={thick}, only marks, mark=square, mark options={scale=0.7}]
        table[x=M, y=y, col sep=comma]{./figs/data/normalshock_itr_alsmeyer.csv};
        
        \addplot[color=black!60, style={thick}, only marks, mark=diamond, mark options={scale=1.1}]
        table[x=M, y=y, col sep=comma]{./figs/data/normalshock_itr_camac.csv};
        
        \addplot[color=black!60, style={thick}, only marks, mark=triangle, mark options={scale=1.1}]
        table[x=M, y=y, col sep=comma]{./figs/data/normalshock_itr_linzer.csv};

	\end{axis}
\end{tikzpicture}

%% file: figs/normalshock_profile_m3p8.tex
\begin{tikzpicture}[spy using outlines={rectangle, height=3cm,width=2.3cm, magnification=3, connect spies}]
	\begin{axis}[name=plot1,
		axis line style={latex-latex},
	    axis x line=left,
        axis y line=left,
		xlabel={$x/\lambda_L$},
    	xmin=-10, xmax=10,
    	ylabel={$\rho^*$, $U^*$, $P^*$},
    	ymin=-.05,ymax=1.05,
    	ytick={0, 0.2, 0.4, 0.6, 0.8, 1.0},
        clip mode=individual,
    	legend style={at={(0.97, 0.1)},anchor=south east,font=\small},
    	legend cell align={left},
    	style={font=\normalsize}]

        \addplot[color=black, 
                style={very thick}]
        table[x=x, y=r, col sep=comma]{./figs/data/normalshock_density_m3p8_bgk.csv} ;
        \addlegendentry{Density};
        
        \addplot[color=black, 
                style={very thick, dotted}]
        table[x=x, y=u, col sep=comma]{./figs/data/normalshock_density_m3p8_bgk.csv} ;
        \addlegendentry{Velocity};
        
        \addplot[color=black, 
                style={very thick, dashed}]
        table[x=x, y=p, col sep=comma]{./figs/data/normalshock_density_m3p8_bgk.csv} ;
        \addlegendentry{Pressure};
        
        \addplot[color=black!80, style={thick}, only marks, mark=square, mark options={scale=0.8}]
        table[x=x, y=r, col sep=comma]{./figs/data/normalshock_density_m3p8_alsmeyer.csv};
        
        \addplot[color=black!80, style={thick}, only marks, mark=triangle, mark options={scale=1.0}]
        table[x=x, y=r, col sep=comma]{./figs/data/normalshock_density_m3p8_bentley.csv};

	\end{axis}
\end{tikzpicture}

%% file: figs/normalshock_profile_m9p0.tex
\begin{tikzpicture}[spy using outlines={rectangle, height=3cm,width=2.3cm, magnification=3, connect spies}]
	\begin{axis}[name=plot1,
		axis line style={latex-latex},
	    axis x line=left,
        axis y line=left,
		xlabel={$x/\lambda_L$},
    	xmin=-10, xmax=10,
    	ylabel={$\rho^*$, $U^*$, $P^*$},
    	ymin=-.05,ymax=1.05,
    	ytick={0, 0.2, 0.4, 0.6, 0.8, 1.0},
        clip mode=individual,
    	legend style={at={(0.97, 0.1)},anchor=south east,font=\small},
    	legend cell align={left},
    	style={font=\normalsize}]

        \addplot[color=black, 
                style={very thick}]
        table[x=x, y=r, col sep=comma]{./figs/data/normalshock_density_m9p0_bgk.csv} ;
        
        \addplot[color=black, 
                style={very thick, dotted}]
        table[x=x, y=u, col sep=comma]{./figs/data/normalshock_density_m9p0_bgk.csv} ;
        
        \addplot[color=black, 
                style={very thick, dashed}]
        table[x=x, y=p, col sep=comma]{./figs/data/normalshock_density_m9p0_bgk.csv} ;
        
        \addplot[color=black!80, style={thick}, only marks, mark=square, mark options={scale=0.8}]
        table[x=x, y=r, col sep=comma]{./figs/data/normalshock_density_m9p0_alsmeyer.csv};
        
        \addplot[color=black!80, style={thick}, only marks, mark=triangle, mark options={scale=1.0}]
        table[x=x, y=r, col sep=comma]{./figs/data/normalshock_density_m9p0_bentley.csv};

	\end{axis}
\end{tikzpicture}

%% file: figs/normalshock_pdf_m3p8.tex
\begin{tikzpicture}[spy using outlines={rectangle, height=3cm,width=2.3cm, magnification=3, connect spies}]
	\begin{axis}[name=plot1,
		axis line style={latex-latex},
	    axis x line=left,
        axis y line=left,
		xlabel={$u$},
    	xmin=-20, xmax=20,
    	ylabel={$f_u (u)$},
    	ymin=-.025,ymax=0.5,
        clip mode=individual,
    	legend style={at={(0.03, 0.97)},anchor=north west,font=\small},
    	legend cell align={left},
    	style={font=\normalsize}]
    	
        \addlegendimage{color=black, style={thin}, mark=o, mark options={solid, thin, scale=0.9}}
        \addlegendentry{Upstream};
        \addlegendimage{color=black, style={thin, dashed}, mark=square, mark options={solid, thin, scale=0.8}}
        \addlegendentry{Shock};
        \addlegendimage{color=black, style={thick, dotted}, mark=triangle, mark options={solid, thin, scale=1}}
        \addlegendentry{Downstream};
        
        \addplot[color=black!80, style={thin}, only marks, mark=o, mark options={scale=0.9}]
        table[x=u, y=fl, col sep=comma]{./figs/data/normalshock_pdf_m3p8_bgk.csv};
        
        \addplot[color=black!80, style={thin}, only marks, mark=square, mark options={scale=0.8}]
        table[x=u, y=fm, col sep=comma]{./figs/data/normalshock_pdf_m3p8_bgk.csv};
        
        \addplot[color=black!80, style={thin}, only marks, mark=triangle, mark options={scale=1}]
        table[x=u, y=fr, col sep=comma]{./figs/data/normalshock_pdf_m3p8_bgk.csv};
        
        \addplot[color=black, style={thin}]
        table[x=u, y=fl, col sep=comma]{./figs/data/normalshock_pdf_m3p8_bgk_splined.csv};
        
        \addplot[color=black, style={thin, dashed}]
        table[x=u, y=fm, col sep=comma]{./figs/data/normalshock_pdf_m3p8_bgk_splined.csv};
        
        \addplot[color=black, style={thick, dotted}]
        table[x=u, y=fr, col sep=comma]{./figs/data/normalshock_pdf_m3p8_bgk_splined.csv};

	\end{axis}
\end{tikzpicture}

%% file: figs/normalshock_pdf_m9p0.tex
\begin{tikzpicture}[spy using outlines={rectangle, height=3cm,width=2.3cm, magnification=3, connect spies}]
	\begin{axis}[name=plot1,
		axis line style={latex-latex},
	    axis x line=left,
        axis y line=left,
		xlabel={$u$},
    	xmin=-30, xmax=30,
    	y tick label style={
    		/pgf/number format/.cd,
        	fixed,
        	precision=2,
    	    /tikz/.cd},
    	ylabel={$f_u (u)$},
    	ymin=-.01,ymax=0.2,
    	ytick={0, 0.05, 0.1, 0.15, 0.2},
        clip mode=individual,
    	legend style={at={(0.03, 0.97)},anchor=north west,font=\small},
    	legend cell align={left},
    	style={font=\normalsize}]

        \addplot[color=black!80, style={thin}, only marks, mark=o, mark options={scale=0.9}]
        table[x=u, y=fl, col sep=comma]{./figs/data/normalshock_pdf_m9p0_bgk.csv};
        
        \addplot[color=black!80, style={thin}, only marks, mark=square, mark options={scale=0.8}]
        table[x=u, y=fm, col sep=comma]{./figs/data/normalshock_pdf_m9p0_bgk.csv};
        
        \addplot[color=black!80, style={thin}, only marks, mark=triangle, mark options={scale=1}]
        table[x=u, y=fr, col sep=comma]{./figs/data/normalshock_pdf_m9p0_bgk.csv};
        
        \addplot[color=black, style={thin}]
        table[x=u, y=fl, col sep=comma]{./figs/data/normalshock_pdf_m9p0_bgk_splined.csv};
        
        \addplot[color=black, style={thin, dashed}]
        table[x=u, y=fm, col sep=comma]{./figs/data/normalshock_pdf_m9p0_bgk_splined.csv};
        
        \addplot[color=black, style={thick, dotted}]
        table[x=u, y=fr, col sep=comma]{./figs/data/normalshock_pdf_m9p0_bgk_splined.csv};

	\end{axis}
\end{tikzpicture}

%% file: figs/sod_density_kn1.tex
\begin{tikzpicture}[spy using outlines={rectangle, height=3cm,width=2.3cm, magnification=3, connect spies}]
	\begin{axis}[name=plot1,
		axis line style={latex-latex},
	    axis x line=left,
        axis y line=left,
		xlabel={$x$},
    	xmin=0, xmax=1,
    	xtick={0, 0.2, 0.4, 0.6, 0.8, 1.0},
    	ylabel={$\rho$,$U$,$p$},
    	ymin=0,ymax=1.1,
    	ytick={0, 0.2, 0.4, 0.6, 0.8, 1.0},
        clip mode=individual,
    	legend style={at={(0.01, 0.45)},anchor=west,font=\small},
    	legend cell align={left},
    	style={font=\normalsize}]
    	
        \addplot[color=gray, style={thick}, only marks, mark=o, mark options={scale=0.5}]
        table[x=x, y=r, col sep=comma]{./figs/data/sod_exact.csv};
        \addlegendentry{Euler};
        
        \addplot[color=gray, style={thick}, only marks, mark=o, mark options={scale=0.5}, forget plot]
        table[x=x, y=u, col sep=comma]{./figs/data/sod_exact.csv};
        
        \addplot[color=gray, style={thick}, only marks, mark=o, mark options={scale=0.5}, forget plot]
        table[x=x, y=p, col sep=comma]{./figs/data/sod_exact.csv};
         
        \addplot[color=black, 
                style={very thick}]
        table[x=x, y=r, col sep=comma]{./figs/data/sod_Ne50_Nv16_Kncell1.csv};
        \addlegendentry{Density};
        
        \addplot[color=black, 
                style={very thick, dotted}]
        table[x=x, y=u, col sep=comma]{./figs/data/sod_Ne50_Nv16_Kncell1.csv};
        \addlegendentry{Velocity};
        
        \addplot[color=black, 
                style={very thick, dashed}]
        table[x=x, y=p, col sep=comma]{./figs/data/sod_Ne50_Nv16_Kncell1.csv};
        \addlegendentry{Pressure};

	\end{axis}
\end{tikzpicture}

%% file: figs/sod_density_kn0p1.tex
\begin{tikzpicture}[spy using outlines={rectangle, height=3cm,width=2.3cm, magnification=3, connect spies}]
	\begin{axis}[name=plot1,
		axis line style={latex-latex},
	    axis x line=left,
        axis y line=left,
		xlabel={$x$},
    	xmin=0, xmax=1,
    	xtick={0, 0.2, 0.4, 0.6, 0.8, 1.0},
    	ylabel={$\rho$,$U$,$p$},
    	ymin=0,ymax=1.1,
    	ytick={0, 0.2, 0.4, 0.6, 0.8, 1.0},
        clip mode=individual,
    	legend style={at={(0.01, 0.45)},anchor=west,font=\small},
    	legend cell align={left},
    	style={font=\normalsize}]
    	
        \addplot[color=gray, style={thick}, only marks, mark=o, mark options={scale=0.5}]
        table[x=x, y=r, col sep=comma]{./figs/data/sod_exact.csv};
        
        \addplot[color=gray, style={thick}, only marks, mark=o, mark options={scale=0.5}, forget plot]
        table[x=x, y=u, col sep=comma]{./figs/data/sod_exact.csv};
        
        \addplot[color=gray, style={thick}, only marks, mark=o, mark options={scale=0.5}, forget plot]
        table[x=x, y=p, col sep=comma]{./figs/data/sod_exact.csv};
        
        \addplot[color=black, 
                style={very thick}]
        table[x=x, y=r, col sep=comma]{./figs/data/sod_Ne50_Nv16_Kncell0.1.csv};
        
        \addplot[color=black, 
                style={very thick, dotted}]
        table[x=x, y=u, col sep=comma]{./figs/data/sod_Ne50_Nv16_Kncell0.1.csv};
        
        \addplot[color=black, 
                style={very thick, dashed}]
        table[x=x, y=p, col sep=comma]{./figs/data/sod_Ne50_Nv16_Kncell0.1.csv};

	\end{axis}
\end{tikzpicture}

%% file: figs/sod_density_kn0p01.tex
\begin{tikzpicture}[spy using outlines={rectangle, height=3cm,width=2.3cm, magnification=3, connect spies}]
	\begin{axis}[name=plot1,
		axis line style={latex-latex},
	    axis x line=left,
        axis y line=left,
		xlabel={$x$},
    	xmin=0, xmax=1,
    	xtick={0, 0.2, 0.4, 0.6, 0.8, 1.0},
    	ylabel={$\rho$,$U$,$p$},
    	ymin=0,ymax=1.1,
    	ytick={0, 0.2, 0.4, 0.6, 0.8, 1.0},
        clip mode=individual,
    	legend style={at={(0.01, 0.45)},anchor=west,font=\small},
    	legend cell align={left},
    	style={font=\normalsize}]
    	
        \addplot[color=gray, style={thick}, only marks, mark=o, mark options={scale=0.5}]
        table[x=x, y=r, col sep=comma]{./figs/data/sod_exact.csv};
        
        \addplot[color=gray, style={thick}, only marks, mark=o, mark options={scale=0.5}, forget plot]
        table[x=x, y=u, col sep=comma]{./figs/data/sod_exact.csv};
        
        \addplot[color=gray, style={thick}, only marks, mark=o, mark options={scale=0.5}, forget plot]
        table[x=x, y=p, col sep=comma]{./figs/data/sod_exact.csv};
        
        \addplot[color=black, 
                style={very thick}]
        table[x=x, y=r, col sep=comma]{./figs/data/sod_Ne50_Nv16_Kncell0.01.csv};
        
        \addplot[color=black, 
                style={very thick, dotted}]
        table[x=x, y=u, col sep=comma]{./figs/data/sod_Ne50_Nv16_Kncell0.01.csv};
        
        \addplot[color=black, 
                style={very thick, dashed}]
        table[x=x, y=p, col sep=comma]{./figs/data/sod_Ne50_Nv16_Kncell0.01.csv};

	\end{axis}
\end{tikzpicture}

%% file: figs/sod_density_50.tex
\begin{tikzpicture}[spy using outlines={rectangle, height=3cm,width=2.3cm, magnification=3, connect spies}]
	\begin{axis}[name=plot1,
		axis line style={latex-latex},
	    axis x line=left,
        axis y line=left,
		xlabel={$x$},
    	xmin=0, xmax=1,
    	xtick={0, 0.2, 0.4, 0.6, 0.8, 1.0},
    	ylabel={$\rho$,$U$,$p$},
    	ymin=0,ymax=1.1,
    	ytick={0, 0.2, 0.4, 0.6, 0.8, 1.0},
        clip mode=individual,
    	legend style={at={(0.01, 0.45)},anchor=west,font=\small},
    	legend cell align={left},
    	style={font=\normalsize}]
    	
        \addplot[color=gray, style={thick}, only marks, mark=o, mark options={scale=0.5}]
        table[x=x, y=r, col sep=comma]{./figs/data/sod_exact.csv};
        \addlegendentry{Euler};
        
        \addplot[color=gray, style={thick}, only marks, mark=o, mark options={scale=0.5}, forget plot]
        table[x=x, y=u, col sep=comma]{./figs/data/sod_exact.csv};
        
        \addplot[color=gray, style={thick}, only marks, mark=o, mark options={scale=0.5}, forget plot]
        table[x=x, y=p, col sep=comma]{./figs/data/sod_exact.csv};

        \addplot[color=black, 
                style={very thick}]
        table[x=x, y=r, col sep=comma]{./figs/data/sod_Ne50_Nv16_Kncell0.1.csv};
        \addlegendentry{Density};
        
        \addplot[color=black, 
                style={very thick, dotted}]
        table[x=x, y=u, col sep=comma]{./figs/data/sod_Ne50_Nv16_Kncell0.1.csv};
        \addlegendentry{Velocity};
        
        \addplot[color=black, 
                style={very thick, dashed}]
        table[x=x, y=p, col sep=comma]{./figs/data/sod_Ne50_Nv16_Kncell0.1.csv};
        \addlegendentry{Pressure};

	\end{axis}
\end{tikzpicture}

%% file: figs/sod_density_100.tex
\begin{tikzpicture}[spy using outlines={rectangle, height=3cm,width=2.3cm, magnification=3, connect spies}]
	\begin{axis}[name=plot1,
		axis line style={latex-latex},
	    axis x line=left,
        axis y line=left,
		xlabel={$x$},
    	xmin=0, xmax=1,
    	xtick={0, 0.2, 0.4, 0.6, 0.8, 1.0},
    	ylabel={$\rho$,$U$,$p$},
    	ymin=0,ymax=1.1,
    	ytick={0, 0.2, 0.4, 0.6, 0.8, 1.0},
        clip mode=individual,
    	legend style={at={(0.01, 0.45)},anchor=west,font=\small},
    	legend cell align={left},
    	style={font=\normalsize}]
    	
        \addplot[color=gray, style={thick}, only marks, mark=o, mark options={scale=0.5}]
        table[x=x, y=r, col sep=comma]{./figs/data/sod_exact.csv};
        
        \addplot[color=gray, style={thick}, only marks, mark=o, mark options={scale=0.5}, forget plot]
        table[x=x, y=u, col sep=comma]{./figs/data/sod_exact.csv};
        
        \addplot[color=gray, style={thick}, only marks, mark=o, mark options={scale=0.5}, forget plot]
        table[x=x, y=p, col sep=comma]{./figs/data/sod_exact.csv};
        
        \addplot[color=black, 
                style={very thick}]
        table[x=x, y=r, col sep=comma]{./figs/data/sod_Ne100_Nv16_Kncell0.1.csv};
        
        \addplot[color=black, 
                style={very thick, dotted}]
        table[x=x, y=u, col sep=comma]{./figs/data/sod_Ne100_Nv16_Kncell0.1.csv};
        
        \addplot[color=black, 
                style={very thick, dashed}]
        table[x=x, y=p, col sep=comma]{./figs/data/sod_Ne100_Nv16_Kncell0.1.csv};

	\end{axis}
\end{tikzpicture}

%% file: figs/sod_density_200.tex
\begin{tikzpicture}[spy using outlines={rectangle, height=3cm,width=2.3cm, magnification=3, connect spies}]
	\begin{axis}[name=plot1,
		axis line style={latex-latex},
	    axis x line=left,
        axis y line=left,
		xlabel={$x$},
    	xmin=0, xmax=1,
    	xtick={0, 0.2, 0.4, 0.6, 0.8, 1.0},
    	ylabel={$\rho$,$U$,$p$},
    	ymin=0,ymax=1.1,
    	ytick={0, 0.2, 0.4, 0.6, 0.8, 1.0},
        clip mode=individual,
    	legend style={at={(0.01, 0.45)},anchor=west,font=\small},
    	legend cell align={left},
    	style={font=\normalsize}]
    	
        \addplot[color=gray, style={thick}, only marks, mark=o, mark options={scale=0.5}]
        table[x=x, y=r, col sep=comma]{./figs/data/sod_exact.csv};
        \addplot[color=gray, style={thick}, only marks, mark=o, mark options={scale=0.5}, forget plot]
        table[x=x, y=u, col sep=comma]{./figs/data/sod_exact.csv};
        
        \addplot[color=gray, style={thick}, only marks, mark=o, mark options={scale=0.5}, forget plot]
        table[x=x, y=p, col sep=comma]{./figs/data/sod_exact.csv};

        \addplot[color=black, 
                style={very thick}]
        table[x=x, y=r, col sep=comma]{./figs/data/sod_Ne200_Nv16_Kncell0.1.csv};
        
        \addplot[color=black, 
                style={very thick, dotted}]
        table[x=x, y=u, col sep=comma]{./figs/data/sod_Ne200_Nv16_Kncell0.1.csv};
        
        \addplot[color=black, 
                style={very thick, dashed}]
        table[x=x, y=p, col sep=comma]{./figs/data/sod_Ne200_Nv16_Kncell0.1.csv};

	\end{axis}
\end{tikzpicture}

%% file: figs/kh_enst.tex
\begin{tikzpicture}[spy using outlines={rectangle, height=3cm,width=2.5cm, magnification=3, connect spies}]
    \begin{semilogyaxis}
    [
        axis line style={latex-latex},
        axis y line=left,
        axis x line=left,
        clip mode=individual,
        xlabel = {$t$},
        ylabel = {$\varepsilon_E$},
        xmin = 0, xmax = 5,
        ymin = 10^1, ymax = 10^4,
        legend cell align={left},
        legend style={font=\scriptsize, at={(1.0, 1.0)}, anchor=north east},
        x tick label style={/pgf/number format/.cd, fixed, fixed zerofill, precision=0, /tikz/.cd},
    ]
        
        \addplot[color=black, style={ultra thick, dotted}] table[x=t, y=enst, col sep=comma]{./figs/data/enst_ke_ns.csv};
        \addlegendentry{Navier--Stokes};
        
        \addplot[color=black, style={thick}] table[x=t, y=enst, col sep=comma ]{./figs/data/enst_ke_bgk.csv};
        \addlegendentry{Boltzmann--BGK};
        
    \end{semilogyaxis}

\end{tikzpicture}

%% file: figs/tgv_64_enst.tex
\begin{tikzpicture}[spy using outlines={rectangle, height=3cm,width=2.5cm, magnification=3, connect spies}]
    \begin{axis}
    [
        axis line style={latex-latex},
        axis y line=left,
        axis x line=left,
        clip mode=individual,
        xmode=linear, 
        ymode=linear,
        xlabel = {$t$},
        ylabel = {$\varepsilon_D$},
        xmin = 0, xmax = 20,
        ymin = 0.00, ymax = 0.015,
        legend cell align={left},
        legend style={font=\scriptsize, at={(1.0, 1.0)}, anchor=north east},
        ytick = {0,0.004,0.008,0.012,0.016},
        x tick label style={/pgf/number format/.cd, fixed, fixed zerofill, precision=0, /tikz/.cd},
        y tick label style={/pgf/number format/.cd, fixed, fixed zerofill, precision=1, /tikz/.cd},	
        scale = 0.9
    ]
        \addplot[ color=gray, style={thick}, only marks, mark=o, mark options={scale=0.8}, mark repeat = 3, mark phase = 0] table[x=t, y=enst, col sep=comma, mark=*]{./figs/data/tgvNS_256_post.csv};
        \addlegendentry{Reference};
        
        \addplot[color=black, style={ultra thick, dotted}] table[x=t, y=enst, col sep=comma]{./figs/data/tgvNS_64_post.csv};
        \addlegendentry{Navier--Stokes};
        
        \addplot[color=black, style={thick}] table[x=t, y=enst, col sep=comma ]{./figs/data/tgvBGK_64_post.csv};
        \addlegendentry{Boltzmann--BGK };
        
    \end{axis}

\end{tikzpicture}

%% file: figs/tgv_128_enst.tex
\begin{tikzpicture}[spy using outlines={rectangle, height=3cm,width=2.5cm, magnification=3, connect spies}]
    \begin{axis}
    [
        axis line style={latex-latex},
        axis y line=left,
        axis x line=left,
        clip mode=individual,
        xmode=linear, 
        ymode=linear,
        xlabel = {$t$},
        ylabel = {$\varepsilon_D$},
        xmin = 0, xmax = 20,
        ymin = 0.00, ymax = 0.015,
        legend cell align={left},
        legend style={font=\scriptsize, at={(1.0, 1.0)}, anchor=north east},
        ytick = {0,0.004,0.008,0.012,0.016},
        x tick label style={/pgf/number format/.cd, fixed, fixed zerofill, precision=0, /tikz/.cd},
        y tick label style={/pgf/number format/.cd, fixed, fixed zerofill, precision=1, /tikz/.cd},	
        scale = 0.9,
    	legend cell align={left},
    	legend columns=3,
    ]

        \addplot[ color=gray, style={thick}, only marks, mark=o, mark options={scale=0.8}, mark repeat = 3, mark phase = 0] table[x=t, y=enst, col sep=comma, mark=* ]{./figs/data/tgvNS_256_post.csv};
        
        \addplot[color=black, style={ultra thick, dotted}] table[x=t, y=enst, col sep=comma]{./figs/data/tgvNS_128_post.csv  };
        
        \addplot[color=black, style={thick}] table[x=t, y=enst, col sep=comma]{./figs/data/tgvBGK_128_post.csv};
        
    \end{axis}

\end{tikzpicture}

%% file: figs/tgv_192_enst.tex
\begin{tikzpicture}[spy using outlines={rectangle, height=3cm,width=2.5cm, magnification=3, connect spies}]
    \begin{axis}
    [
        axis line style={latex-latex},
        axis y line=left,
        axis x line=left,
        clip mode=individual,
        xmode=linear, 
        ymode=linear,
        xlabel = {$t$},
        ylabel = {$\varepsilon_D$},
        xmin = 0, xmax = 20,
        ymin = 0.00, ymax = 0.015,
        legend cell align={left},
        legend style={font=\scriptsize, at={(1.0, 1.0)}, anchor=north east},
        ytick = {0,0.004,0.008,0.012,0.016},
        x tick label style={/pgf/number format/.cd, fixed, fixed zerofill, precision=0, /tikz/.cd},
        y tick label style={/pgf/number format/.cd, fixed, fixed zerofill, precision=1, /tikz/.cd},	
        scale = 0.9
    ]
        
        \addplot[ color=gray, style={thick}, only marks, mark=o, mark options={scale=0.8}, mark repeat = 3, mark phase = 0] table[x=t, y=enst, col sep=comma, mark=*]{./figs/data/tgvNS_256_post.csv};
        
        \addplot[color=black, style={ultra thick, dotted}] table[x=t, y=enst, col sep=comma]{./figs/data/tgvNS_192_post.csv  };
        
        \addplot[color=black, style={thick}] table[x=t, y=enst, col sep=comma]{./figs/data/tgvBGK_192_post.csv};
    \end{axis}

\end{tikzpicture}

%% file: figs/tgv_128_velconv.tex
\begin{tikzpicture}[spy using outlines={rectangle, height=3cm,width=2.5cm, magnification=3, connect spies}]
    \begin{axis}
    [
        axis line style={latex-latex},
        axis y line=left,
        axis x line=left,
        clip mode=individual,
        xmode=linear, 
        ymode=linear,
        xlabel = {$t$},
        ylabel = {$\varepsilon_D$},
        xmin = 0, xmax = 20,
        ymin = 0.00, ymax = 0.015,
        legend cell align={left},
        legend style={font=\scriptsize, at={(1.0, 1.0)}, anchor=north east},
        ytick = {0,0.004,0.008,0.012,0.016},
        x tick label style={/pgf/number format/.cd, fixed, fixed zerofill, precision=0, /tikz/.cd},
        y tick label style={/pgf/number format/.cd, fixed, fixed zerofill, precision=1, /tikz/.cd},	
        scale = 0.9,
    	legend cell align={left},
    ]

        \addplot[ color=gray, style={thick}, only marks, mark=o, mark options={scale=0.8}, mark repeat = 3, mark phase = 0] table[x=t, y=enst, col sep=comma, mark=* ]{./figs/data/tgvNS_256_post.csv};
        
        \addplot[color=black, style={very thick, dashed}] table[x=t, y=enst, col sep=comma]{./figs/data/tgvBGK_128lowv_post.csv};
        
        \addplot[color=black, style={very thick}] table[x=t, y=enst, col sep=comma]{./figs/data/tgvBGK_128highv_post.csv};
        
        \addplot[color=red, style={ultra thick, dotted}] table[x=t, y=enst, col sep=comma]{./figs/data/tgvBGK_128_post.csv};
        
    \end{axis}

\end{tikzpicture}

%% file: figs/tgv_128_conservation.tex
\begin{tikzpicture}[spy using outlines={rectangle, height=3cm,width=2.5cm, magnification=3, connect spies}]
    \begin{axis}
    [
        axis line style={latex-latex},
        axis y line=left,
        axis x line=left,
        clip mode=individual,
        xmode=linear, 
        ymode=log,
        xlabel = {$t$},
        ylabel = {$\epsilon_{m}$},
        xmin = 0, xmax = 20,
        ymin = 1e-14, ymax = 1e-6,
        legend cell align={left},
    	legend style={at={(0.97, 0.97)},anchor=north east,font=\small},
        x tick label style={/pgf/number format/.cd, fixed, fixed zerofill, precision=0, /tikz/.cd},
        scale = 0.9,
    ]
        
        \addlegendimage{only marks, color=gray, style={very thick}, mark=o, mark options={scale=0.8}}
        \addlegendentry{Reference};

        \addplot[color=black, style={very thick, dashed}] table[x=t, y=masserror, col sep=comma]{./figs/data/tgvBGK_128lowv_post.csv};
        \addlegendentry{$N_v = 12 \times 12 \times 6$};
        
        \addplot[color=red, style={ultra thick, dotted}] table[x=t, y expr={\thisrow{masserror}*1.0}, col sep=comma]{./figs/data/tgvBGK_128_post.csv};
        \addlegendentry{$N_v = 16 \times 16 \times 8$};
        
        \addplot[color=black, style={very thick}] table[x=t, y expr={\thisrow{masserror}*1.0}, col sep=comma]{./figs/data/tgvBGK_128highv_post.csv} ;
        \addlegendentry{$N_v = 20 \times 20 \times 10$};
        
    \end{axis}

\end{tikzpicture}

%% file: figs/tgv_64_spectra.tex
\begin{tikzpicture}[spy using outlines={rectangle, height=3cm,width=2.3cm, magnification=3, connect spies}]
	\begin{loglogaxis}[name=plot1,
		xlabel={$k$},
		xmin=1,xmax=256,
		ylabel={$E(k)$},
		ymin=5e-5,ymax=1e-2,
		legend style={at={(0.03,0.03)},anchor=south west,font=\small},
		legend cell align={left},
		style={font=\normalsize},
        ]

        \addplot[ color=gray, style={ultra thick, dashed}] table[x=k, y=E, col sep=comma, mark=* ]{./figs/data/tgvNS_256_spectra.csv};
        \addlegendentry{Reference};
        
        \addplot[color=black, style={ultra thick, dotted}] table[x=k, y=E, col sep=comma]{./figs/data/tgvNS_64_spectra.csv};
        \addlegendentry{Navier--Stokes};
        
        \addplot[color=black, style={thick}] table[x=k, y=E, col sep=comma]{./figs/data/tgvBGK_64_spectra.csv};
        \addlegendentry{Boltzmann--BGK};
			
		\addplot[color=gray, style={solid, thick},forget plot] coordinates{(50, 1e-2) (256, 0.00065748128)};
	    \node [below,color=black] at (axis cs:100,.008) {$-\frac{5}{3}$};
		\end{loglogaxis} 		
	\end{tikzpicture}

%% file: figs/tgv_128_spectra.tex
\begin{tikzpicture}[spy using outlines={rectangle, height=3cm,width=2.3cm, magnification=3, connect spies}]
	\begin{loglogaxis}[name=plot1,
		xlabel={$k$},
		xmin=1,xmax=256,
		ylabel={$E(k)$},
		ymin=5e-5,ymax=1e-2,
		legend style={at={(0.03,0.03)},anchor=south west,font=\small},
		legend cell align={left},
		style={font=\normalsize},
        ]

        \addplot[ color=gray, style={ultra thick, dashed}] table[x=k, y=E, col sep=comma, mark=* ]{./figs/data/tgvNS_256_spectra.csv};
        
        \addplot[color=black, style={ultra thick, dotted}] table[x=k, y=E, col sep=comma]{./figs/data/tgvNS_128_spectra.csv};
        
        \addplot[color=black, style={thick}] table[x=k, y=E, col sep=comma]{./figs/data/tgvBGK_128_spectra.csv};
			
		\addplot[color=gray, style={solid, thick},forget plot] coordinates{(50, 1e-2) (256, 0.00065748128)};
	    \node [below,color=black] at (axis cs:100,.008) {$-\frac{5}{3}$};
		\end{loglogaxis} 		
	\end{tikzpicture}

%% file: figs/tgv_192_spectra.tex
\begin{tikzpicture}[spy using outlines={rectangle, height=3cm,width=2.3cm, magnification=3, connect spies}]
	\begin{loglogaxis}[name=plot1,
		xlabel={$k$},
		xmin=1,xmax=256,
		ylabel={$E(k)$},
		ymin=5e-5,ymax=1e-2,
		legend style={at={(0.03,0.03)},anchor=south west,font=\small},
		legend cell align={left},
		style={font=\normalsize},
        ]

        \addplot[ color=gray, style={ultra thick, dashed}] table[x=k, y=E, col sep=comma, mark=* ]{./figs/data/tgvNS_256_spectra.csv};
        
        \addplot[color=black, style={ultra thick, dotted}] table[x=k, y=E, col sep=comma]{./figs/data/tgvNS_192_spectra.csv};
        
        \addplot[color=black, style={thick}] table[x=k, y=E, col sep=comma]{./figs/data/tgvBGK_192_spectra.csv};
			
		\addplot[color=gray, style={solid, thick},forget plot] coordinates{(50, 1e-2) (256, 0.00065748128)};
	    \node [below,color=black] at (axis cs:100,.008) {$-\frac{5}{3}$};
		\end{loglogaxis} 		
	\end{tikzpicture}

%% file: conclusions.tex
\section{Conclusions}\label{sec:conclusion}

In this work, we have presented a novel numerical approach for the solving the polyatomic Boltzmann equation with a BGK collision model. The proposed approach utilizes a combination of an efficient implementation of a high-order flux reconstruction spatial discretization augmented with a positivity-preserving limiter and a discrete velocity model ensuring the conservation and well-balancing properties of the scheme. As a result, the computation of three-dimensional problems on arbitrary domains becomes entirely feasible on the current generation of high-performance computing architectures. The method was validated on a series of canonical flow problems in both the continuum and rarefied regime, showing high-order spatial accuracy and discrete conservation as well as sub-element resolution of shock waves without the need of any \textit{ad hoc} numerical shock capturing techniques. Furthermore, the results of the simulation of a three-dimensional compressible Taylor--Green vortex, which is, to the authors' knowledge, the first such case of a three-dimensional turbulent flow simulated via directly solving the Boltzmann equation, demonstrate the ability of the proposed approach as a method for direct numerical simulation of compressible turbulence. 

Extensions of this method will be performed on the analysis and validation of wall boundary conditions, the use of moment-based velocity discretizations, implementation of implicit time integration to address stiffness in the collision term, and applications to high-temperature and rarefied flows. We believe that the demonstrated ability of directly solving the Boltzmann equation opens a new path for the study of transition to turbulent and shock boundary layer interactions in high speed flow since the evolution of the probability density function encodes the dynamics of the flow in a form that is not accessible from the solution of the compressible Navier--Stokes equations. Furthermore, this approach allows for a unified framework for simulating flows in both the continuum and rarefied regime.